\documentclass[12pt]{article}

\usepackage{graphicx}
\usepackage{subfigure}
\usepackage{amssymb,amsmath,amsthm,epsfig}
\usepackage{latexsym, enumerate}
\usepackage{eepic}
\usepackage{epic}
\usepackage{graphicx}
\usepackage{color}
\usepackage{amsmath}
\usepackage{multirow}

\DeclareMathOperator*{\argmin}{arg\,min}

\usepackage{ifpdf}
\usepackage{booktabs}

\usepackage{dsfont}
\usepackage{multirow}
\usepackage{bm}
\usepackage{caption}
\usepackage{subfig}
\usepackage{subfloat}

\theoremstyle{plain}

\theoremstyle{definition}

\theoremstyle{remark}

\begin{document}
\title{\large\bf An improved implicit sampling for Bayesian inverse problems of multi-term time fractional multiscale diffusion models}
\author{
Xiaoyan Song\thanks
{College of Mathematics and Econometrics, Hunan University, Changsha 410082, China.
Email: xiaoyansong@hnu.edu.cn}
\and
Lijian Jiang\thanks
{School  of Mathematical Sciences, Tongji University, Shanghai 200092, China. Email: ljjiang@tongji.edu.cn.}
\and
Guang-Hui Zheng\thanks
{College of Mathematics and Econometrics, Hunan University, Changsha 410082, China.
Email: zhgh1980@163.com.}
}
\date{}
\maketitle
\begin{center}{\bf ABSTRACT}
\end{center}\smallskip
This paper presents an improved implicit sampling method for hierarchical Bayesian inverse problems.
A widely used approach  for sampling  posterior distribution is based on Markov chain Monte Carlo (MCMC).
However, the samples generated by MCMC are usually  strongly  correlated.  This  may lead to a small size of effective samples from a  long Markov chain and the resultant  posterior estimate
may be inaccurate. An implicit sampling method proposed in \cite{PNAS-Chorin} can generate independent samples and capture  some inherent  non-Gaussian features of the posterior based on the weights of samples.
In the implicit sampling method, the posterior samples are generated by constructing a map and distribute around the MAP point.
However, the weights of implicit sampling in previous works may cause excessive concentration of samples and lead to ensemble collapse. To overcome this issue,
we propose a new weight formulation and make resampling based on the new weights.
In practice, some parameters in prior density are often unknown and a hierarchical Bayesian inference  is necessary for posterior exploration.
To this end, the hierarchical Bayesian  formulation is  used to estimate the MAP point  and
integrated in  the implicit sampling framework. Compared to  conventional  implicit sampling, the proposed implicit sampling method can significantly improve the posterior estimator and
the  applicability for high dimensional inverse problems.
The improved implicit sampling method is applied to the Bayesian inverse problems of multi-term time fractional diffusion models in heterogeneous media.
To effectively capture the heterogeneity effect,  we present a mixed generalized  multiscale finite element method (mixed GMsFEM) to
solve the time fractional diffusion models in a coarse grid, which can substantially speed up the Bayesian inversion.  Using
the improved implicit sampling, we carry out a few numerical examples  to effectively identify  various unknown inputs  in these models.

\smallskip
{\bf keywords}: Bayesian inversion, improved implicit sampling, mixed GMsFEM, multi-term time fractional diffusion models

\section{Introduction}
The fractional diffusion equations have been widely used in modelling of the anomalous phenomenon, such as fractal media, chemistry, and material physics \cite{frac2, frac3, frac4, frac5}. The time fractional diffusion equations are obtained from the standard diffusion equations through replacing the integer-order time derivative with a fractional derivative, which can be derived from  random walk models or generalized master equations.
The diffusion process for inhomogeneous anisotropic media
modelling by time fractional diffusion equations is also very common in physics
\cite{Carcione-2013, Fomin-2011, Berkowitz-2002}.
In order to make the modelling more precise, fractional diffusion equations with multi-term time fractional derivatives have been investigated in recent years \cite{frac6, frac7}. If all parameters in the fractional models are given, there are many methods to solve the forward model \cite{forward-fem-1, forward-fem-2}.
However, in practice, there are many inputs unknown in the fractional model, such as the multi-fractional derivatives, the diffusion field and the reaction field. Thus, identification of the unknown inputs in fractional diffusion models lead to inverse problems. To the authors' best knowledge, there are only few works  focused on the inverse problems for multi-term time fractional diffusion model, such as Jiang et al. \cite{jiang2017} for determining the source term, Sun et al. \cite{sun2019} for identifying reaction coefficient, and Li et al. \cite{li2016} for recovering fractional orders and reaction coefficient simultaneously.

The estimation of unknown parameters from limited and noisy measurement data is a big challenge in science and engineering. Uncertainty in model inputs
leads to the uncertainty in model predictions, which can in turn have a great impact on experimental design and decision-making. In general, the inverse  problems are often ill-posed because of the limited indirect and noisy data.
One of the most popular approaches to estimate the unknown parameters is to minimize a cost function involving the misfit between the measurement data and the simulated response based on a regularization term. This is usually called Tikhonov regularization \cite{CR.Vogel-2002}. Iterative methods are  used to solve the minimization problems. However, these deterministic approaches usually only generate a point estimate of the unknown parameters, without quantifying the uncertainty in parameters. To treat this situation, the Bayesian
statistical approach \cite{mcmc, Hiera1, Jiang, kle-book} is often used to solve  the inverse problems.

The Bayesian inference provides a natural mechanism for learning from the measurement data by incorporating some prior information and has distinct features over classical deterministic regularization methods. One can obtain
the marginalizing, moments ,or confidence intervals through characterizing the posterior distribution. However, only in some special cases, the posterior may be in a closed form. Even when the prior and likelihood are Gaussian, the posterior may be non-Gaussian due to the nonlinearity of the forward operator. The most common choice for exploring the posterior is Markov chain Monte Carlo (MCMC) \cite{mcmc}, which generates a serial  of samples to evaluate the statistical information. However, MCMC  has some drawbacks.
For instance, MCMC usually requires millions of samples to estimate the statistical information. This means that a large number of  forward model simulations need to be computed. In addition,  the resultant samples may be correlated and  the  size of effective samples is small.  It is also difficult to diagnose the convergence of MCMC in practice.

For computationally expensive  forward models, sampling using MCMC may be infeasible.
Hence, it is necessary  to find  some other sampling methods  with a small computation  burden. The ensemble-based method is one of the approaches and give an  approximation of the posterior.
The typical ensemble methods include  linearization around the maximum a posterior distribution (LMAP) \cite{mcmc}, randomized maximum likelihood (RML) \cite{SpRML}, and ensemble Kalman filter (ENKF) \cite{Data-assimilation, YumingBa-2018}, etc..
These methods may be inaccurate because they  can be regarded as Gaussian approximations of the posterior.
In addition, there exist some other methods based on optimization  to improve the performance of sampling in the frame of Bayesian inference, such as  randomize-then-optimize (RTO) \cite{RTO} and optimal map \cite{optimal-maps}. RTO produces candidate samples from solving a randomly perturbed optimization problem, which acts as a Metropolis proposal. And \cite{RTO-L1} extends RTO to  some cases with non-Gaussian prior. The optimal map needs to construct a map that  push the prior measure to posterior measure. However, the construction of the map is often computationally expensive.

An alternative to the above approaches is importance sampling \cite{Data-assimilation, mcmc}, which is a sampling method without any Gaussian assumption.
 The idea is to draw samples from another easy-sampling importance function with a weight of each sample instead of
drawing  samples from the target distribution itself, which  is usually difficult to explore directly.
But the effective sample size \cite{Data-assimilation, mcmc} may  be very small in the conventional importance sampling  if the variance of the weights are large. This
 situation  is called filter degeneracy or ensemble collapse \cite{nonlinear-data-assimilation}.
It is critical  to carefully select the important function, which can in turn influence the weights.
Recently, a method called implicit sampling  was proposed in \cite{PNAS-Chorin} and studied for inverse problems \cite{IS_data-assimilation, IS_random, implicit_sampling, IS_seq}.
The implicit sampling method acts  as a special formulation of importance sampling to improve sample performance by  providing an importance function based on optimization.
In order to make samples concentrate on the region with
high probability, implicit sampling firstly solve a minimization problem to obtain the MAP point.
Then a serial  of samples can be generated through solving a nonlinear algebraic equation with a random right hand term.
The implicit sampling shares some similarities with LMAP. But  the implicit sampling can capture some non-Gaussian information of the posterior
based on the weight of each sample
compared with LMAP. Moreover, the implicit sampling produces independent samples, which greatly reduce the autocorrelation of samples and avoid the long burn-in process in  MCMC.
However, there still exist many shortcomings in the conventional  implicit sampling.
For example, the conventional  implicit sampling may still cause  sample degeneracy. The implicit sampling is stilled considered for Gaussian prior in most of the previous works \cite{PNAS-Chorin, IS_data-assimilation, IS_random, implicit_sampling, IS_seq}.

In this paper, we present an improved implicit sampling method and its application with some non-Gaussian priors.
When the dimension of unknown parameter is high and  the model  structure is complicated, the weights of the conventional implicit sampling may cause excessive concentration of samples and lead to ensemble collapse, which may give a poor estimation of the posterior. To overcome the issue, we propose a new weight formulation, which can properly balance the weights of samples.
 Meanwhile, this new formulation of weight can still retain  the original weight order.
In practice, the prior information often involves unknown parameters.  Hierarchical Bayesian model \cite{Hiera1} is necessary to treat this situation.
In this paper, we consider  the hierarchical Bayesian model when making the MAP estimate  in implicit sampling.
In addition, we use  a Laplace prior for the inputs with heterogeneous property depending on the knowledge of sparsity regularization \cite{SpRML, Li.L-2010}.

For the multi-term time fractional diffusion models in heterogeneous media, there may exist some multiple scales and high contrast feature in the diffusion field.
Direct simulating the multiscale models may be  very computationally  expensive. But  Bayesian inversion requires to simulate  the forward model many times.
Hence, it is necessary to use a suitable  model reduction  for multiscale models to significantly improve the computation efficiency.
It is well known that  Generalized Multiscale Finite Element Method (GMsFEM) \cite{GMs3, GMs1, GMs2} can efficiently solve the multiscale models in a coarse grid.
GMsFEM  aims at dividing the computation into offline and online steps.
One entails constructing the offline space in order to provide a good approximation of the solution with fewer basis functions.
The offline space is established depending on the spectral decomposition of the snapshot space, which can be obtained by solving some local problems.
In this paper, we utilize  mixed GMsFEM \cite{mixed-GMs1, mixed-fem1} to solve  the multi-term time  fractional multiscale diffusion models.
The mixed GMsFEM can give accurate  flux information.  The measurement data are
the flux information collected at the boundary of the spatial domain.

The paper is organized as follows. In Section 2, we introduce the implicit sampling in Bayesian framework, which includes the approaches to compute  the MAP point for different prior information.
In Section 3, we introduce the multi-term time fractional diffusion model and then construct the reduced model using  mixed GMsFEM.
Some numerical results are presented in Section 4 to show the performance of the proposed method. The paper ends with some conclusions and comments in Section 5.

\section{Bayesian inference based on implicit sampling}
We begin with describing the Bayesian framework and denoting some notions for the paper. We consider the following general model
\begin{eqnarray}
\label{Model-general}
\mathcal{B}(\theta,u)=0 \ \ \text{in} \ \ \Omega,
\end{eqnarray}
which can represent a PDE modeling some physical problems.
If the model inputs $\theta$ are available, the model output $u$ can be obtained by numerical simulation. Let  $\theta\in \mathds{R}^{m}$.
In Bayesian inversion, the unknown parameters and measurement data are usually regarded as random variables.
Let $(\mathcal{S},\mathcal{F},\mathds{P})$ be a measure space, where $\mathcal{S}$ is the sampling space, $\mathcal{F}$ is the $\sigma\text{-field}$ and $\mathds{P}$ is the probability measure.
Then  $p(\theta)=d\mu_{0}/d\theta$ is the density of $\theta$ with respect to Lebesgue measure, where $\mu_{0}$ is a prior measure.
According to Bayes's rule, the posterior probability can be obtained by
\begin{eqnarray}
\label{Baye's}
p(\theta|d)=\frac{p(d|\theta)p(\theta)}{\beta},
\end{eqnarray}
where the function $p(d|\theta)$ is the  likelihood function and
$\beta$ is the evidence or marginal likelihood $\beta=\int p(d|\theta)p(\theta)d\theta$.
Assume that $\theta^{*}\in \mathds{R}^{m}$ is the reference parameter.
We assume that  the measurement data is  generated  by
\begin{eqnarray}
\label{data}
d=H(\theta^{*})+\epsilon,
\end{eqnarray}
where the function $H:\mathds{R}^{m}\rightarrow \mathds{R}^{n}$ denotes the parameter-to-observation  (or forward) map according to model (\ref{Model-general}) and $\epsilon$ is a random variable independent of $\theta$. Generally, we assume $\epsilon \sim N(0,\sigma^2 \text{I})$, where $\text{I}\in \mathds{R}^{n\times n}$ is the identity matrix.
In Bayesian inverse problems, one intends to explore the posterior, which  relies  on the sparse measurement data $d$ and some prior information $p(\theta)$.

\subsection{Hierarchical Bayesian formulation}\label{sec-2.1}
Under the above assumptions, it can be concluded that the likelihood function has the following form

\[
p(d|\theta)\propto (\sigma^2)^{-n/2}\exp(-\frac{\|H(\theta)-d\|^{2}_{2}}{2\sigma^2}).
\]
The prior $p(\theta)$ depends  on the knowledge of the unknown inputs and  plays a critical role on exploring the posterior. A versatile choice for $p(\theta)$ is Markov random field (MRF), i.e.
\[
p(\theta)\propto \lambda^{\frac{m}{2}}\exp(-\frac{\lambda}{2}\|\theta\|^{2}_{2}).
\]
If $\lambda$ is known, the posterior density $p(\theta|d)$ can be obtained as

\[
p(\theta|d)\propto\lambda^{\frac{m}{2}}
\exp(-\frac{\|H(\theta)-d\|^{2}_{2}}{2\sigma^2})
\exp(-\frac{\lambda}{2}\|\theta\|^{2}_{2}).
\]
The posterior density $p(\theta|d)$ provides the complete distribution of $\theta$ relying on the observations $d$.
We can compute the MAP point, ${\theta}_{MAP}:=\arg \max_{\theta}p(\theta|d)$, which is equivalent to the following minimization problem
\[
\theta_{MAP}=\arg \min\big\{\|H(\theta)-d\|^{2}_{2}+\mu \|\theta\|^{2}_{2}\big\},
\]
where $\mu=\lambda\sigma^2$ plays the role of regularization parameter in  inverse problems.
 However, as we know, it is not an easy task to ascertain the regularization parameter, especially in nonlinear inverse problems. Fortunately, the hierarchical Bayesian model \cite{mcmc,Hiera1} can provide an effective and flexible way to overcome the difficulty.
 In hierarchical Bayesian framework, $\lambda$ can be regarded as a hyper-parameter, i.e., we can choose a hyper-prior $p(\lambda)$ for it. Then the posterior density can be written as
\[
p(\theta,\lambda)\propto p(d|\theta)p(\theta)p(\lambda).
\]
A common method to choose hyper-prior $p(\lambda)$ is using conjugate priors. If we use Gamma distribution as the hyper-prior for $\lambda$, then the posterior density becomes
\begin{eqnarray}
\label{hyper_post}
p(\theta,\lambda|d)\propto
\lambda^{\frac{m}{2}}
\exp(-\frac{\|H(\theta)-d\|^{2}_{2}}{2\sigma^2})
\exp(-\frac{\lambda}{2}\|\theta\|^{2}_{2})
\lambda^{a-1}\exp(-b\lambda),
\end{eqnarray}
where $\lambda\thicksim Gamma(a,b)$ \cite{mcmc}.
Then the MAP estimate for (\ref{hyper_post}) can be computed by minimizing the following functional
\begin{eqnarray}
\label{func_map}
\mathcal{M}(\theta,\lambda)=\frac{\|H(\theta)-d\|^2_{2}}{2\sigma^2}+
\frac{\lambda}{2}\|\theta\|^{2}_{2}+b\lambda-(\frac{m}{2}+a-1)\ln\lambda.
\end{eqnarray}
However, when solving the forward model, it  usually depends  on the spatially discretized grid system. Thus if the random vector $\theta$ is grid-based, its dimension may be very high. This
 will result in the instability and inefficiency of the numerical algorithms. Therefore, the reduction of the dimension of parameter space is necessary. Let  the parameter $\theta$ can be represented  by
\begin{eqnarray}
\label{Bases}
\log \theta=\Phi \textbf{v}+\kappa,
\end{eqnarray}
where $\textbf{v}=[v_{1},v_{2},\cdots, v_{l}]$ are the coefficients under the base matrix  $\Phi\in  \mathbb{R}^{m\times l}$ and $\kappa$ is a fixed constant vector. Thus the parameter-to-observation  map $H(\theta)$ can be expressed by $H(\textbf{v})$.
For convenience, we use $\textbf{v}^{*}$ to indicate the reference parameter. We can use the  augmented-Tikhonov method \cite{Hiera1}  to get the MAP estimate.
Table \ref{a-Tikh_algorithm} describes the augmented-Tikhonov method when the prior for $\textbf{v}$ is assumed to be Gaussian.

\begin{table}[htbp]
\centering
 \caption{Augmented-Tikhonov algorithm}\label{a-Tikh_algorithm}
 \begin{tabular}{l}
  \toprule
 \textbf{Input}: The noise level $\sigma$,
  the maximum iterations $R$ and the precision $eps$ .\\
  \textbf{Output}: $ \textbf{v}_{k}=(v_{1}^{k},...,v_{l}^{k})^{T}.$\\
  \midrule
  1.  $\bf{Initialize}$ Set k=0 and give the initial value $\textbf{v}_{0}$ and $\mu_{0}$;\\
  2.  Settle the forward problem and obtain the additional data  $d=H(\textbf{v}^{*})+\epsilon$;
  \\
  3.  $\bf{While}$ $k<R$\\
  4. Compute  $
  h_{k}=(\bar{H}^{T}\bar{H}+\mu_{k-1}I)^{-1}\bar{H}^{T}(d-H(\textbf{v}_{k-1}))
  $,
  where $\bar{H}$ is the sensitivity matrix;\\
  5. Update $\textbf{v}_{k}$ by
   $\textbf{v}_{k}=\textbf{v}_{k-1}+h_{k};$
 \\
  6.  Update the hyper-parameter $\lambda_{k}$ by
  $\lambda_{k}=\frac{\frac{l}{2}+a-1}{\frac{1}{2}\|\textbf{v}_{k}\|^{2}+b}$ and set $\mu_{k}=\lambda_{k}\sigma^2$;
  \\
  7.  $k\leftarrow k+1$;\\
  8.  $\bf{Terminate}$ if $\|h_{k}\|<eps$;\\
  9.  $\bf{end}$\\
  \bottomrule
 \end{tabular}
 \end{table}

\subsection{MAP estimate with Laplace prior}\label{sec-2.2}
In this subsection,  we consider a non-Gaussian prior case, Laplace prior, which is often used in describing  realistic geologic system \cite{SpRML, Li.L-2010}.
Assume the grid-based parameter $\log \theta$ has an approximation formulation as (\ref{Bases}). 
A common measure of sparsity of the vector $\textbf{v}$ involves the number of nonzero entries in $\textbf{v}$, i.e.,  $\|\textbf{v}\|=\sharp \{i, v_{i}\neq 0\}$. The sparse reconstruction aims at recovering the sparse vector $\textbf{v}$ based on the measurement data $d$. The original sparse inversion problem can be represented as
\[
\min \|\textbf{v}\|_{0} \ \ \  s.t.\ \  H(\textbf{v})=d.
\]
The exact solution to the above $l_{0}$-$norm$ minimization problem is usually NP-hard. As an alternative,
one can find a sparse solution under some mild condition by solving the following problem
\begin{eqnarray}\label{P_{1}}
\min \|\textbf{v}\|_{1} \ \ \  s.t.\ \  H(\textbf{v})=d.
\end{eqnarray}
The constrained optimization problem  (\ref{P_{1}}) can be represented by the following  unconstrained optimization problem
\begin{eqnarray}\label{l1-min}
\min \mathcal{J}(\textbf{v})=\|H(\textbf{v})-d\|_{2}^{2}+\mu\|\textbf{v}\|_{1},
\end{eqnarray}
where $\mu$ acts as the regularization parameter.

For  Bayesian inverse problems, if we select the Laplace distribution as a prior for $\textbf{v}$, i.e.,  $p(\textbf{v})=\frac{\lambda}{2} \exp(-\lambda \|\textbf{v}\|_{1})$, then the posterior density is
\begin{eqnarray}
p(\textbf{v}|d)\propto
\exp(-\mathcal{M}(\textbf{v})), \quad \text{where} \quad \mathcal{M}(\textbf{v})=\frac{\|H(\textbf{v})-d\|_{2}^{2}}{2\sigma^2}+\lambda\|\textbf{v}\|_{1}.
\end{eqnarray}
Let $\mu=2\lambda\sigma^2$, then we can find the MAP estimate of $p(\textbf{v}|d)$ by solving the
minimization problem (\ref{l1-min}).
However, since the derivative of the function $\mathcal{J}(\textbf{v})$ is not defined at the origin,
we use the iteratively reweighted approach \cite{Li.L-2010, SpRML} to solve (\ref{l1-min}), which is widely used in compressive sensing.
As mentioned above, the cost function in the iteratively reweighted approach
is
\begin{eqnarray}\label{itera-reweighted}
\min
\mathcal{J}(\textbf{v})=\|H(\textbf{v})-d\|_{2}^{2}+\mu\|\textbf{v}\|_{W}^{2},
\end{eqnarray}
where $W$ is a diagonal weighting matrix generated by the vector $w=[(v_{1}^{2}+\varepsilon)^{-\frac{1}{2}},\cdots, (v_{l}^{2}+\varepsilon)^{-\frac{1}{2}}]$ with a small constant $\varepsilon > 0$, and  $\|\textbf{v}\|^{2}_{W}=\textbf{v}^{T}W\textbf{v}$.
Note that when $\varepsilon=0$, the term $\|\textbf{v}\|^{2}_{W}$ actually equals to $\|\textbf{v}\|_{1}$.
We have  obtained  a similar  formalization of the minimization function for the Gaussian prior and Laplace prior.
Based on (\ref{itera-reweighted}), the iteratively reweighted approach is listed in Table \ref{itera-reweight algorithm}.

\begin{table}[htbp]
\centering
 \caption{The iteratively reweighted approach}\label{itera-reweight algorithm}
 \begin{tabular}{l}
  \toprule
 \textbf{Input}: The noise level $\sigma$, the regularization parameter $\mu$,
  the maximum iterations $R$,\\ and the positive constant $\varepsilon$.\\
  \textbf{Output}: $ \textbf{v}_{k}=(v_{1}^{k},...,v_{l}^{k})^{T}.$\\
  \midrule
  1.  $\bf{Initialize}$ Set k=0 and give the initial value $\textbf{v}_{0}$ and $\mu_{0}$;\\
  2.  Settle the forward problem and obtain the additional data  $d=H(\textbf{v}^{*})+\epsilon$;
  \\
  3.  $\bf{While}$ $k<R$\\
  4. Update the matrix $W$;\\
  5. Compute  $
  \textbf{v}_{k}=(\bar{H}^{T}\bar{H}+\mu W)^{-1}\bar{H}^{T}(d-H(\textbf{v}_{k-1})
  +\bar{H}\textbf{v}_{k-1})
  $,\\
  where $\bar{H}$ is the sensitivity matrix;\\
  6.  $\bf{end}$\\
  \bottomrule
 \end{tabular}
 \end{table}
The above approaches focus on  the point estimate (MAP)  of unknown inputs from the perspective of variation formulation.
However, it is still not enough to characterize the uncertainty of the unknown parameters and the uncertainty propagation of model response.
As we know, the most popular method to explore the posterior by samples is Markov chain Monte Carlo (MCMC) \cite{mcmc}. Based on large number of samples generated by MCMC, some statistical information can be estimated to approximate the target  distribution. However,  MCMC  usually needs a lot of  samples to estimate the parameter and it is also difficult to diagnose if the chain has converged to the posterior distribution.
To overcome these difficulties, we introduce an improved implicit sampling method to efficiently explore the posterior samples.


\subsection{Improved implicit sampling}

In this subsection, we present implicit sampling, which  shares some idea with importance sampling.  To this end, we  first describe the basic idea of importance sampling \cite{mcmc, nonlinear-data-assimilation}. Importance sampling methods can be regarded as multi-step sampling generation techniques.
Suppose that we want to sample from a  probability density function (pdf) $\rho(\theta)$. However, it is not always an easy task to draw from $\rho(\theta)$ directly. In importance sampling, one first needs to select another easy-sampling pdf $q(\theta)$ and then provides some sort of correction mechanism to draw the samples from $q(\theta)$ to
approximate the interest distribution $\rho(\theta)$. As an example, we  consider the integration of the form
\[
I=\int f(\theta)\rho(\theta)d\theta.
\]
If the samples $\{\theta_{i}\}_{i=1}^{N_{s}}$ can be obtained from $\rho(\theta)$ directly, then we can use the simple Monte Carlo (MC) estimator \cite{mcmc} to calculate the integral:
\[
I\approx \frac{1}{N_{s}}\sum_{i=1}^{N_{s}}f(\theta_{i}) \ \ \text{as} \ \ N_{s} \rightarrow \infty.
\]
If $\rho(\theta)$ is not easy for sampling, the integral can be rewritten in the following form with the auxiliary density $q(\theta)$:
\[
I=\int f(\theta)\rho(\theta)dx=\int f(\theta)\frac{\rho(\theta)}{q(\theta)}q(\theta)d\theta.
\]
Draw samples $\theta_{i}\sim q(\theta)$, and then
\[
I\approx \frac{1}{N_{s}}  \sum_{i=1}^{N_{s}}f(\theta_{i})\omega_{i} \ \ \text{as} \ \ N_{s} \rightarrow \infty,
\]
where the weight $\omega_{i}=\frac{\rho(\theta_{i})}{q(\theta_{i})}$.
The auxiliary density must be chosen carefully.
In particular, when the interest distribution is the posterior $p(\theta|d)$ and the auxiliary density is chosen as the prior $p(\theta)$, then according to Baye's rule,
\begin{equation}
\label{import-sam}
E(f(\theta))=\int f(\theta)p(\theta|d)d\theta \approx \sum_{i=1}^{N_{s}}\omega_{i}f(\theta_{i}),
\end{equation}
where $\theta_{i}\sim p(\theta)$ and $\omega_{i}=\frac{p(d|\theta_{i})}{\sum_{i=1}^{N}p(d|\theta_{i})}$.
However, when the prior and likelihood have disjoint support or the prior is too broad, the  estimate in  (\ref{import-sam}) is not accurate.
Implicit sampling can effectively improve the conventional importance sampling in this situation.
Based on \cite{implicit_sampling}, we give a brief introduction to the implicit sampling, which concentrate  on the region with  the high  posterior probability.

Let
\[
F(\theta)=-\log(p(d|\theta)p(\theta)).
\]
First, we need to find the MAP estimate by minimizing the negative logarithm of the posterior.
Suppose that
\[
\theta_{MAP}=\argmin_{\theta}F(\theta), \quad \phi_{F}=\min F(\theta).
\]
 Then we intend to select an auxiliary function to efficiently draw samples from high probability region, i.e., near the MAP point. To this end, we can use  standard Gaussian distribution
\[
\bar{q}(\xi)\varpropto \exp(-\frac{1}{2}\xi^{T}\xi)
\]
to play the role of the easy-sampling density. Denote $G(\xi)=\frac{1}{2}\xi^{T}\xi$ and $\phi_{G}=\min G(\xi)$, then we solve the following nonlinear algebraic equation:
\begin{equation}
\label{Is-equation}
F(\theta)-\phi_{F}=G(\xi)-\phi_{G}.
\end{equation}
It can be concluded that the left-hand of (\ref{Is-equation}) tends to be small since samples drawn from $\bar{q}(\xi)$ can make $G(\xi)$ close to $\phi_{G}$. By solving the above algebraic equation for a series of samples $\xi$, we can obtain the samples located around the MAP point.
Based on this idea, we define a map
$\mathcal{A}: \xi \rightarrow \theta$ to find the solution.
The pdf $q(\theta)$ of $\theta$ can be obtained by
\[
q(\theta)=\bar{q}(\mathcal{A}^{-1}(\theta))|\det(\frac{\partial \mathcal{A}^{-1}(\theta)}{\partial \theta})|,
\]
where $\frac{\partial \mathcal{A}^{-1}(\theta)}{\partial \theta}$ is the Jacobian matrix for $\mathcal{A}^{-1}(\theta)$ with respect to $\theta$.
If we select $q(\theta)$ as the auxiliary function in importance sampling, then the weight is
\begin{equation}
\label{ratio}
\omega\propto \frac{p(\theta|d)}{q(\theta)}\propto\exp(G(\mathcal{A}^{-1}(\theta))-F(\theta))
|\det(\frac{\partial \mathcal{A}^{-1}(\theta)}{\partial \theta})|^{-1},
\end{equation}
 The second-order Taylor  expansion of  $F(\theta)$ is given by
\[
\hat{F}(\theta)=\phi_{F}+\frac{1}{2}(\theta-\theta_{MAP})^{T}\mathcal{H}(\theta-\theta_{MAP}),
\]
where $\mathcal{H}$ is the Hessian matrix of $F(\theta)$ at the minimum $\theta_{MAP}$.
Then the above method is usually called  Gaussian  approximation \cite{mcmc} for the posterior $p(\theta|d)$,
 i.e., $\theta\sim N(\theta_{MAP},\mathcal{H}^{-1}).$
Denote $\bar{H}$ to represent the Fr$\acute{e}$chet derivative of $H(\theta)$ evaluated  at $\theta=\theta_{MAP}$, i.e.,  linearization for $H(\theta)$.
In particular, if the prior is Gaussian,  we can obtain
\[
\mathcal{H}^{-1}\approx C-C\bar{H}^{T}(\bar{H}C\bar{H}^{T}+\Gamma)^{-1}\bar{H}C,
\]
where $C$ and $\Gamma$ are the prior covariance matrix and likelihood covariance matrix, respectively.
If the prior is Laplace distribution, $\mathcal{H}$ can be represented as
\[
\mathcal{H}\approx \bar{H}^{T}\Gamma^{-1}\bar{H}+2\lambda W,
\]
where $W= W(\theta_{MAP})$.
By substituting  $\hat{F}(\theta)$ with $F(\theta)$ in (\ref{Is-equation}), we  solve the following equation
\begin{equation}
\label{euation-linear}
\hat{F}(\theta)-\phi_{F}=\frac{1}{2}\xi^{T}\xi.
\end{equation}
Then we can define  the map $\mathcal{A}$  as
\begin{equation}
\label{linear map}
\mathcal{A}(\xi)=\theta_{MAP}+L^{T}\xi,
\end{equation}
where $\xi\sim N(0,\text{I})$ and L is the Cholesky factorization of $\mathcal{H}^{-1}$, i.e.,  $\mathcal{H}^{-1}=L^{T}L$.
It is easy to see that (\ref{linear map}) is one of the solutions of (\ref{euation-linear}).
Then the weights of samples in (\ref{ratio}) are
\begin{equation}
\label{weight-propotion}
\omega\propto \exp (\hat{F}(\theta)-F(\theta)).
\end{equation}
Once the randomly perturbed samples in (\ref{linear map}) are produced, we can implement the resampling technology \cite{nonlinear-data-assimilation}  based on the weights in (\ref{weight-propotion}) to characterize the posterior. The basic idea of resampling is that samples with low weights are abandoned, while multiple copies of samples with high weights are kept.
It can be interpreted that the higher the weight of a sample, the more copies
of the sample are produced. The resampling can be performed in many methods and in this paper we utilize  stochastic universal sampling \cite{nonlinear-data-assimilation}.

However, it will come out an extreme situation, where  the weight of some  sample is concentrated near 1 and the weights  of  the remaining samples are almost 0. In this case, most samples generated by implicit sampling are the copies of
the sample with high weight, which can be regarded as the point estimate and will give  a poor distribution for  the posterior. This is usually called filter degeneracy or ensemble collapse.
In order to avoid this issue,  we relax the weights in  (\ref{weight-propotion}) by
\begin{equation}
\label{weight-relax1}
\omega\propto \exp (\frac{\hat{F}(\theta)-F(\theta)}{\vartheta}),
\end{equation}
where $\vartheta$ is a positive constant greater than 1.
\begin{figure}
\centering
\includegraphics[width=0.45\textwidth, height=0.4\textwidth]{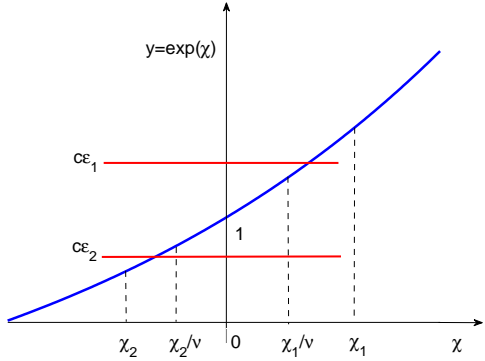}
\caption{Illustration of  the improved  weights in (\ref{weight-relax1}); The positive constant $\vartheta$ can squeeze  the weights $\frac{1}{c}\exp(\frac{\hat{F}(\theta)-F(\theta)}{\vartheta})$ into $[\epsilon_{1}, \epsilon_{2}]$, where $c$ is the normalization constant.}\label{figure-weight}
\end{figure}
By (\ref{weight-relax1}), we find that the new weights   not only keep the original order in (\ref{weight-propotion}) but also alleviate the degeneracy phenomenon, where one of weights is  too close to 1.
In other words, a proper  selection of $\vartheta$ will decrease  the large weight and increase the small weight  simultaneously, which can squeeze  the weights into  the interval $[\epsilon_{1}, \epsilon_{2}]$ (see Figure \ref{figure-weight}).
In particular,  the weight of each sample will be close to $\frac{1}{N_{s}}$ when $\vartheta$ tends to infinity, where $N_{s}$ is  the number of  samples from the  implicit sampling.
A simple method for dynamically selecting scale parameters $\vartheta$ is  based on effective sample size (ESS) \cite{mcmc}, which can be regarded as the size of  random samples with the same variance¡£
The algorithm  is given in Table \ref{Choice-rescale-constant}.
As a special case of importance sampling, we can compute the ESS of improved implicit sampling according to the following formula \cite{Ess_weight2, ESS_weight}

\begin{equation}
\label{ess_w}
\text{ESS}=\frac{1}{\sum_{i=1}^{N}\omega_{i}^2},
\end{equation}
where $\omega_{i}$ is the weight in improved implicit sampling.
We first give a preset value $\mathcal{S} $ of the target ESS  and the maximum iterations $N_{\vartheta}$. If the ESS is too small, we can increase $\vartheta$ until it reaches the threshold value $\mathcal{S}$.
\begin{table}[htbp]
\centering
 \caption{A strategy for selecting scale parameter $\vartheta$}\label{Choice-rescale-constant}
 \begin{tabular}{l}
  \toprule
 \textbf{Input}: the preset value of $\mathcal{S}$ and the maximum iteration $N_{\vartheta}$;\\
  \textbf{Output}: $\vartheta_{k}$\\
  \midrule
  1.  $\bf{Initialize}$ Set $k=0$ and give the initial value $\vartheta_{0}=1$;\\
  2.  $\bf{While}$ $k<N_{\vartheta}$\\
  3. Compute ESS by $\rm ESS=\frac{1}{\sum_{i=1}^{N}\omega_{i}^2}$;\\
  4. If $\text{ESS}<\mathcal{S}$\\
  \ \ Update $\vartheta_{k}=\vartheta_{k}+1$ and set $k=k+1$;\\
  5.  $\bf{end}$\\
  \bottomrule
 \end{tabular}
 \end{table}
We can find  that the samples generated by improved implicit sampling are not only independent of  each other but can concentrate on the region with high probability. In improved implicit sampling, one can avoid the case where the support of posterior density is too wide and catch some non-Gaussian feature of the posterior simultaneously.
Improved implicit sampling shares the merit of importance sampling and implies a high probability region of the posterior by solving a nonlinear algebraic equation.
Based on the Gaussian prior and Laplace prior described  in section \ref{sec-2.1} and \ref{sec-2.2}, the algorithm of the improved implicit sampling is listed in Table \ref{IS-algorithm}, where the scale parameter $\nu$ is selected by Table \ref{Choice-rescale-constant}.

\begin{table}[htbp]
\centering
 \caption{Improved implicit sampling}\label{IS-algorithm}
 \begin{tabular}{l}
  \toprule

  1. Compute the MAP point $\theta_{MAP}$ by augmented-Tikhonov algorithm or\\
  the iteratively reweighted approach;\\
  2. Compute the Cholesky factorization $L$ of matrix $\mathcal{H}^{-1}$;\\
  3. $\bf{For}$ $1 \leqslant i\leqslant N_{s}$\\
  4. Generate samples by $\theta_{i}=\theta_{MAP}+L^{T}\xi_{i}$, where $\xi_{i}\thicksim N(0,\text{I})$;\\
  5. Compute the weights of samples by
  $\omega_{i}=\exp (\frac{\hat{F}(\theta_{i})-F(\theta_{i})}{\vartheta})$, then normalization;\\
  6. Resampling based on weights;\\
  7.  $\bf{end}$\\
  \bottomrule
 \end{tabular}
 \end{table}

%

\section{Reduced model based on mixed GMsFEM}
We will apply the improved implicit sampling to the inverse problems of  the multi-term time fractional diffusion equation, which is described as
\begin{eqnarray}
\label{model-fpde}
\begin{cases}
(\gamma_{1}\cdot{}_0 D_t^{\alpha_{1}}+\gamma_{2}\cdot{}_0 D_t^{\alpha_{2}})u(x,t)-div(k(x)\nabla u(x,t))+q(x)u(x,t)=f(x,t)  \ \ \text{in} \ \  \Omega \times (0,T], \\
u(x,0)=0  \ \ \text{in} \ \ \Omega, \\
u(x,t)=g(x,t) \ \ \text{on} \ \  \partial \Omega \times (0,T],
\end{cases}
\end{eqnarray}
where $\gamma_{1}$ and $\gamma_{2}$ are some positive constants,  and $\alpha_{1},\alpha_{2} \in (0,1)$ are the fractional orders of the derivative in time. Here ${}_0 D_t^\alpha u $ refers to the Caputo derivative \cite{frac4} with respect to $t$, i.e.,
\begin{equation}\label{caputa-def}
{}_0 D_t^\alpha u =\frac {1}{\Gamma (1-\alpha)}\int_0^t \frac {\partial u(x,s)}{\partial s}\frac {ds}{(t-s)^\alpha},
\end{equation}
where $\Gamma$ is the Gamma function.
 In practice, the diffusion field $k(x)$ is  heterogeneous and spatially varies in different scales.

\begin{figure}
\centering
\includegraphics[width=0.55\textwidth, height=0.45\textwidth]{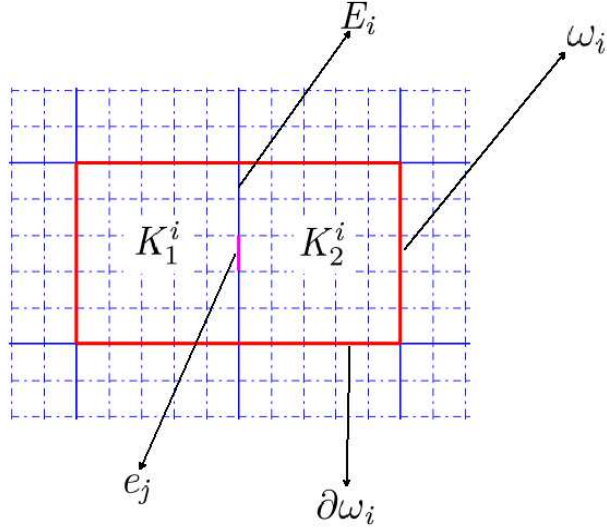}
\caption{An illustration of a neighborhood $\omega_{i}=K_{1}^{i}\cup K_{2}^{i}$ based on a coarse grid edge $E_{i}$.}
\label{mixedGMS-grid}
\end{figure}

In this section, we present a mixed GMsFEM for solving multi-term time fractional multiscale  diffusion equation. Following the mixed GMsFEM framework in \cite{mixed-GMs1},
 the main idea of  mixed GMsFEM is to divide the computation into offline stage  and online stage.  It  produces a coarsen computational model.
  For mixed GMsFEM, we first need to construct the snapshot space.  Let $E_{i}$ and $e_{j}$ be the coarse grid edge and the find grid edge,  respectively (see Figure \ref{mixedGMS-grid}).
  We define
\begin{eqnarray}
\label{deta-function}
\delta_{j}^{i}=
\begin{cases}
\ 1&\ \  \text{on $e_{j}$},\\
\ 0&\ \ \text{on other fine grid edges on $E_{i}$}.
\end{cases}
\end{eqnarray}
For the local snapshot  $\Psi^{i,snap}_{j}:=v_{j}^{i}$, we need to solve the following problem on the coarse neighborhood $\omega_{i}$ associated with the edge $E_{i}$,
\begin{eqnarray}
\label{basis-function}
\begin{cases}
k^{-1}v_{j}^{(i)}+\triangledown p_{j}^{i}=0\ \ \ in \ \ \omega_{i},\\
div(v_{j}^{(i)})=\alpha_{j}^{(i)}\ \ \ in \ \ \omega_{i},\\
v_{j}^{(i)}\cdot n_{i}=0 \ \ \ on \ \ \  \partial\omega_{i},
\end{cases}
\end{eqnarray}
where $n_{i}$ denotes the outward unit-normal vector on $\partial\omega_{i}$
and the constant $\alpha_{j}^{(i)}$ satisfies the compatibility condition $\int_{K_{l}}\alpha_{j}^{(i)}=\int_{E_{i}}\delta_{j}^{i}$.
We note that the coarse edge $E_{i}$ can be represented by a collection of the fine grid edges $e_{j}$, i.e.,  $E_{i}=\sum_{j=1}^{J_{i}}e_{j}$, where $J_{i}$ is the total number of the fine grid edges on $E_{i}$, see
Figure \ref{mixedGMS-grid}.
Then we define  the snapshot space $V_{snap}$ by
\begin{eqnarray*}
V_{snap}={\rm span}\{\Psi^{i,snap}_{j}: \ \ \ 1\leq j \leq J_{i},\ \ \ 1\leq i \leq N_{e}\}.
\end{eqnarray*}
Let $M_{snap}=\sum_{i=1}^{N_{e}}J_{i}$. Then the snapshot space can be represented as
\begin{eqnarray*}
V_{snap}={\rm span}\{\Psi^{snap}_{j}: \ \ \ 1\leq j \leq M_{snap}\}.
\end{eqnarray*}
Therefore we can define a matrix as
\begin{eqnarray*}
R_{snap}=[\psi_{1}^{snap},\cdots,\psi_{M_{snap}}^{snap}].
\end{eqnarray*}
To capture the important information in the snapshot space and reduce the dimension of the snapshot space, we consider the local spectral problem: find the eigen-pair $\{\lambda\}$ such that
\begin{eqnarray}
\label{spectral problem}
a_{i}(v,w)=\lambda s_{i}(v,w)   \ \ \ \forall w \in V_{snap}^{i},
\end{eqnarray}
where $a_{i}(v,w)$ and $s_{i}(v,w)$ are defined as following
\begin{eqnarray*}
\begin{cases}
a_{i}(v,w)=\int_{E_{i}}k^{-1}(v\cdot m_{i})(w\cdot m_{i})\\
s_{i}(v,w)=\int_{\omega_{i}}k^{-1}v\cdot w+\int_{\omega_{i}}div(v)div(w).
\end{cases}
\end{eqnarray*}
Then we can obtain a matrix representation as
\begin{eqnarray*}
A_{snap}^{(i)}Z_{k}^{i}=\lambda_{k}^{i}S_{snap}^{(i)}Z_{k}^{i},
\end{eqnarray*}
where
\begin{eqnarray*}
A_{snap}^{(i)}=R_{snap}^{T}A_{fine}^{i}R_{snap},\ \
S_{snap}^{(i)}=R_{snap}^{T}S_{fine}^{i}R_{snap}.
\end{eqnarray*}
Taking the first $L_{b}$ smallest eigenvalues and the corresponding  eigenfunctions, the global offline space can be constructed as
\begin{eqnarray*}
V_{H}=\span\{\Psi_{k}^{i}: \Psi_{k}^{i}=\sum_{j=1}^{J_{i}}Z_{kj}^{i}\Psi_{j}^{i,snap}, 1\leq k \leq L_{b}, 1\leq i \leq N_{e}\}.
\end{eqnarray*}
This space will be used to approximate the velocity in  mixed GMsFEM.
If we represent each $\Psi_{k}^{i}$ by the basis functions supported on fine grid and use a single index notation, then we can obtain the following matrix form
\begin{eqnarray*}
R_{\text{off}}=[\psi_{1},\psi_{2},\cdots,\psi_{M_{t}}],
\end{eqnarray*}
where $M_{t}=\sum _{i=1}^{Ne}L_{b}$ denotes the total number of multiscale offline basis functions. Once $R_{\text{off}}$ is computed,  it can  be repeatedly used for online simulation.

Let $u^{n}=u(.,t_{n})$.   In this paper, we use the standard  $L_{1}$ method \cite{caputo-deri-2} to approximate the Caputo derivative by
\begin{equation}\label{appro-caputo}
  {}_0 D_t^\alpha u^{n}=\frac{1}{(\Delta t)^{\alpha}\Gamma(2-\alpha)}\sum_{k=1}^{n}(u^{k}-u^{k-1})\times[(n+1-k)^{1-\alpha}-(n-k)^{1-\alpha}],
\end{equation}
where $0=t_{0}<t_{1}<\cdots<t_{M}=T, \ t_{n}=n\triangle t$.
Let $\eta=-k(x)\nabla u$, and $Q_{H}$ and $Q_{h}$ be the space of piecewise constant functions  with respect to the coarse grid and fine grid, respectively.  Then we can obtain the weak formulation for (\ref{model-fpde}),
\begin{equation*}
\begin{cases}
\int_{\Omega}k^{-1}(x)\eta \cdot v-\int_{\Omega}div(v)u=-\int_{\partial \Omega}g(x,t)v\cdot n, \ \ \forall v \in V_{H}\\
\int_{\Omega}((\gamma_{1}\cdot{}_0 D_t^{\alpha_{1}}+\gamma_{2}\cdot{}_0 D_t^{\alpha_{2}}))u\cdot p+\int_{\Omega}div(\eta)\cdot p+\int_{\Omega}q(x)u\cdot p=\int_{\Omega}f(x,t)p, \ \ \forall p \in Q_{H}.
\end{cases}
\end{equation*}
Suppose that $\eta^{n}=\sum _{i=1}^{M_{t}}\sigma_{i}^{n}\psi_{i}$ and $u^{n}=\sum_{k=1}^{J}\beta_{k}^{n}\phi_{k}$, where $\{\phi_{k}\}_{k=1}^{J}\subseteq Q_{H}$.

Let $s_{i}=(\triangle t)^{\alpha_{i}}\Gamma(2-\alpha_{i})\ (i=1,2)$, $s=\frac{s_{1}s_{2}}{\gamma_{1}s_{2}+\gamma_{2}s_{2}}$ and
\[
\begin{aligned}
&b_{n}=\sum_{i=1}^{2}[n^{1-\alpha_{i}}-(n-1)^{1-\alpha_{i}}]
\cdot\frac{\gamma_{i}s_{i+1}}{\gamma_{1}s_{2}+\gamma_{2}s_{1}}, \ \ 1\leqslant n \leqslant M,\\
&c_{k}=\sum_{i=1}^{2}[2k^{1-\alpha_{i}}-(k+1)^{1-\alpha_{i}}-
(k-1)^{1-\alpha_{i}}]\cdot\frac{\gamma_{i}s_{i+1}}{\gamma_{1}s_{2}+
\gamma_{2}s_{1}},\ \ 1\leqslant k \leqslant M-1,
\end{aligned}
\]
where $s_{3}=s_{1}$.
Then the above weak formulation can be transformed into the following matrix form,
\[
n=1, \ \
\begin{pmatrix}
A & B \\
sB^{T} & C+sD
\end{pmatrix}
\begin{pmatrix}
\sigma^{1} \\
\beta^{1}
\end{pmatrix}
=
\begin{pmatrix}
G(:,1) \\
sF(:,1)+b_{1}C\beta^{0}
\end{pmatrix}
,\]
\[
n\geq 2, \ \
\begin{pmatrix}
A & B \\
sB^{T} & C+sD
\end{pmatrix}
\begin{pmatrix}
\sigma^{n} \\
\beta^{n}
\end{pmatrix}
=
\begin{pmatrix}
G(:,n) \\
sF(:,n)+C[c_{1}\beta^{n-1}+c_{2}\beta^{n-2}+\cdots+c_{n-1}\beta^{1}]+b_{n}\beta^{0}
\end{pmatrix}
,\]
where the matrices $A, B, C, D$ and $F,G$ are the  corresponding stiffness matrix, mass matrix and load vectors with respect to the multiscale basis functions. If we denote the analogous fine scale matrices by subscript $(\cdot)_{f}$ in the fine grid basis functions, then we can obtain the following relationship between coarse-grid system and fine-grid system,
\begin{equation}
\label{matrix_form}
\begin{aligned}
&A=R_{\text{off}}'A_{f}R_{\text{off}},\ \ \ B=R_{\text{off}}'B_{f}B_{\text{off}},\ \ \ C=G_{\text{off}}C_{f}G_{\text{off}}'\\
&D=G_{\text{off}}D_{f}G_{\text{off}}',\ \  G=R_{\text{off}}'G_{f},\ \ \ \ \ \ \ \ \ F=G_{\text{off}}F_{f}.
\end{aligned}
\end{equation}
Here  $G_{\text{off}}$ is the restriction operator from $Q_{H}$ to $Q_{h}$.
The fine-scale solution and the coarse-scale solution are connected through the matrix $R_{off}$.
From (\ref{matrix_form}), we can see that the size of matrix $A$ is $M_{t}\times M_{t}$, which is much smaller than $A_{f}\in \mathbb{R}^{N_{h}\times N_{h}} (M_{t}\ll N_{h})$. This
shows the multiscale  model reduction in terms of algebraic system.

\section{Numerical results}

In this section, we present a few  numerical experiments to identify the unknown parameters in the muti-term time fractional diffusion equations to evaluate the performance of the improved  implicit sampling. The mixed GMsFEM is used to build a surrogate model for the multiscale diffusion models to speed up  sampling posterior. In subsection \ref{Ex_1},  we identify the multi-term time fractional derivative orders simultaneously. In subsection \ref{Ex_2}, we recover the diffusion coefficient $k(x)$ based on the Gaussian prior when the two occasions   occur: known fractional derivatives and  unknown fractional derivatives. In subsection \ref{Ex_3}, the reaction coefficient $q(x)$ is estimated  based on   Laplace prior.

For the numerical examples, we set the spatial domain $\Omega=[0,1]\times[0,1]$ and the final time $T=1$. Measurement data are generated by the fine grid FEM with time step $\Delta t=0.02$, and the noise level is set to be $\sigma=0.01$.
The step size to compute the sensitivity matrix $\bar{H}$ using the difference method is equal to  $0.5$.

\subsection{Inversion for multi-term time fractional derivative orders}\label{Ex_1}
\begin{figure}
\centering
\subfigure[]{
\includegraphics[width=0.55\textwidth, height=0.35\textwidth]{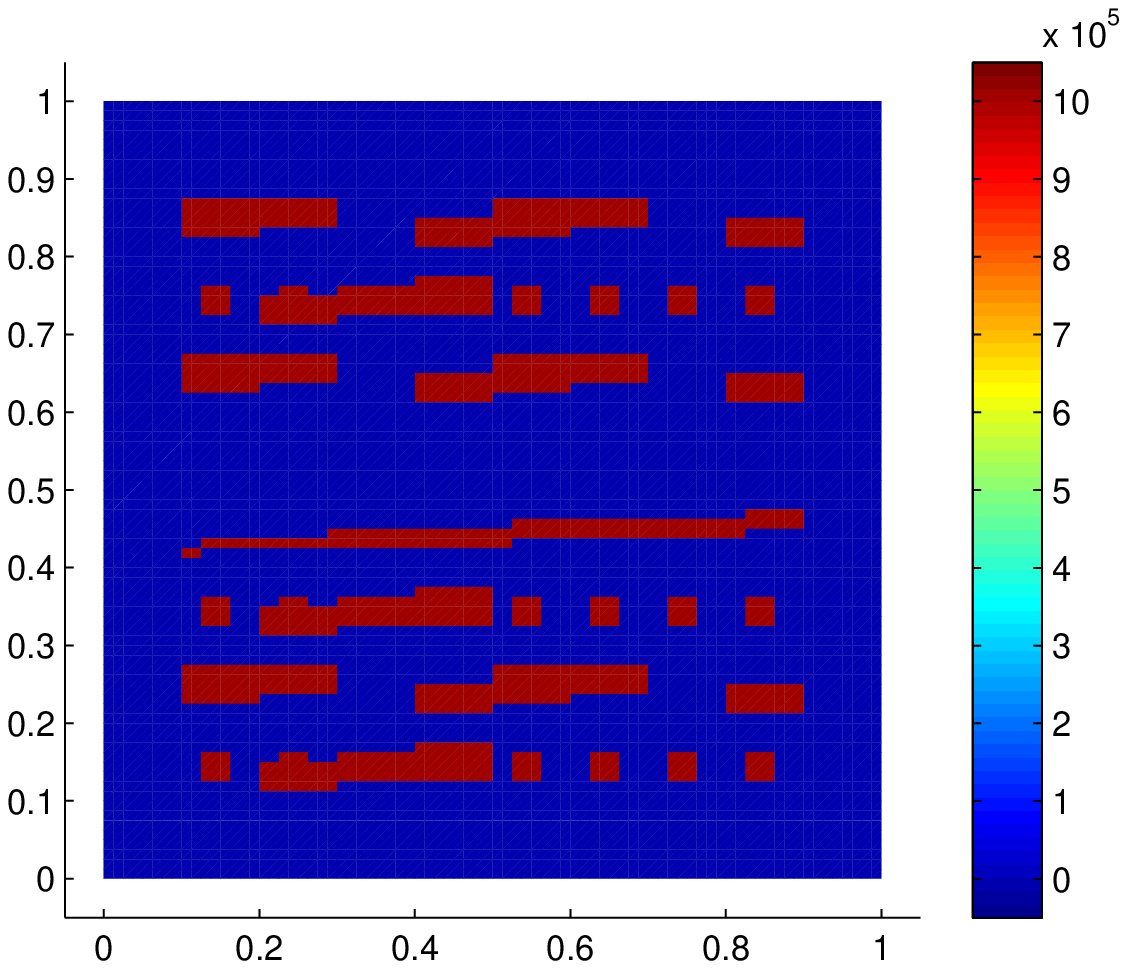}
}
\subfigure[]{
\includegraphics[width=0.37\textwidth, height=0.34\textwidth]{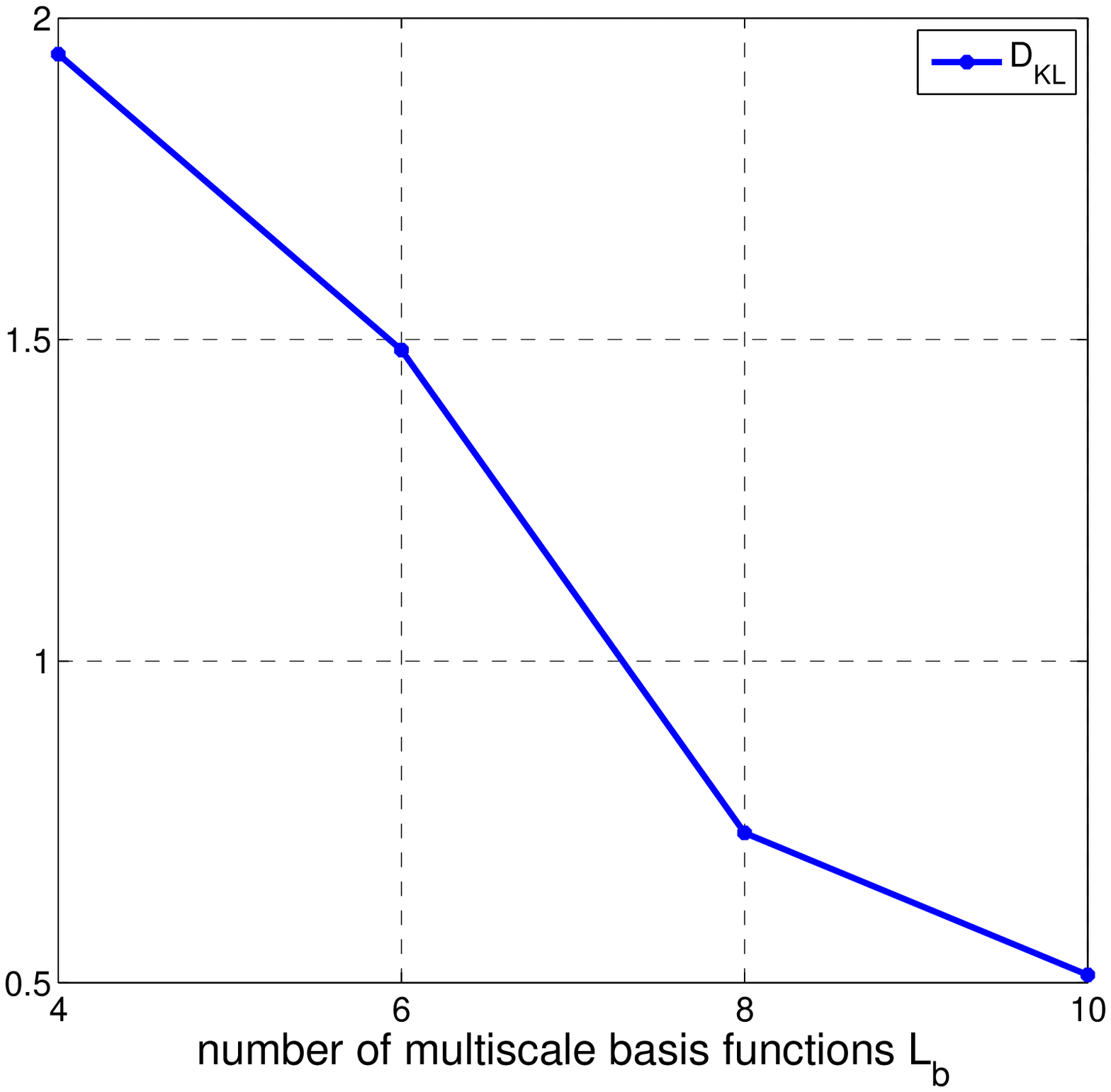}
}
\caption{(a) The spatial distribution of permeability field; (b) KL divergence between the approximate posterior density $\tilde{p}(\theta|d)$ and reference posterior density $p(\theta|d)$ with respect to the number of multiscale basis funcitons $L_{b}$.}
\label{permeability}
\end{figure}

%
In this subsection, our objective is to identify the multi-term time fractional derivatives  in model (\ref{caputa-def}), where Dirichlet boundary condition $g(x,t)\equiv 1$ on the four boundary sides, the source term $f(x,t)=10$ and the reaction coefficient $q(x)=1$. The positive constants  $\gamma_{1}=0.2$ and $\gamma_{2}=0.8$.
The forward model is discretized on a $80\times 80$ uniform fine grid, and we set the coarse grid $8\times 8$ for mixed GMsFEM computation. We consider permeability field $k(x)$ as heterogeneous media illustrated in Figure \ref{permeability} (a). We take $n=240$ Neumann observations on the left and right sides at $t=0.4$ and $t=1$.
The reference fractional derivative orders are given as $\alpha_{1}=0.3$ and $\alpha_{2}=0.6$. In the example, the fractional orders are restricted in the interval  $(0,1)$ (subdiffusion). To this end, we take a transform \cite{Jiang} as following
\begin{eqnarray*}
\label{trans}
h(x)=\frac{1}{2}+\frac{1}{\pi}\arctan(x), \ x\in \mathbb{R}.
\end{eqnarray*}
By the transformation, it can make the fractional orders lie  in the interval $(0,1)$.
We impose a Guassian prior for the unknown parameter vector $[\alpha_{1}, \alpha_{2}]$ under the framework of hierarchical Bayesian and the variance $\lambda$ of the Gaussian prior is unknown and  regarded as the hyper-parameter.
Since the dimension of unknown parameters is relatively low, it is appropriate to set $\vartheta=1$ in this example.
The initial guess for $(\alpha_{1},\alpha_{2})$ is given by $(0.5,0.5)$ to compute the MAP point by augmented-Tikhonov algorithm described  in Table \ref{a-Tikh_algorithm}.
In the example, the approximated posterior $\tilde{p}(\theta|d)$ refers to the posterior obtained by using mixed GMsFEM and the reference posterior $p(\theta|d)$ is obtained by mixed FEM in the sampling process.
To illustrate the convergence of KL divergence $D_{KL}(\tilde{p}(\theta|d) || p(\theta|d))$ \cite{Gibbs-2002, Jiang}, we plot the KL curve  with respect to the number of multiscale basis functions $L_{b}$ by using 5000 samples generated by improved implicit sampling in Figure \ref{permeability} (b).
From this figure, it can be found that the $D_{KL}$ substantially decreases as the number of multiscale basis functions increases.
In particular, we also list the CPU time in improved implicit sampling by using mixed GMsFEM and mixed FEM in Table \ref{CPU-time}. From this table, it can be seen that the CPU time decreases  significantly when mixed GMsFEM is applied in the improved implicit sampling.

\begin{table}
  \centering
\begin{tabular}{|c|c|c|c|c|}
\hline
$L_{b}$ & 4 & 6 & 8 & 10\\
\hline
$t_{off}$ & $2.2472s$ & $2.3359s$ & $2.3870s$ & $2.2800s$\\
\hline
$t_{on}$ & $3.2088\times10^{3}s$ & $7.0975\times10^{3}s$ & $1.2994\times10^{4}s$ & $2.0652\times10^{4}s$\\
\hline
$t_f$ & \multicolumn{4}{c|}{$2.7453\times10^{4}s$}\\
\hline
\end{tabular}
\caption{Comparison of CPU time of improved implicit sampling by using mixed FEM and mixed GMsFEM with respect to the number of multiscale basis functions $L_{b}$.
$t_{off}$ and $t_{on}$ are the offline time and online time by using mixed GMsFEM respectively and $t_{f}$ is the total time by using mixed FEM.
}\label{CPU-time}
\end{table}

\begin{figure}
\centering
\includegraphics[width=0.4\textwidth, height=0.35\textwidth]{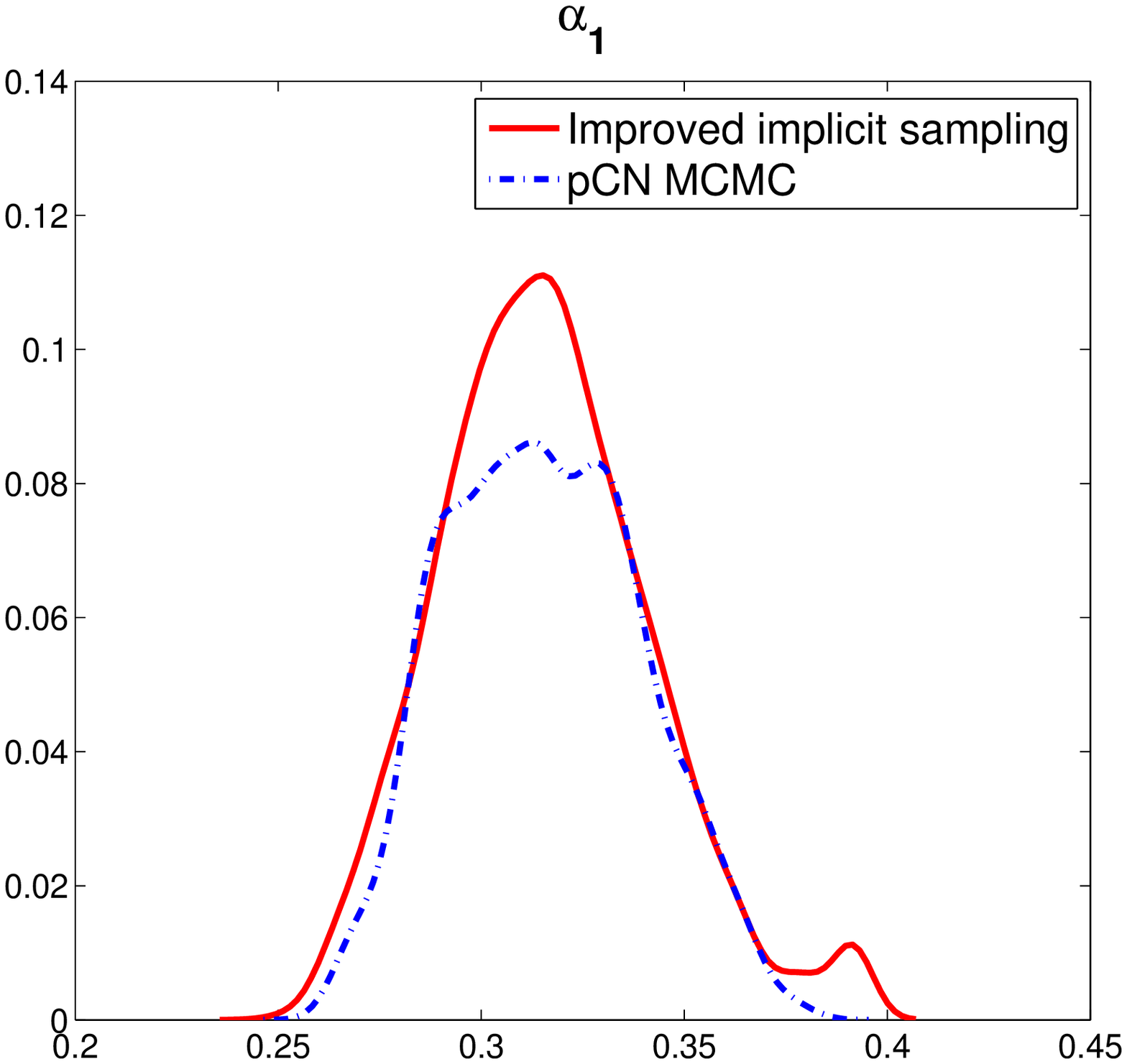}
\includegraphics[width=0.4\textwidth, height=0.35\textwidth]{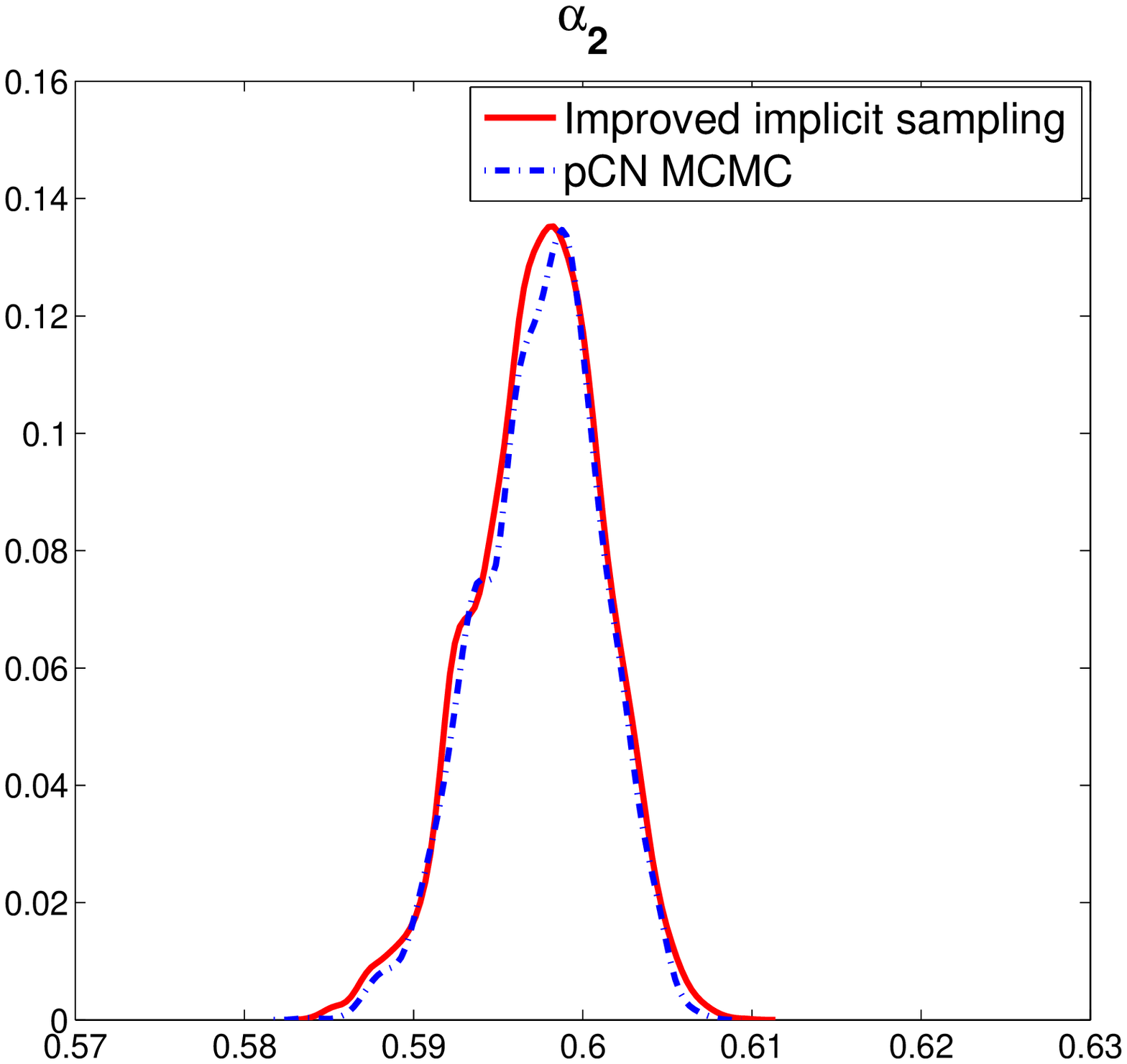}
\caption{The posterior marginal distribution for $\alpha_{1}$ (left) and  $\alpha_{2}$ (right) by improved
implicit sampling and pCN MCMC method.}\label{pdf2}
\end{figure}

\begin{figure}
\centering
\includegraphics[width=0.41\textwidth, height=0.35\textwidth]{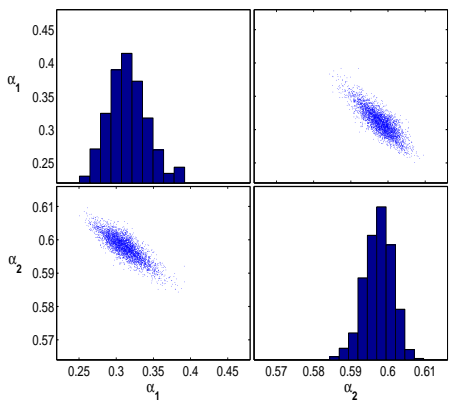}
\caption{The matrix plot for $(\alpha_{1}, \alpha_{2})$  by improved implicit sampling.}
\label{matrix_frac}
\end{figure}

In order to compare  the performance by using the improved implicit sampling  and MCMC sampling, we also use pCN MCMC method \cite{Evalution-Gaussian,nonlinear-data-assimilation} to solve the inverse problem.
We select $L_{b}=6$ multiscale basis functions at each coarse block to construct the offline space.
The tuning parameter in pCN MCMC method is set at $0.05$ with an acceptance probability $30.22\%$ to make the chain efficiently explore the parameter space.
The number of samples used in these methods is 5000. The first 1000 steps of the pCN MCMC samples are taken out as the burn-in period and the remaining 4000 samples are used to compute the corresponding statistical quantities.
In Figure \ref{pdf2}, we plot the posterior marginal distribution of the fractional orders by improved implicit sampling and pCN MCMC method, respectively. From this figure, we can see that the posterior marginal distribution of $\alpha_{2}$ produced by the two methods are similar to each other. While for $\alpha_{1}$, it seems that the samples generated by improved implicit sampling have the ability of focusing on the region with high probability than pCN MCMC method, this is  because the weights in improved implicit sampling can make the samples be close to the MAP point.
 By the support of these curves, it can be seen that the uncertainty of $\alpha_{1}$ is larger than that of $\alpha_{2}$. Furthermore, in order to explore the correlation between the two parameters, we present the matrix plot in Figure \ref{matrix_frac}. It is clearly  to see that the two parameters are strongly negatively correlated.

Next, in order to compare the sampling  performance of these two methods, we  calculate the effective sample size (ESS).
Based on (\ref{ess_w}), we can compute the ESS in improved implicit sampling. There are more than 1712 effective samples in 5000 samples, which is sufficient to verify the validity of improved implicit sampling. In Table \ref{table_mean_std}, we also show the posterior marginal mean and standard deviation of these two methods, both of which give almost the same results.
However, the above formula could not be applied to calculating  the ESS for  MCMC \cite{rmap_Bui-Thanh}. Alternatively,  we adopt the following formula to do some qualitative analysis, i.e.,
\begin{equation}
\label{ess_mcmc}
\text{ESS}=\frac{N_{s}}{\bar{\rho}},
\end{equation}
where $N_{s}$ is the total number of samples and $\bar{\rho}$ is the averaged integrated auto-correlation time (IACT).
For a model with $m$ unknown parameters, IACT $\bar{\rho}$ can be calculated  by
\[
\bar{\rho}=\frac{1}{m}\sum_{i=1}^{m}(1+2\sum_{k=1}^{\infty}\rho_{k}),
\]
 where  $\rho_{k}$ is the auto-correlation function (ACF) at lag $k$.  For a time series $M_{t}$ with mean $\mu_{0}$ and variance $\sigma_{0}^2$, $\rho_{k}$ is defined as
 \[
 \rho_{k}=\frac{E[(M_{t}-\mu_{0})^{T}(M_{t+k}-\mu_{0})]}{\sigma_{0}^2}.
 \]

\begin{figure}
\centering
\includegraphics[width=0.41\textwidth, height=0.35\textwidth]{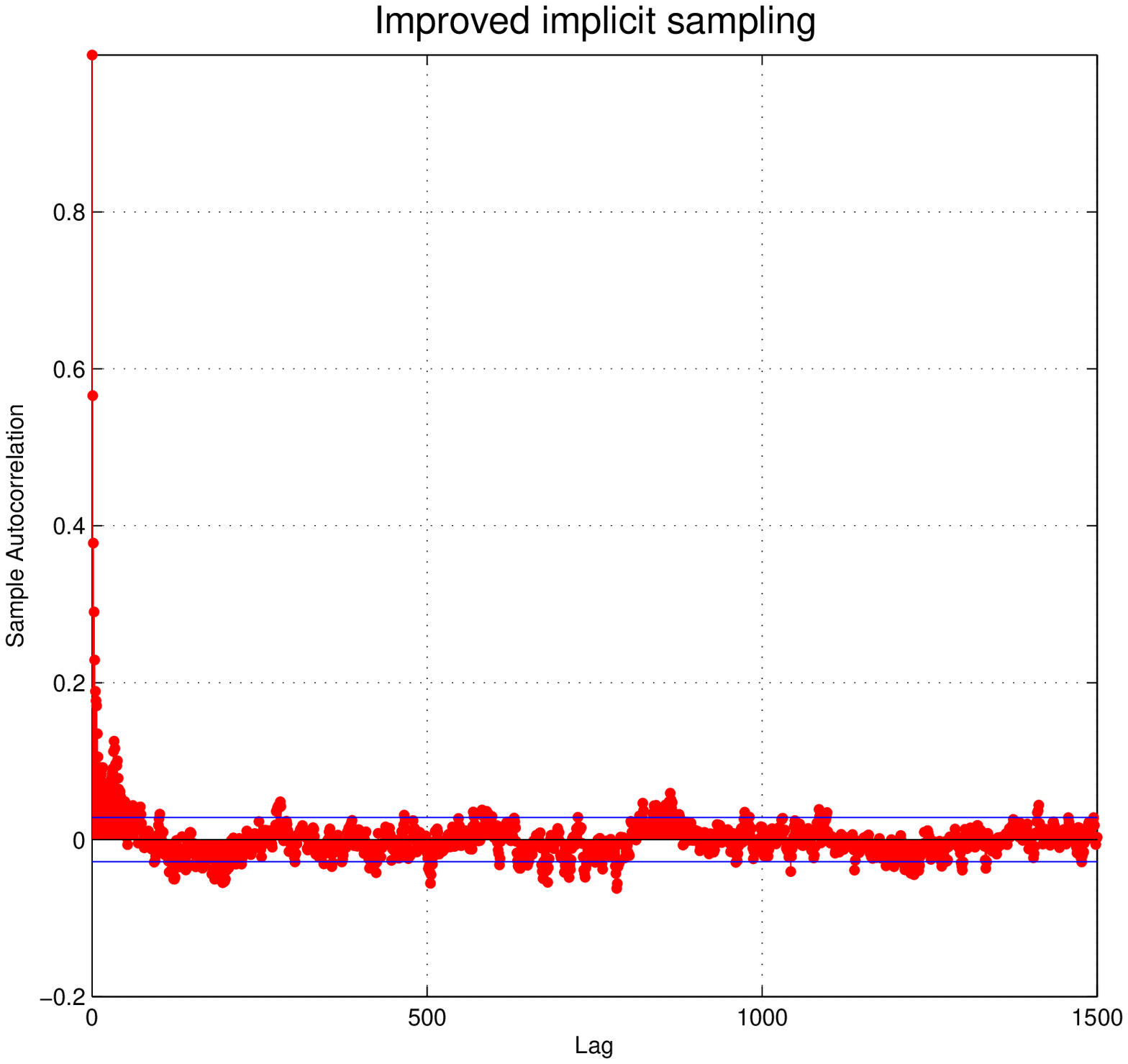}
\includegraphics[width=0.41\textwidth, height=0.35\textwidth]{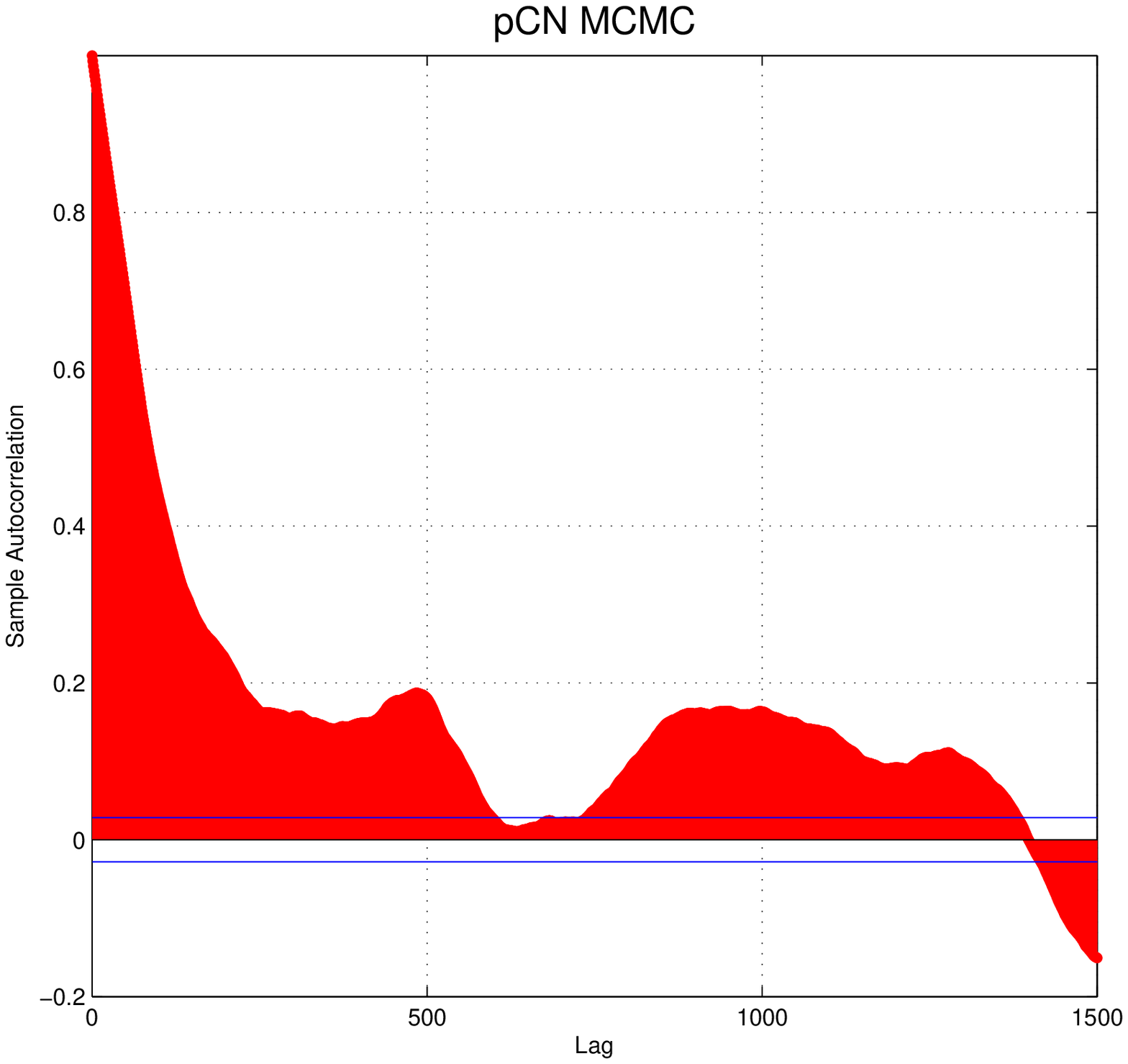}
\caption{The ACF of $\alpha_{1}$  for improved  implicit sampling (left) and pCN MCMC method (right).}\label{acf}
\end{figure}

The difference between  ESS and the actual sample size is due to the large autocorrelation structure imposed by the MCMC method. In Figure \ref{acf}, we plot the ACF of these two methods. We can see that the ACF of improved impilcit sampling gets  closer to 0 as the lag increases. However, it can be seen that the ACF of MCMC  method is much larger than that of improved implicit sampling and has no convergence tendency.
The figure implies that the ESS of MCMC method is much smaller than improved implicit sampling based on (\ref{ess_mcmc}), though  the exact number of the ESS of MCMC method is not sure.

\begin{figure}
\centering
\includegraphics[width=0.41\textwidth, height=0.35\textwidth]{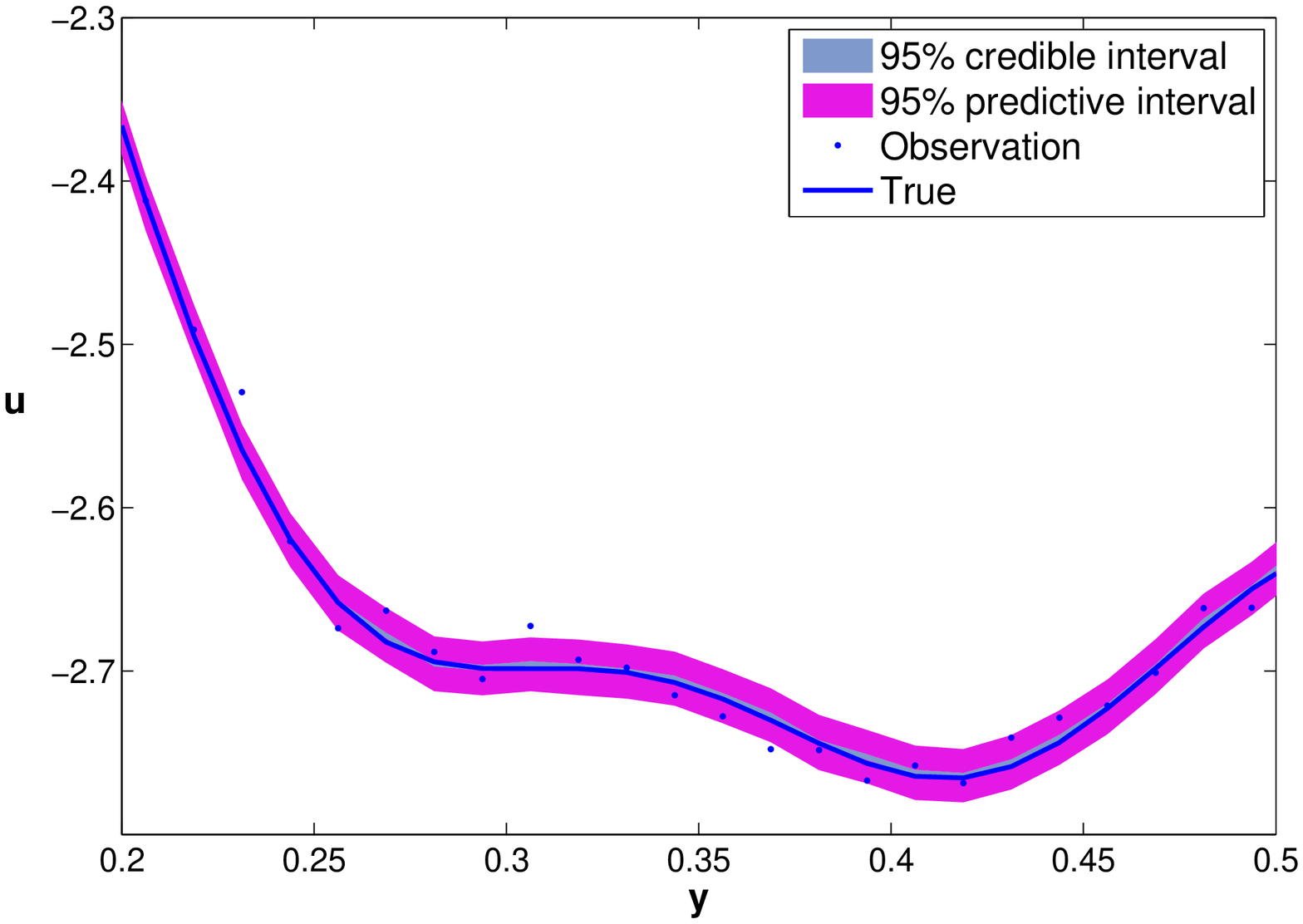}
\includegraphics[width=0.41\textwidth, height=0.35\textwidth]{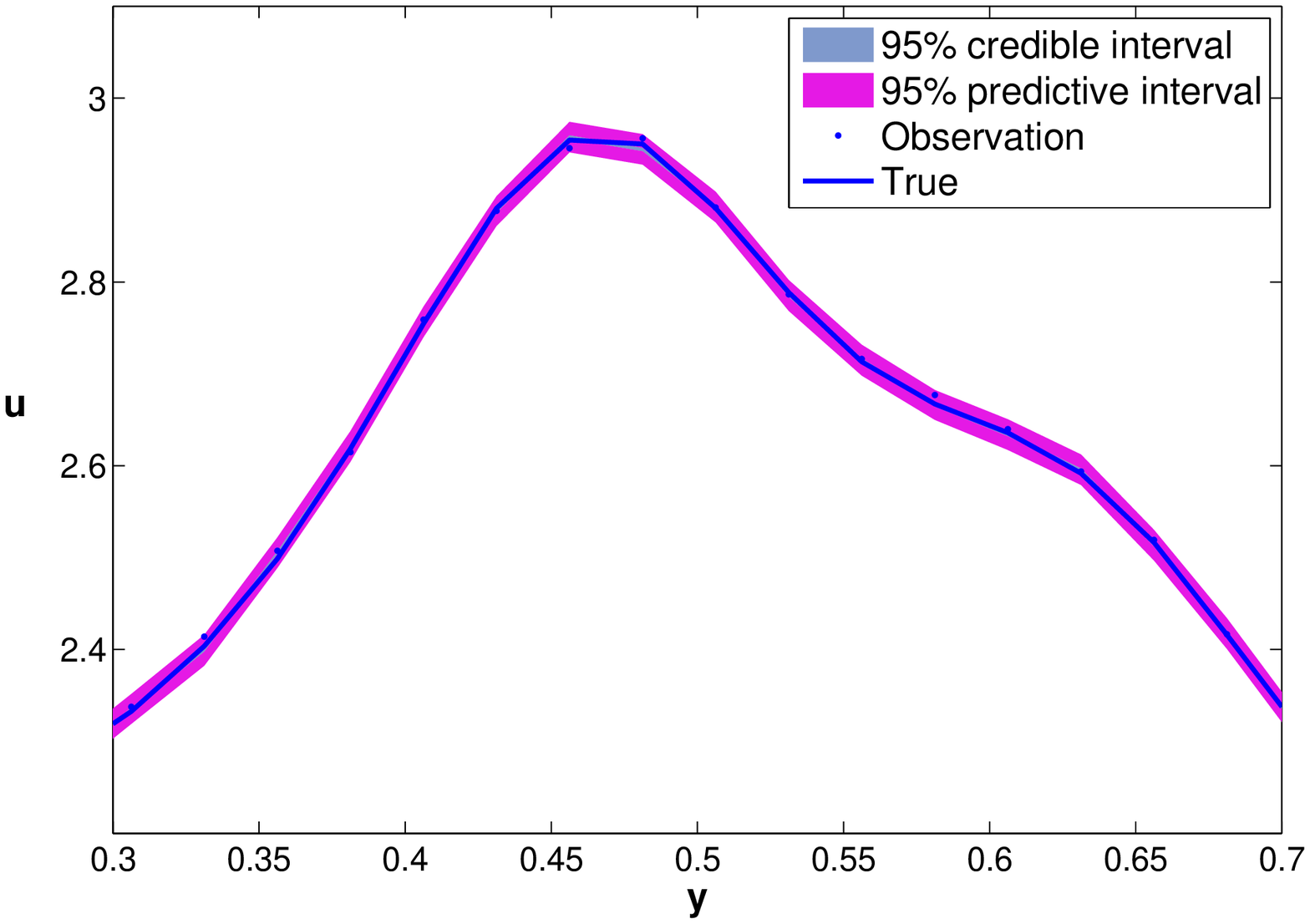}
\caption{Data, true, and $95\%$ credible interval  and prediction interval  by improved implicit sampling for $u(0,y; 0.4)$ (left) and $u(1,y; 0.4)$ (right).}
\label{cre_frac}
\end{figure}

The observations and the actual model response are plotted in Figure \ref{cre_frac}. To construct $95\%$ credible and prediction intervals for model response at $u(0,y; 0.4)$ and $u(1,y;0.4)$,   we use $5000$ samples to produce realizations of the model. For $u(0,y; 0.4)$, we note that the uncertainty focus on the region $[0.3,0.5]$, which is mainly due to the Dirichlet boundary condition imposed in the forward model.
It can be also seen from these two figures that most of observations lie in  the predictive interval. This implies  that the improved implicit sampling properly describes  the uncertainty propagation  for the model.

\begin{table}
  \centering
\begin{tabular}{|c|c|c|}
\hline
 & Improved  implicit sampling  &  MCMC method\\
\hline
 Mean & (0.3160, 0.5974) & (0.3155, 0.5974) \\
\hline
Standard deviation & (0.0262, 0.0038) & (0.0235, 0.0036) \\
\hline
\end{tabular}
\caption{The posterior marginal mean, standard derivation for $[\alpha_{1}, \alpha_{2}]$ of improved implicit sampling and pCN MCMC method.}\label{table_mean_std}
\end{table}

\subsection{Inversion for diffusion coefficient}\label{Ex_2}
In this subsection, we apply the improved  implicit sampling to identify the diffusion coefficient $k(x)$ in model (\ref{model-fpde}). Assume that $\gamma_{1}=\gamma_{2}=1$ and the fractional orders $\alpha_{1}=0.3, \ \alpha_{2}=0.6$ are fixed.
The  reaction coefficient $q(x)=1$ is fixed  in this example.
For this example, we  set the source term $f=10$ and the Dirichlet boundary condition $g=1$.
In order to reduce the dimension of the grid-based unknown input $k(x)$,
we regard it as a random field $k(x,\varrho)$ and assume the log permeability field $\log[k(x,\varrho)]$ can be represented by the Karhunen-Lo$\grave{e}$ve expansion \cite{kle-book}.
In particular, we consider a correlated second-order Guassian random field $k(x,\varrho)$ \cite{kle-book} defined for $x\in\Omega$ with the expectation  $\text{E}[k(x,\varrho)]$ and covariance function $C(x;y)$. For $k(x,\varrho)$, it can be represented as
\[
\log k(x,\varrho)=\text{E}[k(x,\varrho)]+\sum_{j=1}^{\infty}\sqrt{\lambda_{j}}
\zeta_{j}(x)v_{j}(\varrho),
\]
where $\lambda_{j}$ and $\zeta_{j}$ are the eigenvalues and orthonormal eigenfunctions of the covariance function $C(x;y)$ satisfying
\[
\int_{\Omega}C(x;y)\zeta_{j}(y)dy=\lambda_{j}\zeta_{j}(x)
\]
for $x\in \Omega$.
We note  that the random variables $v_{j}(\varrho)$ are centered and uncorrelated.
In addition, we assume  the covariance function
\[
C(x_{1},x_{2};y_{1},y_{2})=\rho^2\exp(-\frac{|x_{1}-y_{1}|^{2}}{2 l_{1}^2}
-\frac{|x_{2}-y_{2}|^{2}}{2 l_{2}^2}).
\]
In order to reduce the dimension of the unknown parameters, the truncated expansion based on the decay rate of the eigenvalues $\lambda_{j}$ is usually employed for numerical simulation, which can be expressed as
\begin{equation}
\label{kle}
\log k(x,\varrho)=\text{E}[k(x,\varrho)]+\sum_{j=1}^{l}\sqrt{\lambda_{j}}
\zeta_{j}(x)v_{j}(\varrho)
\end{equation}
for $l\ll m$. Finally, one can rewrite (\ref{kle}) in a matrix formulation as
\[
\log k(x,\varrho)=\Phi\textbf{v}+\kappa,
\]
where $\kappa=\text{E}[\theta(x,\varrho)]$ and $\Phi=[\sqrt{\lambda_{1}}\zeta_{1},\cdots,\sqrt{\lambda_{l}}\zeta_{l}]$.
As a result, the unknown parameter $\textbf{v}\in \mathbb{R}^{l}$ needs to be   estimated under the  Bayesian framework.

\begin{figure}
\centering
\includegraphics[width=0.41\textwidth, height=0.35\textwidth]{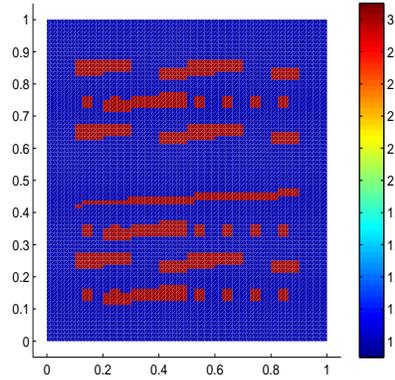}
\caption{The spatial distribution of $\text{E}(\log k(x,\varrho))$.}\label{prior_mean}
\end{figure}

  As a prior information, $\text{E}(\log [k(x,\varrho])$ is illustrated in Figure \ref{prior_mean}. The corresponding parameters in the covariance function  are given by $\rho=0.5$, and $l_{1}=l_{2}=0.1$. To reduce the dimensional of the unknown coefficient, we retain the first 12 terms based on the decay of the eigenvalues. The forward model is computed on a uniform $80\times 80$ fine grid and the coarse grid is  $5\times5$ for mixed GMsFEM. We choose $L_{b}=6$  multiscale basis functions
to construct the offline space. When we compute the matrix $R_{\text{off}}$, the number of training parameters used in this example is $12$. We choose the whole boundary Neumann data at $t=0.6$ as measurements.

\begin{figure}
\centering
\subfigure[]{
\includegraphics[width=0.3\textwidth, height=0.3\textwidth]{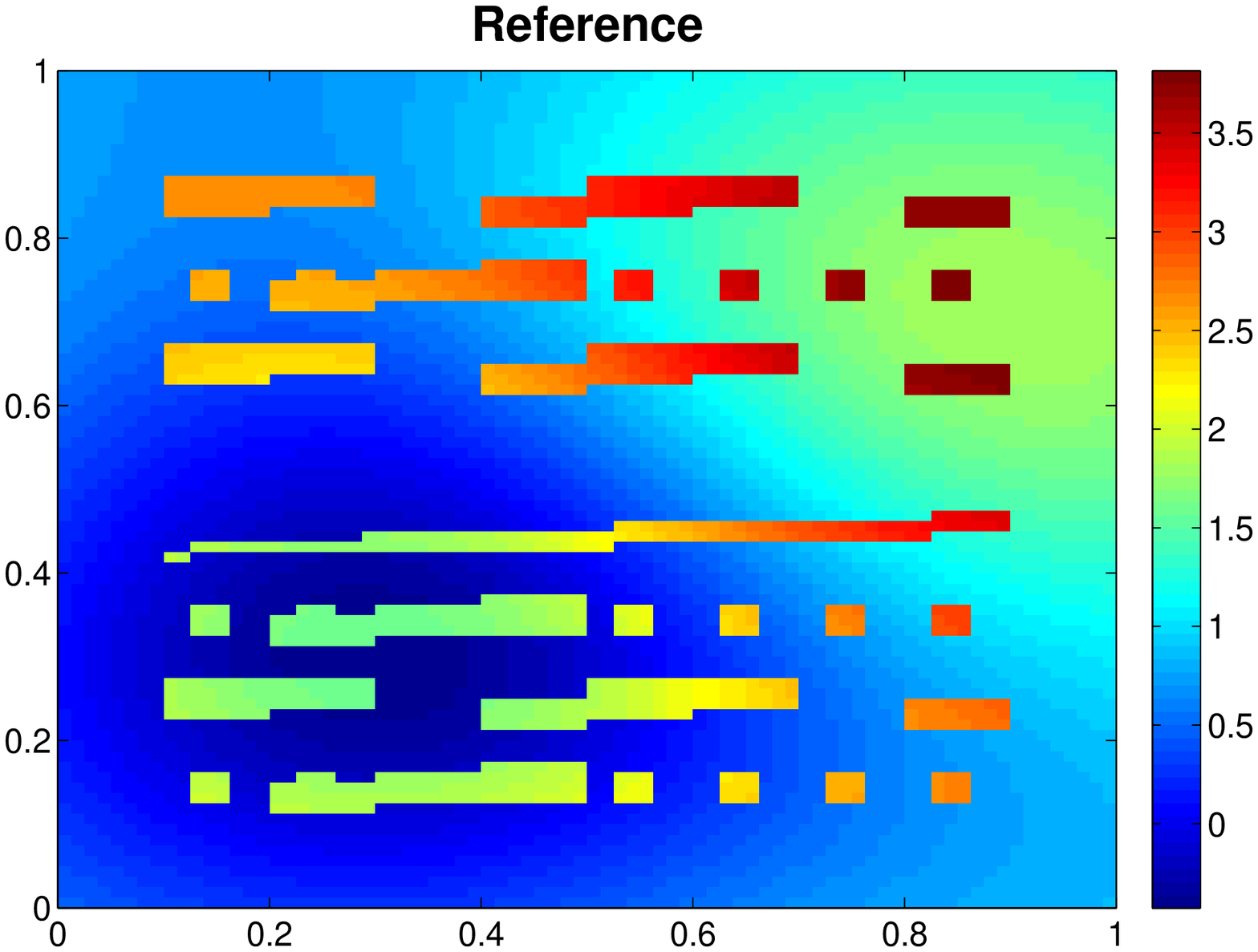}
}
\subfigure[]{
\includegraphics[width=0.3\textwidth, height=0.3\textwidth]{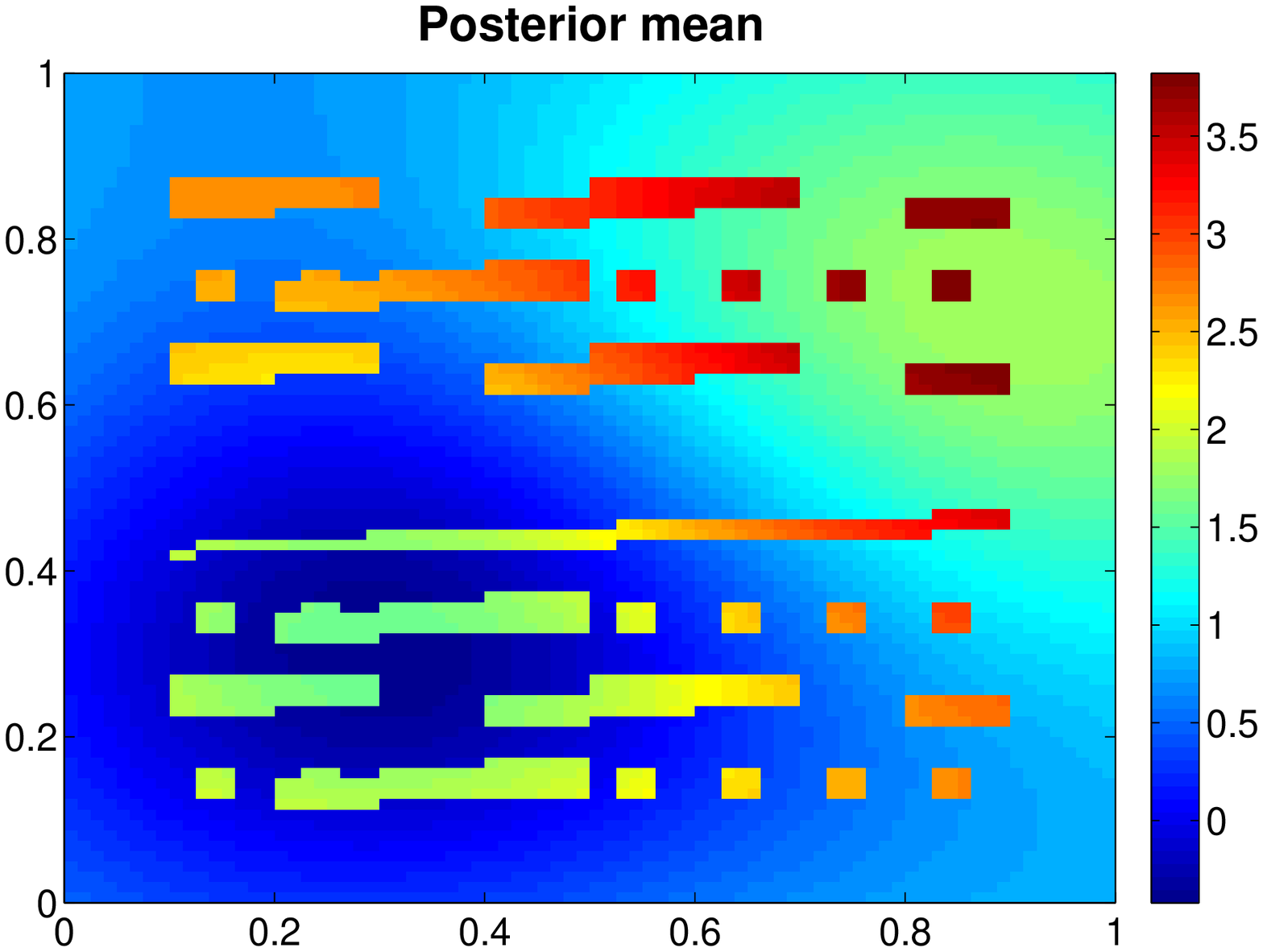}
}
\subfigure[]{
\includegraphics[width=0.3\textwidth, height=0.3\textwidth]{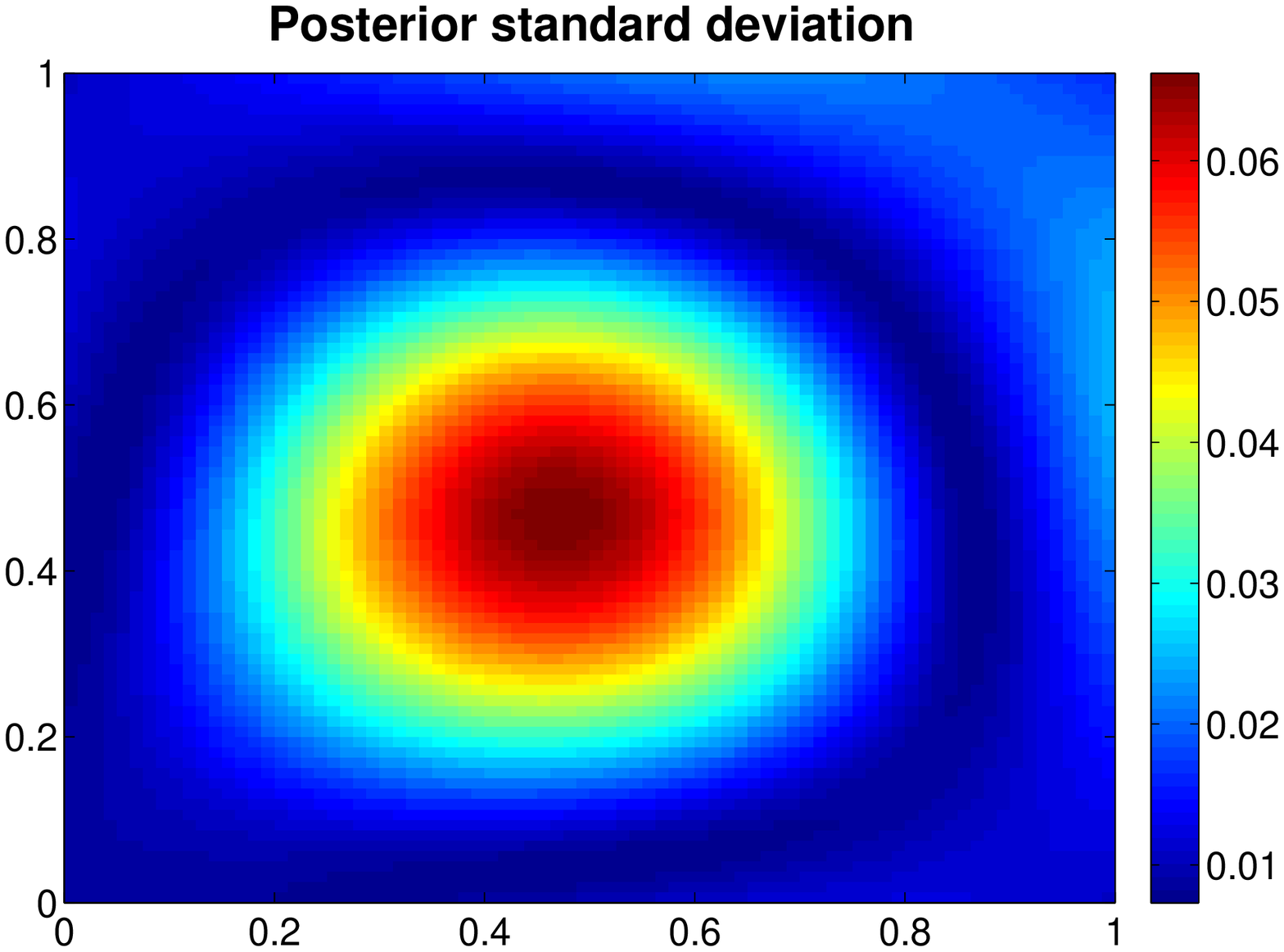}
}
\caption{The reference solution (left), the posterior mean (middle) and  posterior standard deviation (right) of $\log k(x,\varrho)$ by improved  implicit sampling.}\label{ref}
\end{figure}

In Figure \ref{ref}, we plot the spatial distribution of the reference, the corresponding posterior mean and the standard deviation in a logarithmic scale by the improved implicit sampling with $5000$ samples. The scale parameter $\vartheta$ is fixed at 1.
From these figures, it can be seen that the posterior mean can approximate the reference well. We note  that the standard deviation at the center of the physical domain is large. This may be  because that only the  boundary Neumann data are used to recover the unknown coefficient.

\begin{figure}
\centering
\subfigure[]{
\includegraphics[width=0.3\textwidth, height=0.3\textwidth]{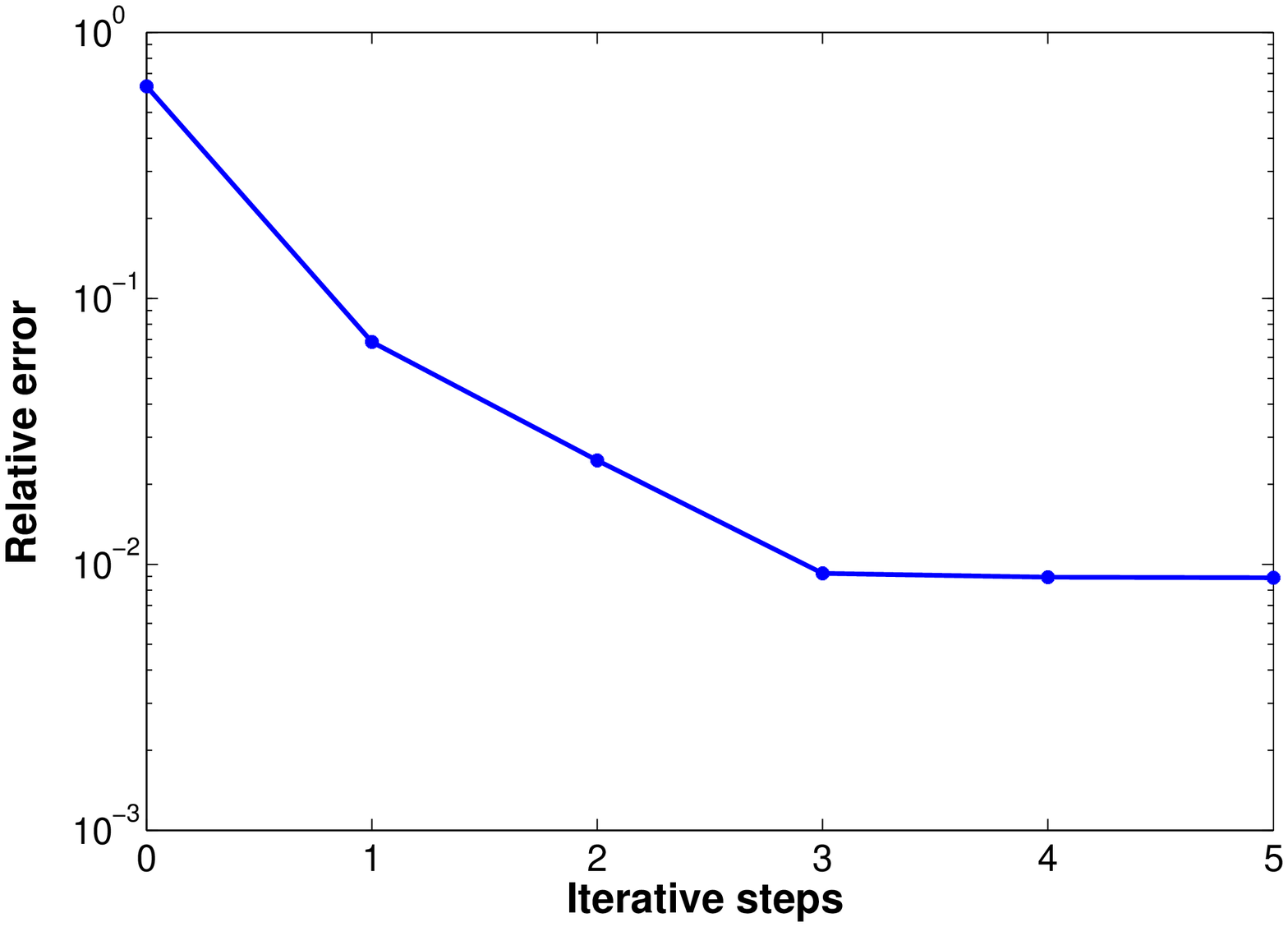}
}
\subfigure[]{
\includegraphics[width=0.3\textwidth, height=0.3\textwidth]{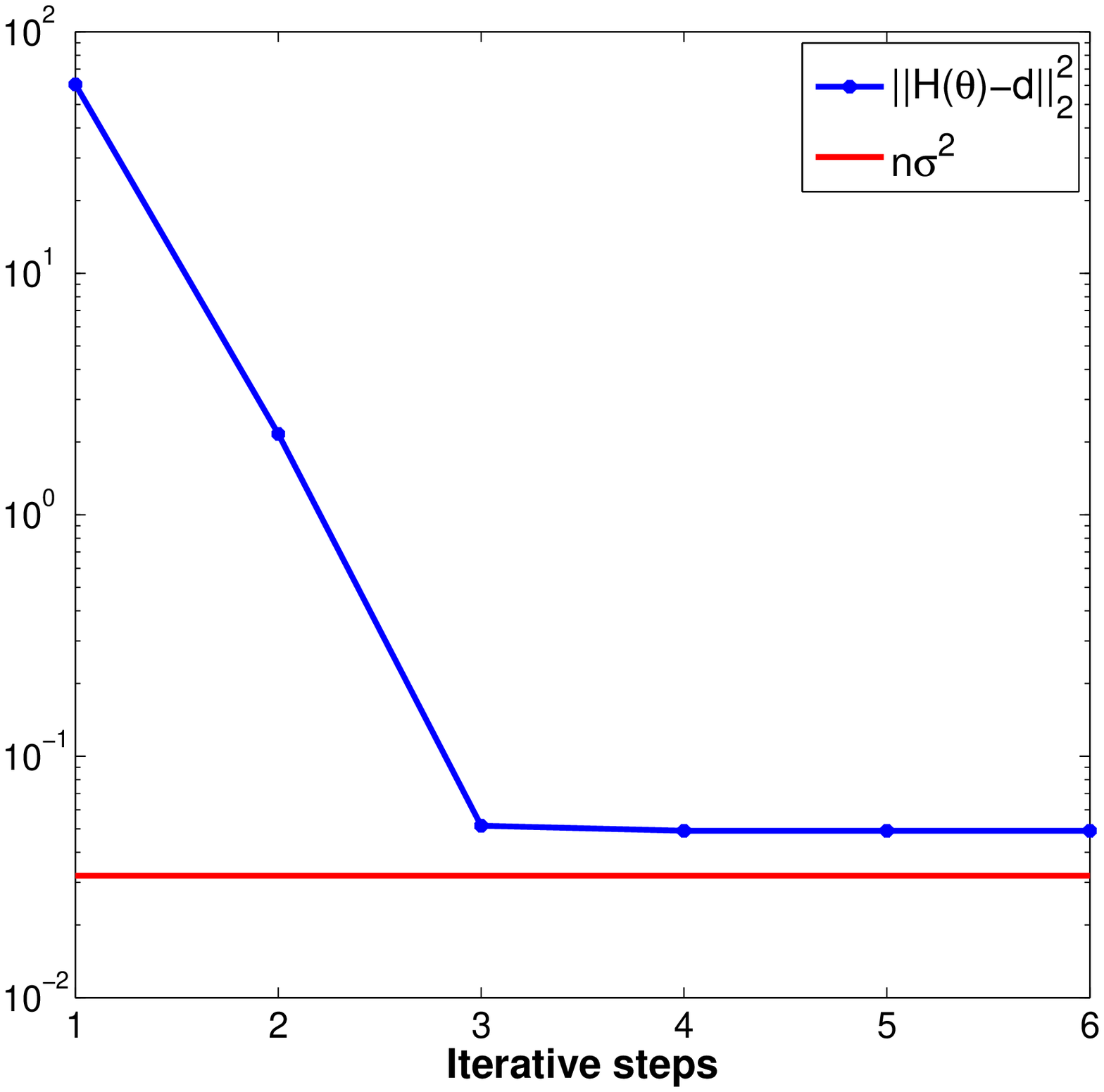}
}
\subfigure[]{
\includegraphics[width=0.3\textwidth, height=0.3\textwidth]{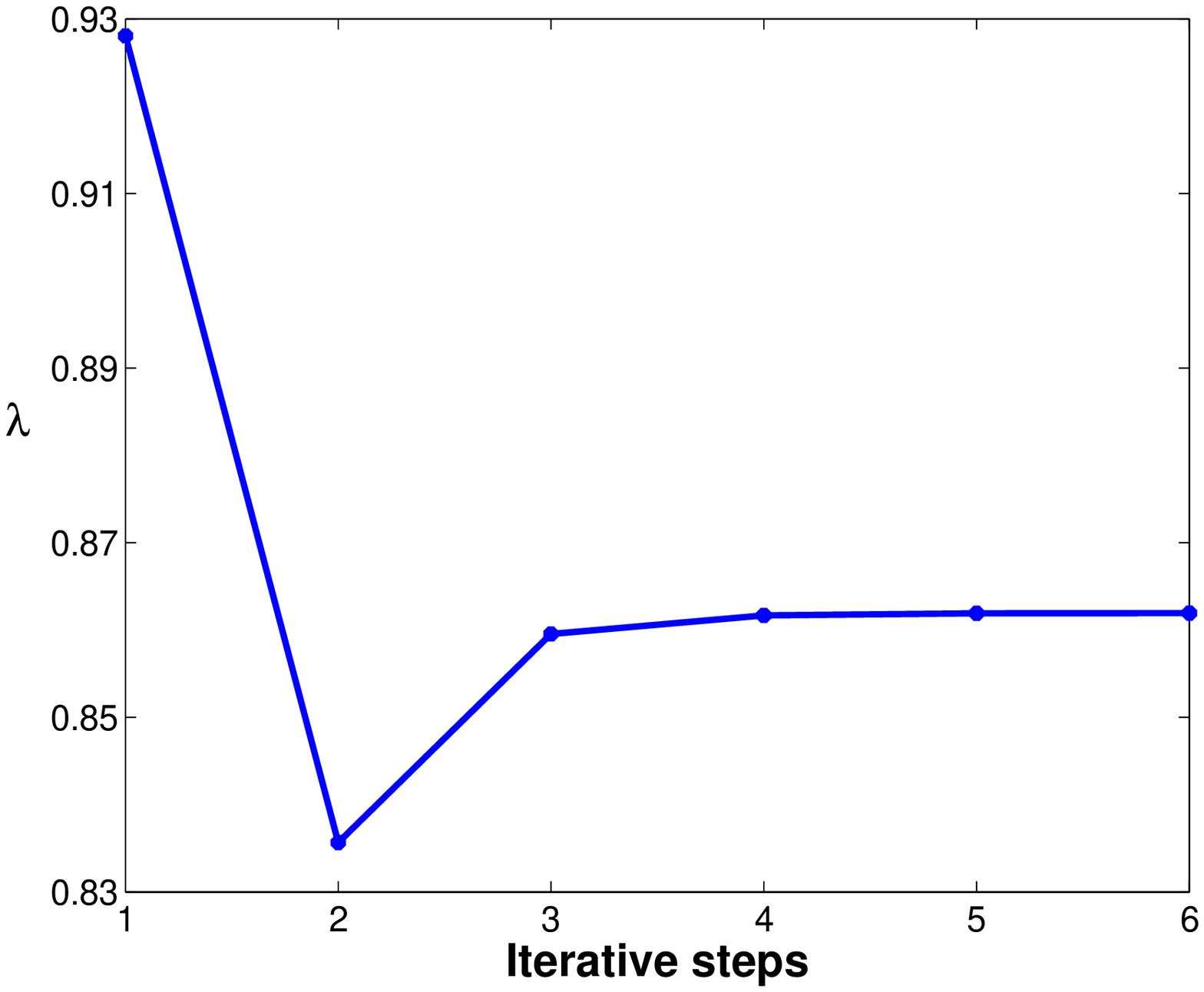}
}
\caption{Numerical results  of augmented-Tikhonov method to compute the MAP point. }\label{MAP}
\end{figure}

We use augmented-Tikhonov method to compute MAP point.  The convergence of the method  is shown  in Figure \ref{MAP}. From Figure \ref{MAP} (a), it can be seen that the relative errors  rapidly decrease  in the first three iterations. Meanwhile, from Figure \ref{MAP} (c), we can find that the hyper-parameter $\lambda$ exhibits a stable convergence trend after a few iterations, which implies the effectiveness  of augmented-Tikhonov method. As we know, the regularization parameter plays a crucial role in inverse problems, and  there are a few  strategies to select it. In order to compare the hierarchical Bayesian model  with discrepancy principle \cite{CR.Vogel-2002}, a method  of  selecting regularization parameters, we plot $\|H(\textbf{v})-d\|^{2}_{2}$ and $n\sigma^2$ versus the iterative steps in Figure \ref{MAP} (b). As we expect, $\|H(\textbf{v})-d\|^{2}_{2}$  decreases fast firstly  and then becomes stable after the third  iterative step. Though the final value of $\|H(\textbf{v})-d\|^{2}_{2}$ is not exactly equal to $n\sigma^2$, the regularization parameter $\mu=\lambda\sigma^2$ determined by the hierarchical Bayesian model can be also regarded as an approximate solution of the following nonlinear equation
\[
\|H(\textbf{v}_{\mu})-d\|^{2}_{2}-n\sigma^2=0,
\]
which is the practical implementation of discrepancy principle.

\begin{figure}[]
\centering
\subfigure[]{
\includegraphics[width=0.3\textwidth, height=0.3\textwidth]{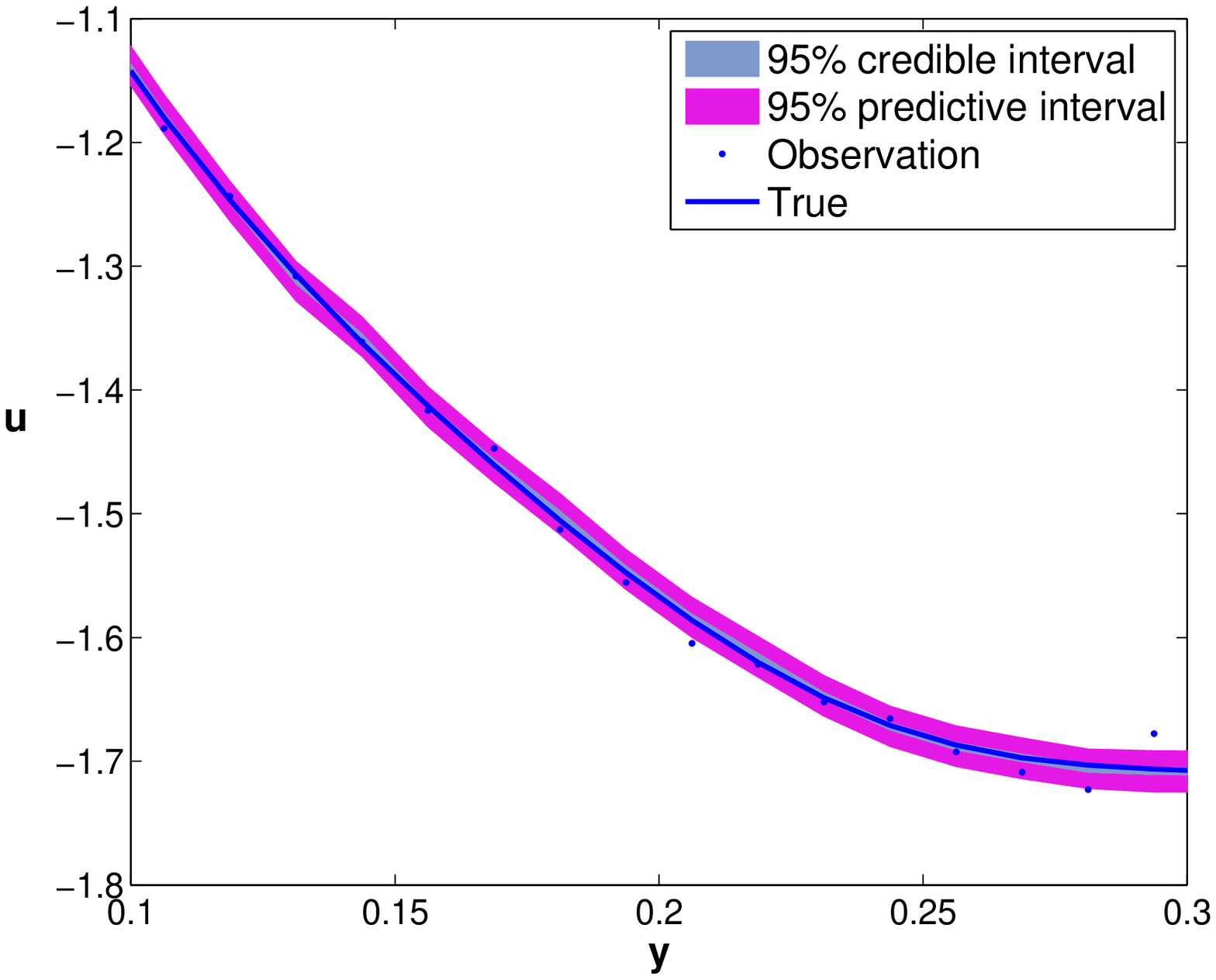}}
\quad
\subfigure[]{
\includegraphics[width=0.3\textwidth, height=0.3\textwidth]{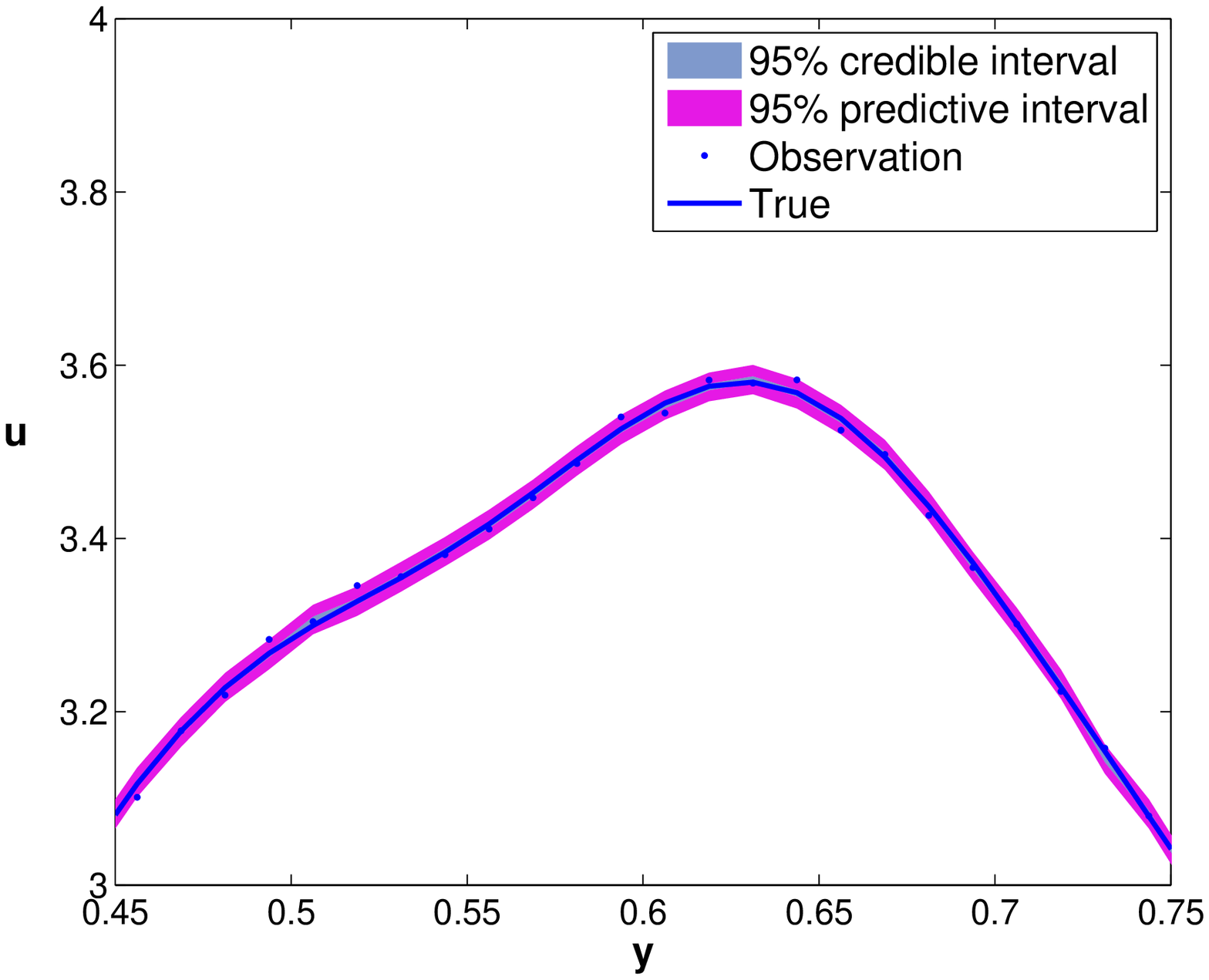}
}
\quad
\subfigure[]{
\includegraphics[width=0.3\textwidth, height=0.3\textwidth]{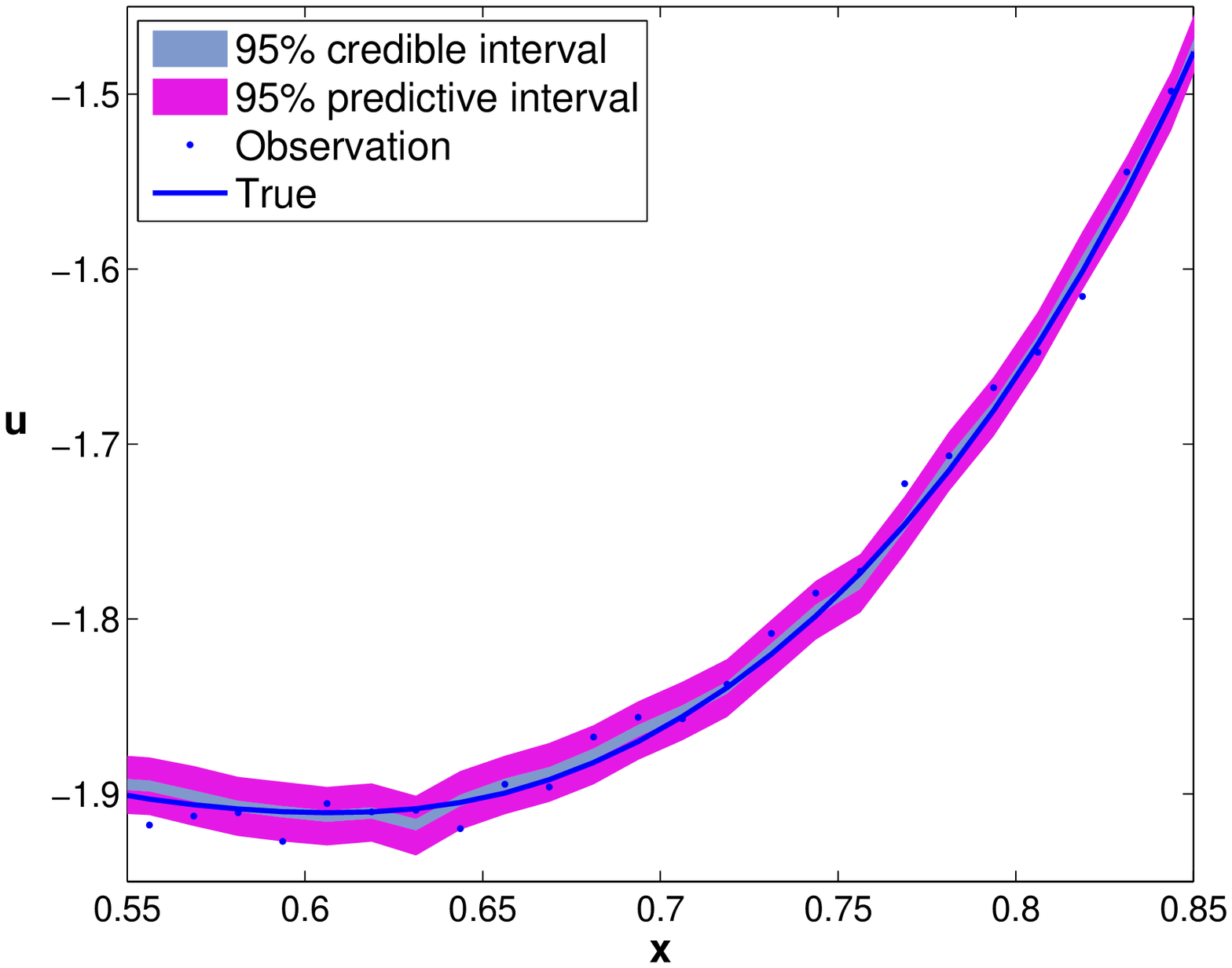}
}
\quad
\subfigure[]{
\includegraphics[width=0.3\textwidth, height=0.3\textwidth]{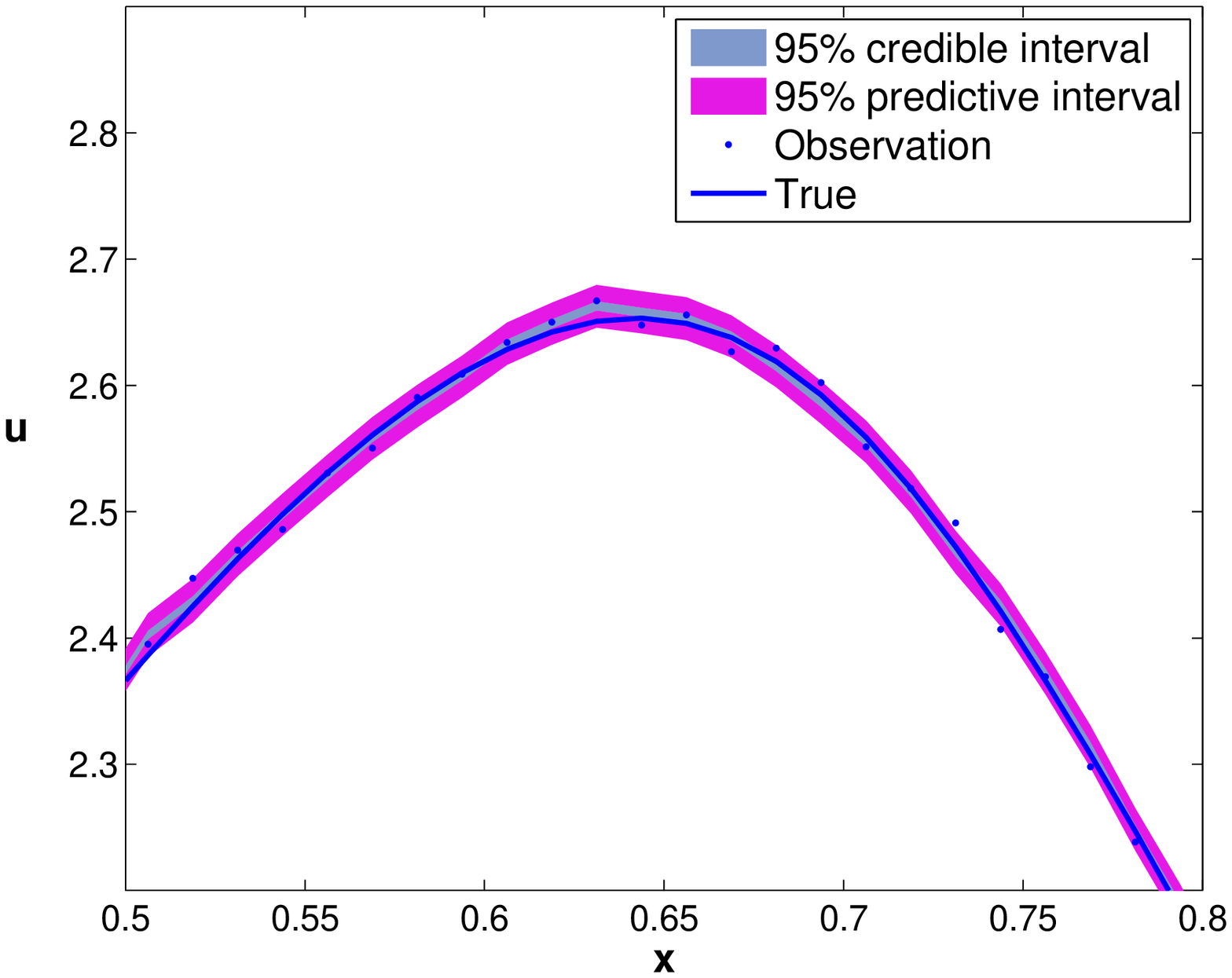}
}
\caption{Data, true, and $95\%$ credible interval  and prediction interval  by improved implicit sampling for (a) $u(0,y; 0.6)$, (b) $u(1,y; 0.6)$, (c) $u(x,0; 0.6)$, (d) $u(x,1; 0.6)$.}
\label{cre_2}
\end{figure}

To construct $95\%$ credible intervals for the model response, $5000$  model realizations  generated by improved implicit sampling are used here. The measurement error is incorporated to construct the prediction intervals. The credible interval  and predictive interval, together with the true model response and data, are plotted in Figure \ref{cre_2}. From this plot, we can see that the synthesized data are almost contained  in the prediction intervals.
 As illustrated in Figure \ref{cre_2}, the prediction interval gets wider as $y$ increases for $u(0,y; 0.6)$ in $[0.1, 0.3]$. This is because of the imposed  Dirichlet boundary  condition.

\begin{figure}
\centering
\subfigure[]{
\begin{minipage}[b]{0.17\textwidth}
\includegraphics[width=1\textwidth,height=1\textwidth]{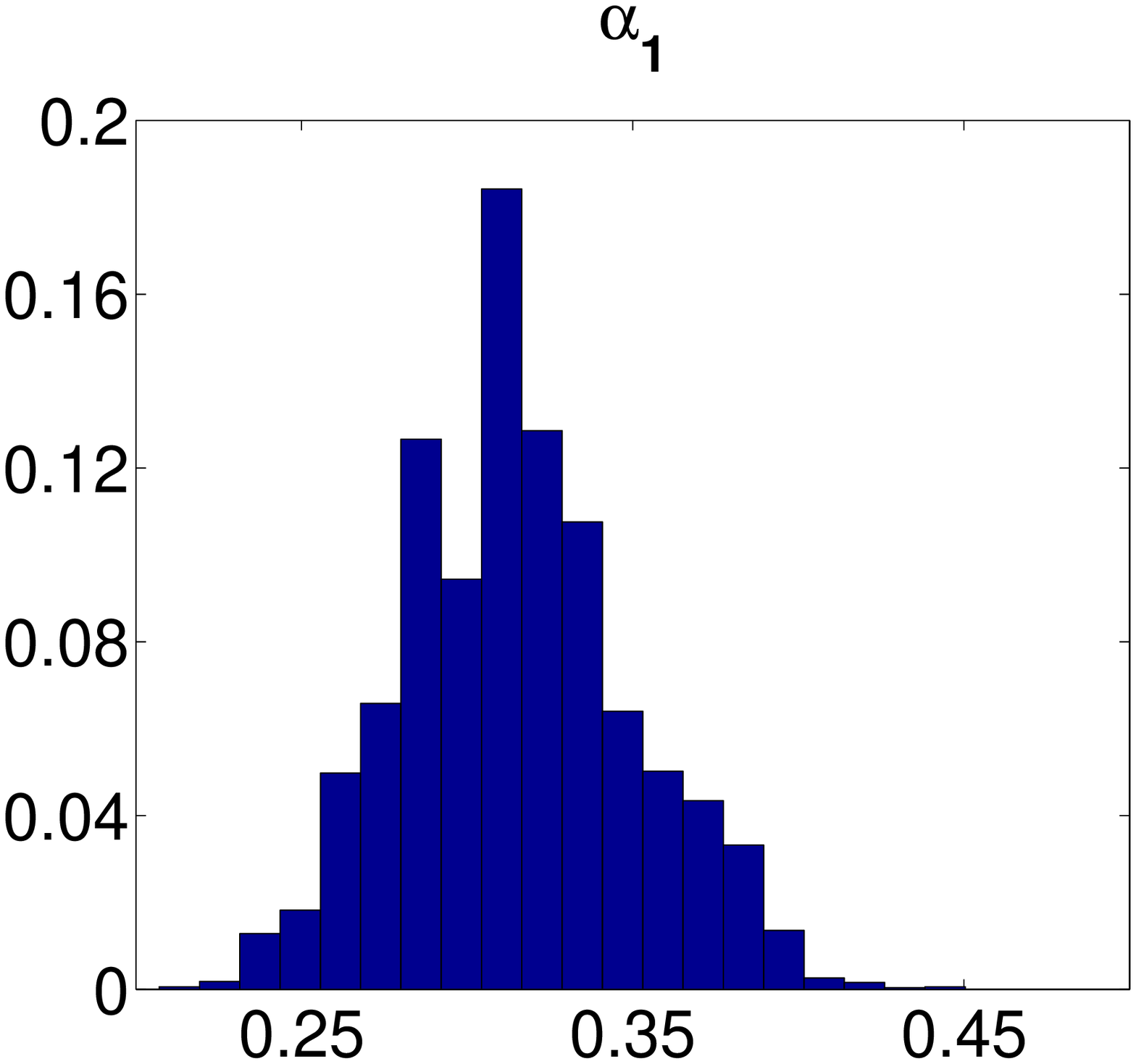} \\
\includegraphics[width=1\textwidth,height=1\textwidth]{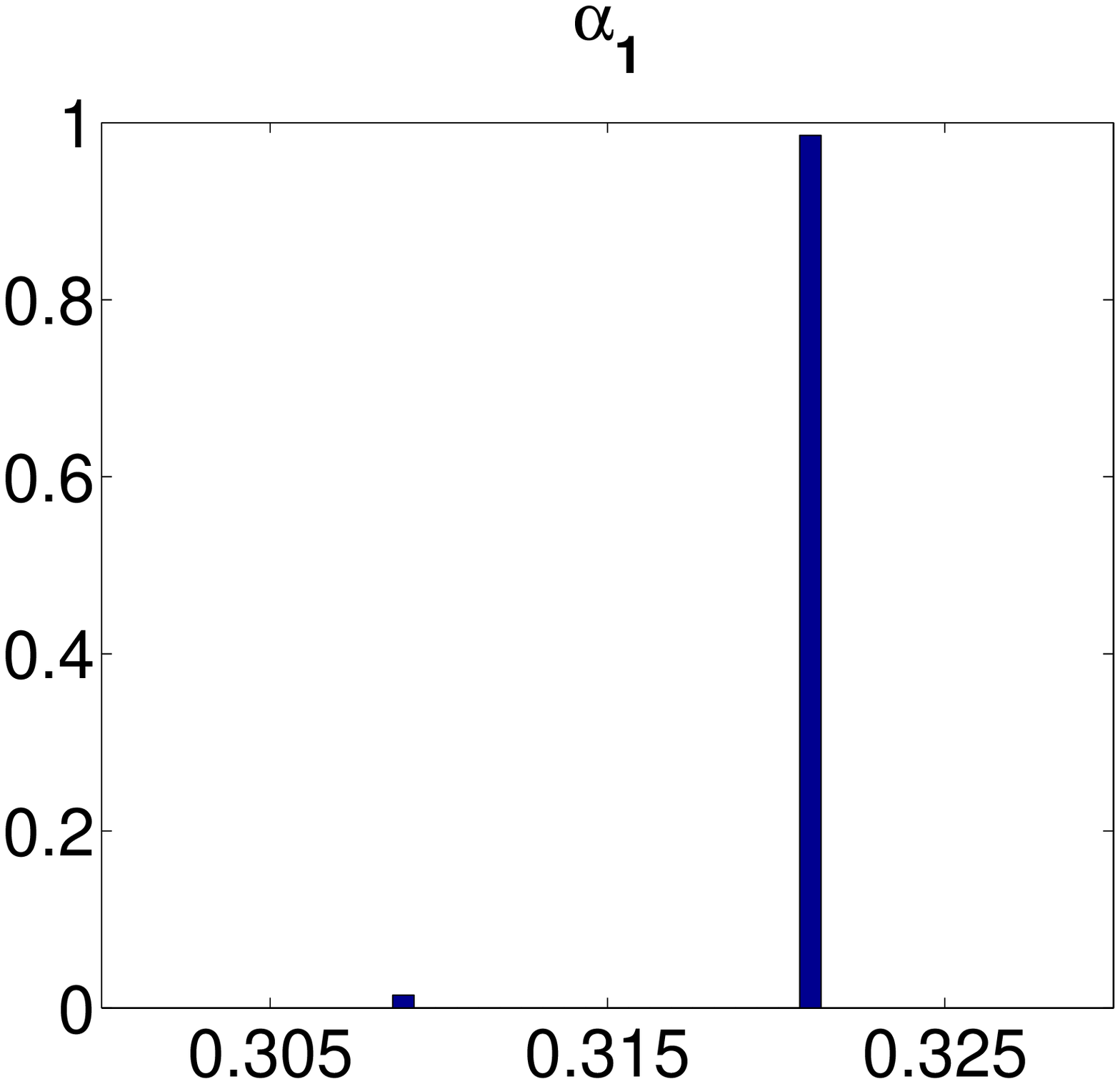}
\end{minipage}
}
\subfigure[]{
\begin{minipage}[b]{0.17\textwidth}
\includegraphics[width=1\textwidth,height=1\textwidth]{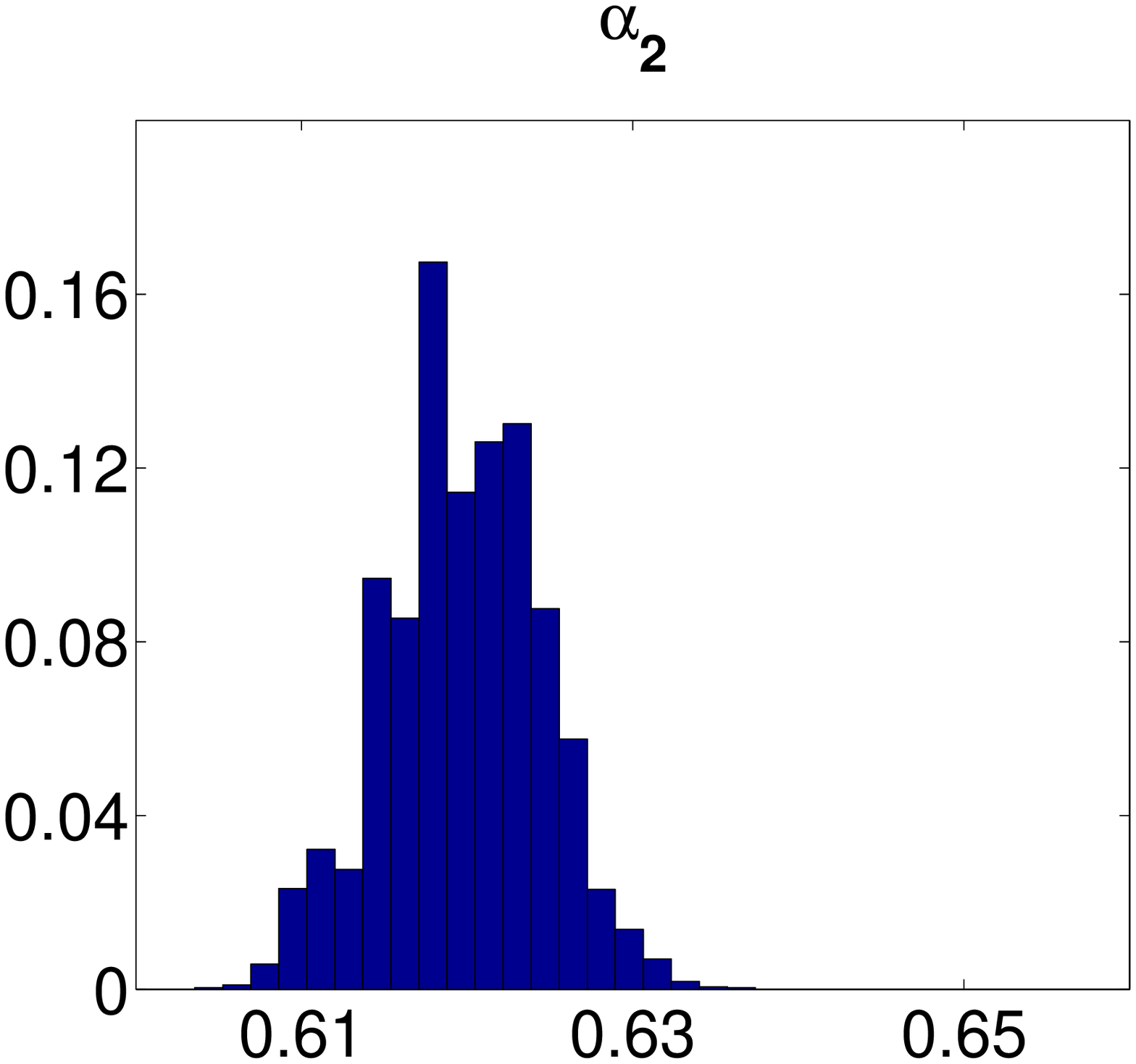} \\
\includegraphics[width=1\textwidth,height=1\textwidth]{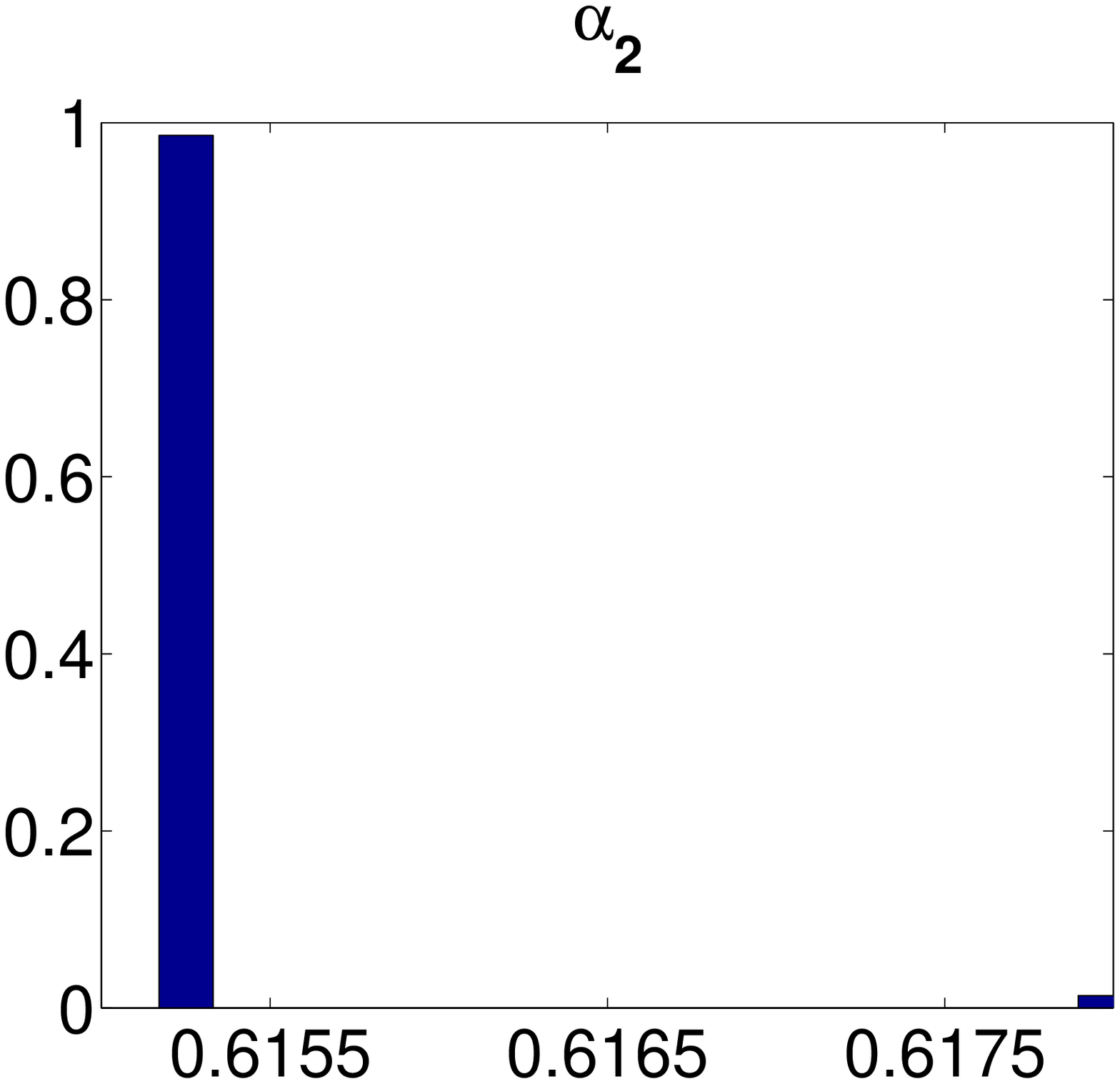}
\end{minipage}
}
\subfigure[]{
\begin{minipage}[b]{0.17\textwidth}
\includegraphics[width=1\textwidth,height=1\textwidth]{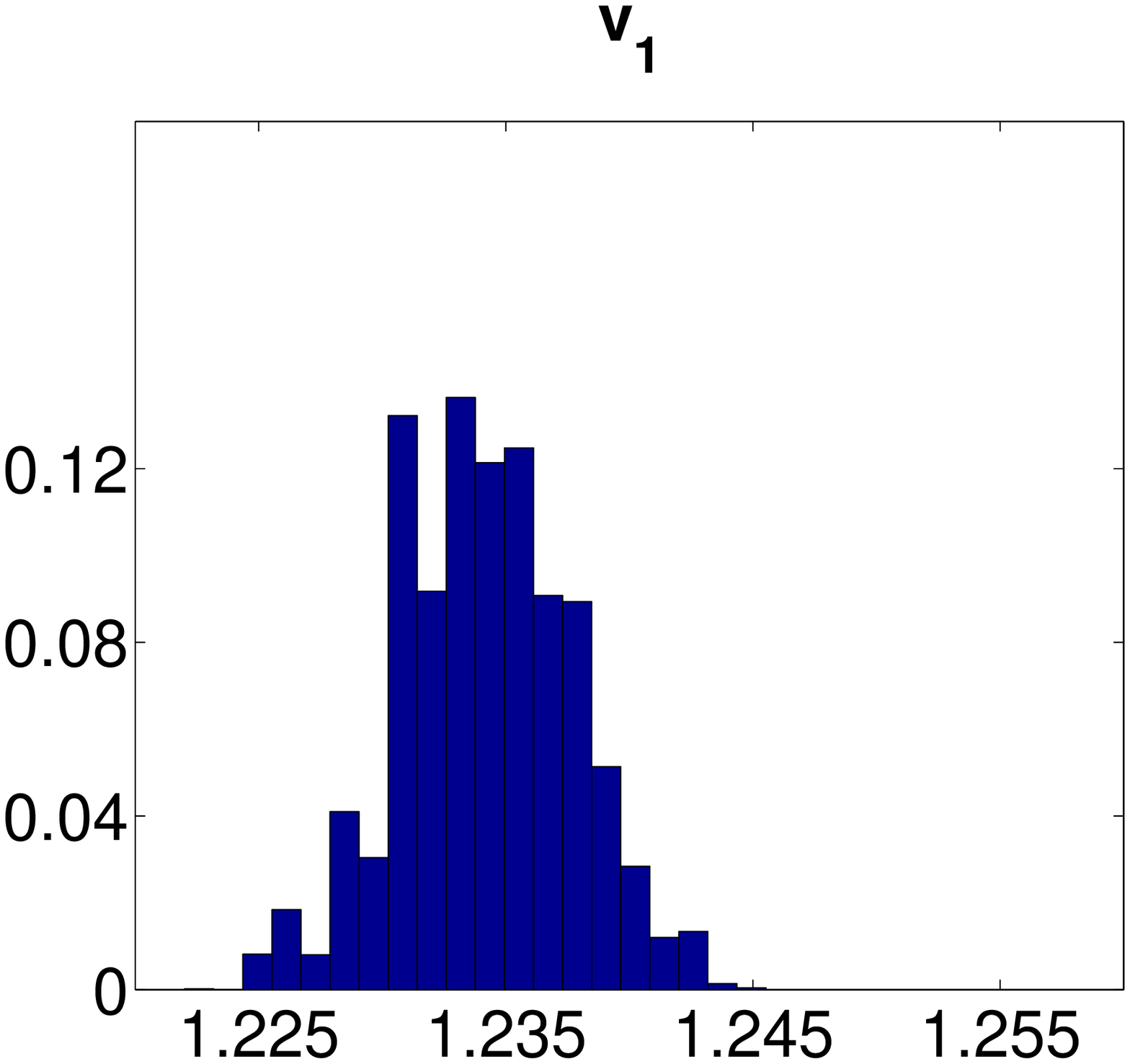} \\
\includegraphics[width=1\textwidth,height=1\textwidth]{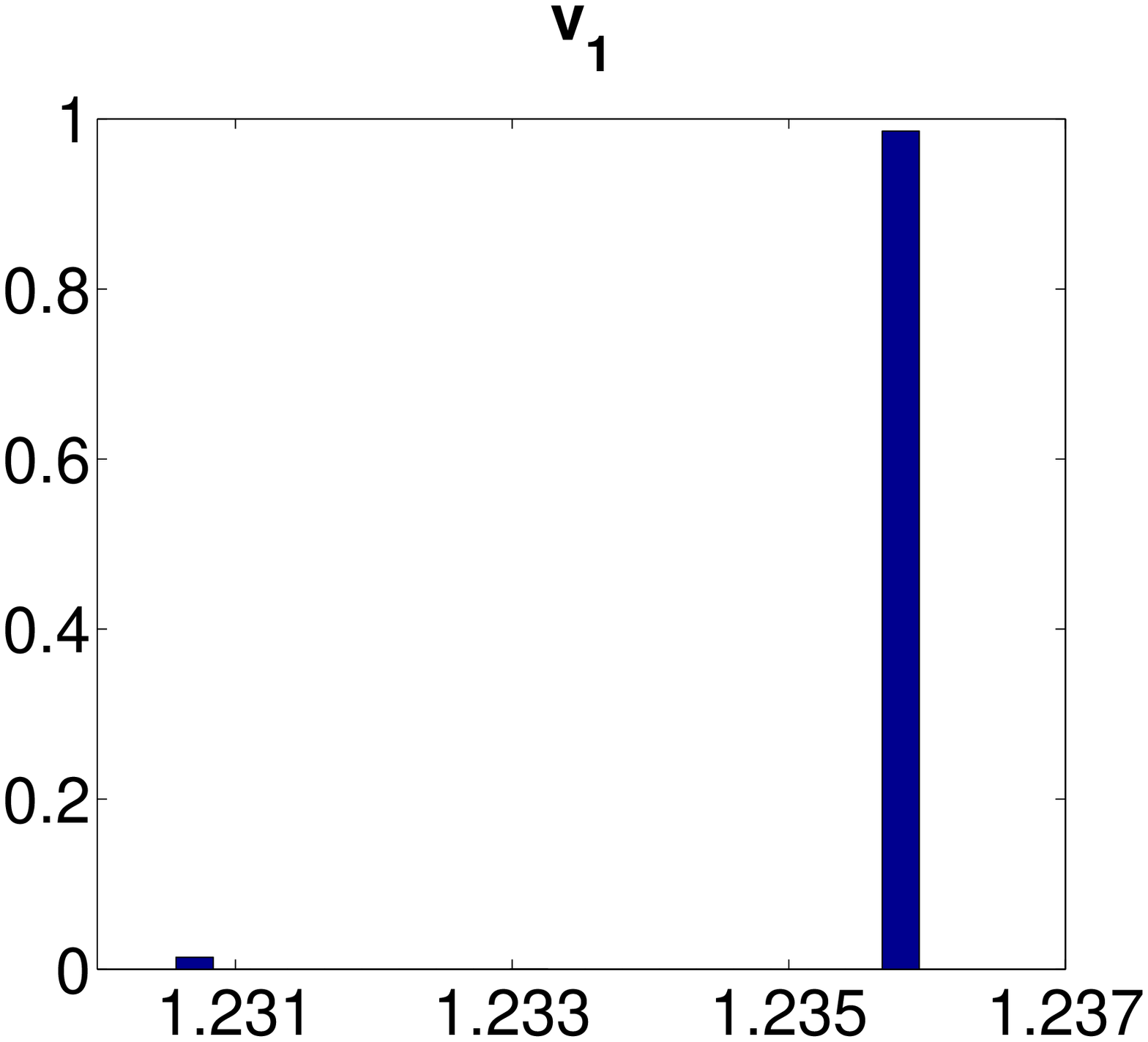}
\end{minipage}
}
\subfigure[]{
\begin{minipage}[b]{0.17\textwidth}
\includegraphics[width=1\textwidth,height=1\textwidth]{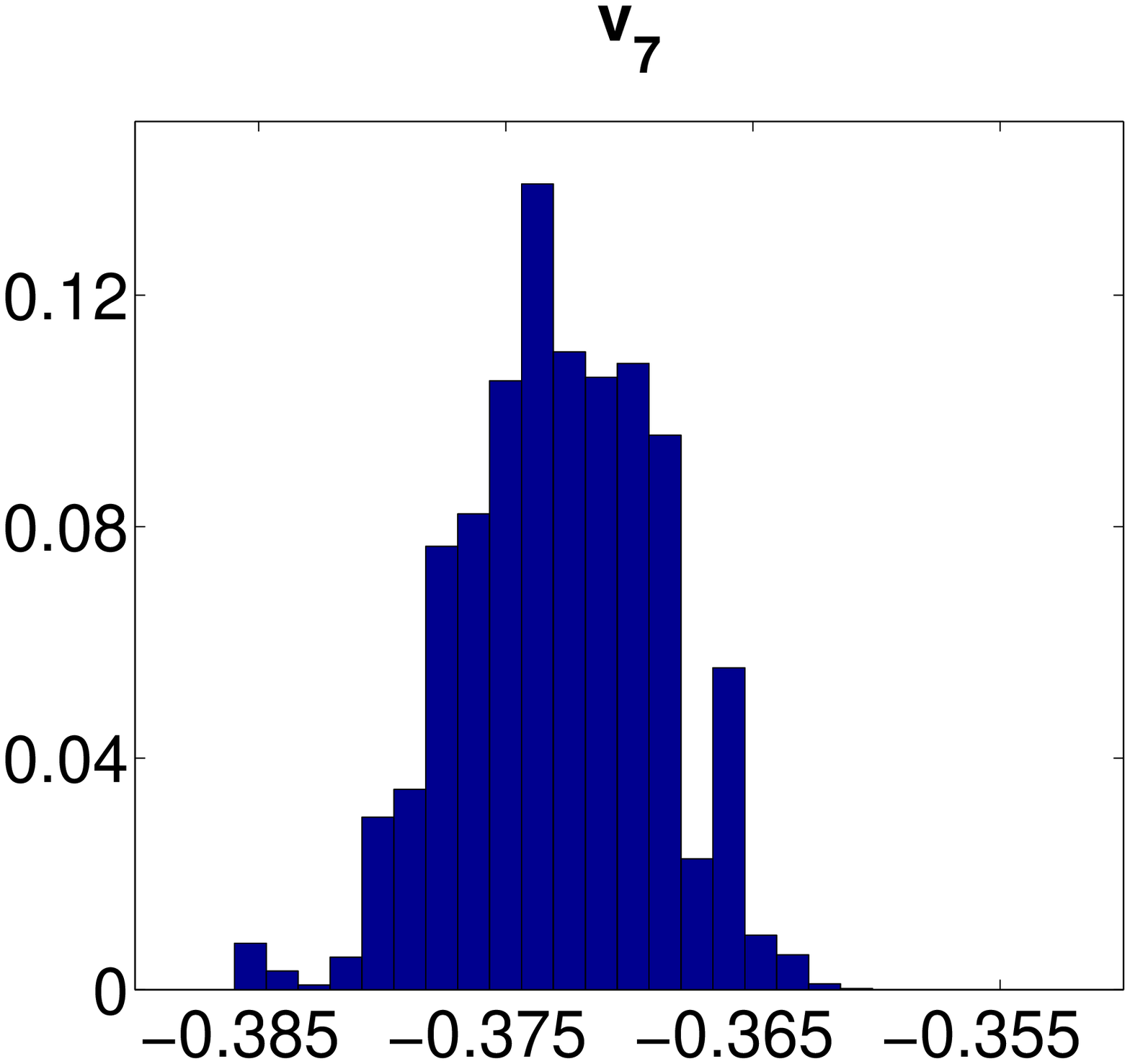} \\
\includegraphics[width=1\textwidth,height=1\textwidth]{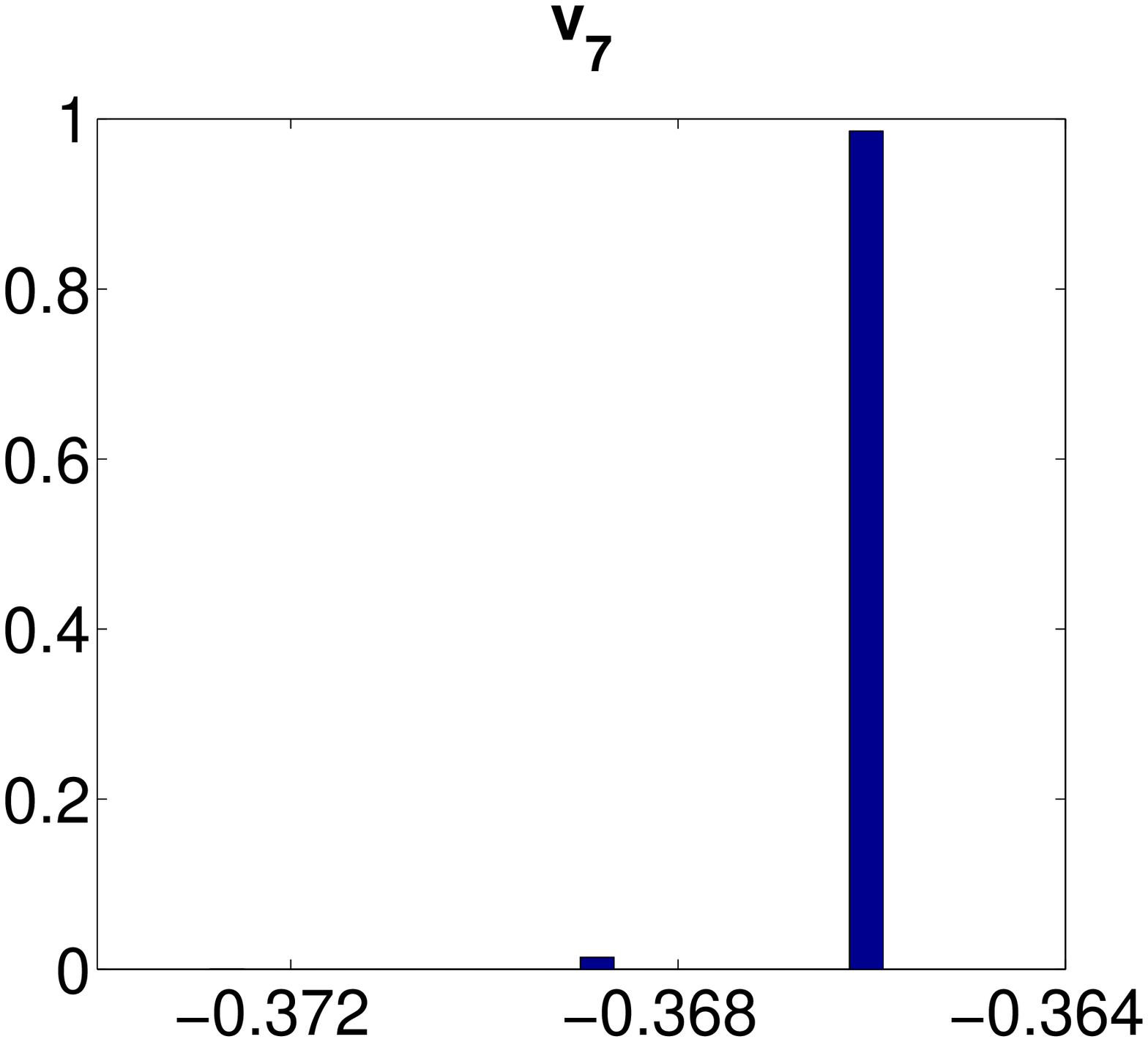}
\end{minipage}
}\subfigure[]{
\begin{minipage}[b]{0.17\textwidth}
\includegraphics[width=1\textwidth,height=1\textwidth]{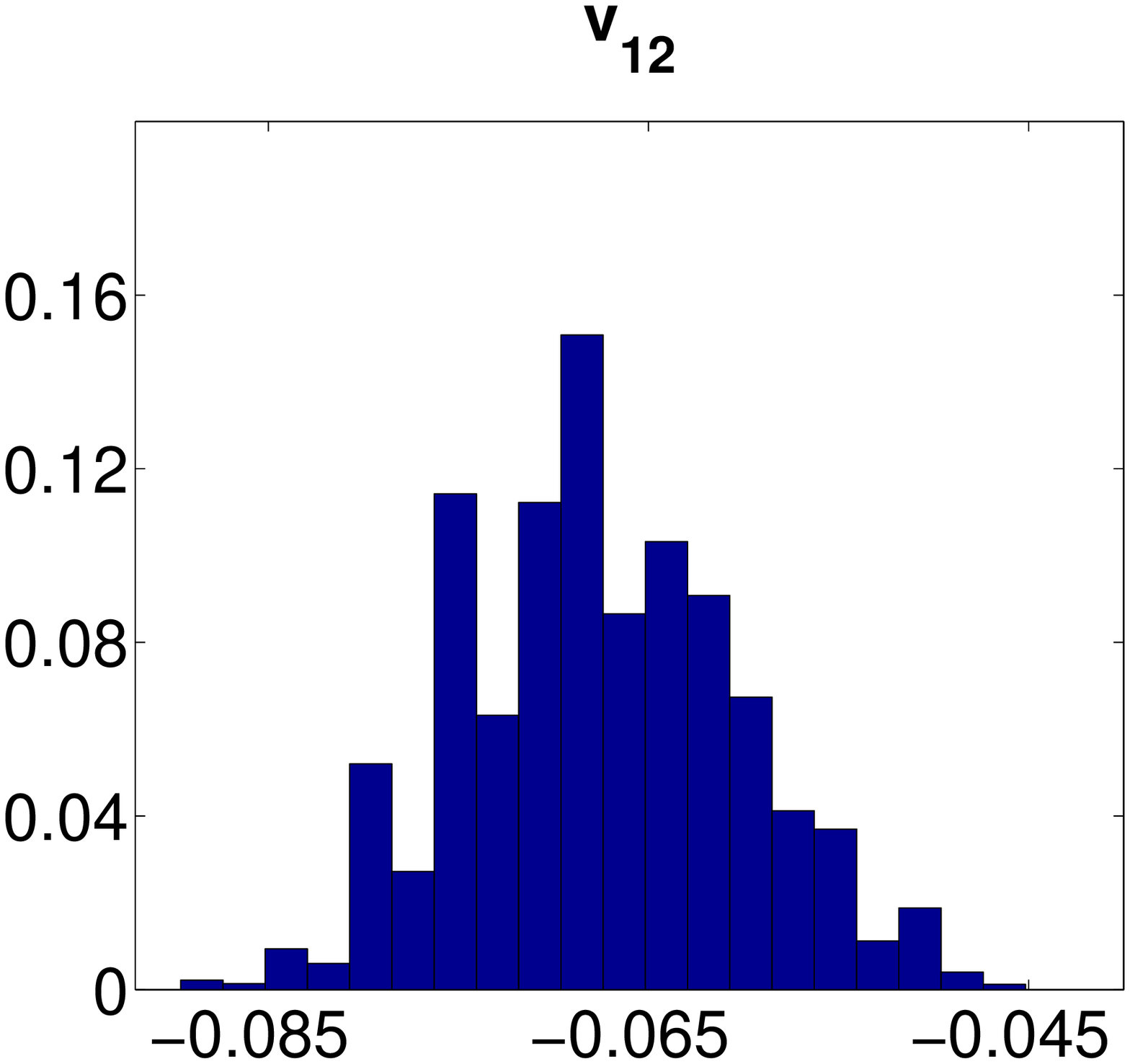} \\
\includegraphics[width=1\textwidth,height=1\textwidth]{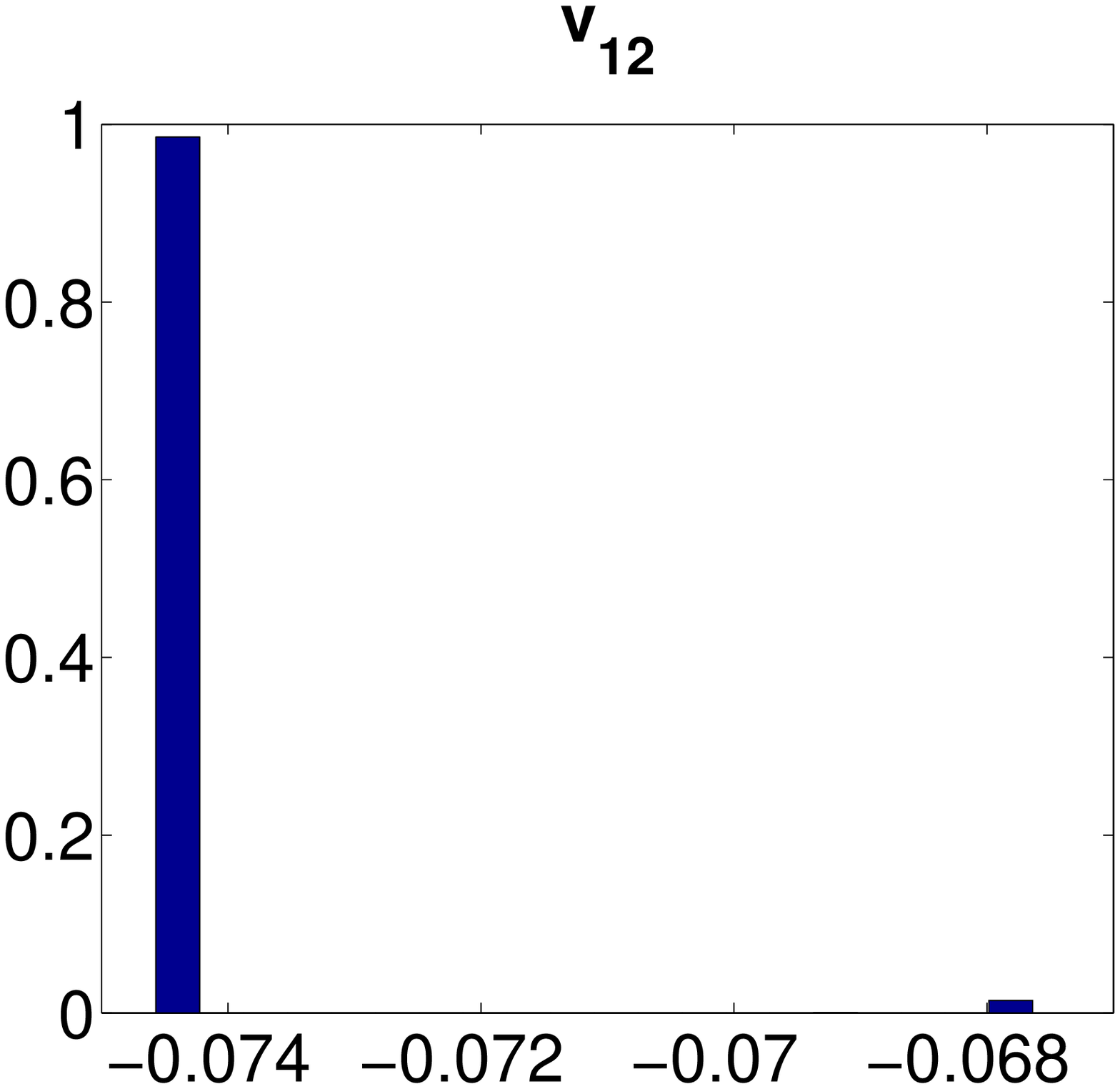}
\end{minipage}
}
\caption{The histogram of marginal posterior by improved implicit sampling  (the first row) and conventional implicit sampling  (the second row) for $[\alpha_{1}, \alpha_{2}, v_{1}, v_{7}, v_{12}]$.}
\label{hist-2}
\end{figure}

Next we  identify the multi-term time fractional derivative  orders and the diffusion field  simultaneously. For this simulation, we reset the boundary condition $g=1+x_{1}$ and the source term $f=2x_{1}+2$. The positive constants and the reference multi-term time fractional orders are chosen the same as in Section \ref{Ex_1}. In addition, the reference diffusion field is the same as Figure \ref{ref} (a). The mesh partition of the forward model is the same as before.
Thus we need to recover the unknown parameter $\textbf{v}=[\alpha_{1},\alpha_{2},v_{1}, \cdots, v_{12}]\in \mathbb{R}^{14}$ in this case. We first implement the augmented-Tikhonov algorithm to obtain the MAP estimate.
Then we use conventional implicit sampling and improved implicit sampling to make comparison.   We first  compute the weight of each sample generating from the conventional implicit sampling  method.  In this method,  the excessive concentration of weights occurs. That is to say, one weight of the particle is close to $1$ and other particles have no influence since their weights are too small, which will result in a poor representation of the posterior.
In order to avoid filter collapse, we adopt the formulation in (\ref{weight-relax1}) to improve the weights.
We note that the weights in the improved implicit sampling obey the formula  (\ref{weight-relax1}) and the weights in the conventional  implicit sampling obey the formula (\ref{weight-propotion}).
To show the distribution of the weights for the two implicit sampling methods,  we count  the number of samples when associated weights lie in the disjoint  intervals and list the results  in Table \ref{table-weights-2}.
 Here we set the scale parameter  $\vartheta=15$ for the improved implicit sampling.
From this table, it can be found that most of the samples have extremely small weights in the conventional implicit sampling, which will cause the ensemble collapse.   Thus, the conventional implicit sampling  may
 underestimate the uncertainty of the unknown parameters. But for the improved implicit sampling, we can see that
the distribution of weights is relatively flat while it retain the original order. This is a significant improvement.
In order to clearly illustrate the effect of the two implicit sampling  on the posterior distribution,   we plot the marginal posterior histogram of $[\alpha_{1}, \alpha_{2}, v_{1}, v_{7}, v_{12}]$ in Figure \ref{hist-2}.
It can be seen that samples are over-concentrate by the conventional implicit sampling  and the improved implicit sampling effectively alleviates the  situation.  This  is consistent with the results in  Table \ref{table-weights-2}.

\renewcommand\arraystretch{1.5}
\begin{table}
\centering
\resizebox{\textwidth}{!}{
   \begin{tabular}{cccccccc}
       \toprule
       \multirow{2}{*}{Method}
        & \multicolumn{7}{c}{Interval}\\
       \cmidrule{2-8}
       & [0, $10^{-6}$) & [$10^{-6}$, $10^{-5}$) & [$10^{-5}$,$10^{-4}$) & [$10^{-4}$,$10^{-3}$) & [$10^{-3}$,$10^{-2}$) & [$10^{-2}$,$10^{-1}$) & [$10^{-1}$,1]\\
       $\text{\uppercase\expandafter{I-IS}}$&3004&423&593&750&219&11&0\\
       $\text{\uppercase\expandafter{C-IS}}$&4993&3&1&1&0&1&1\\
       \bottomrule
  \end{tabular}
  }
  \caption{The  weights distribution  for improved implicit sampling (I-IS) and conventional implicit sampling (C-IS)}\label{table-weights-2}
\end{table}

%

\begin{figure}
\centering
\subfigure[]{
\includegraphics[width=0.3\textwidth, height=0.3\textwidth]{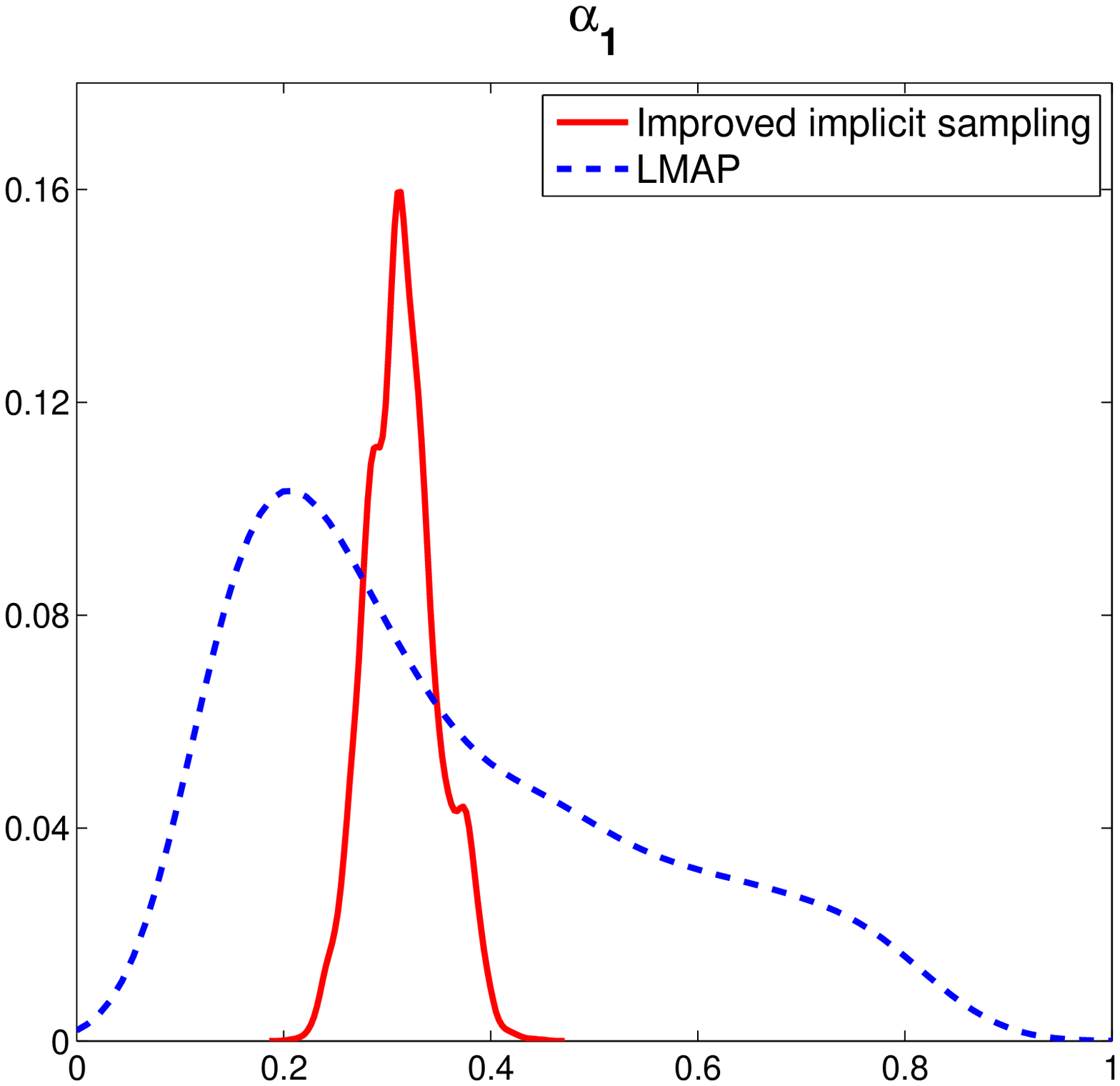}
}
\subfigure[]{
\includegraphics[width=0.3\textwidth, height=0.3\textwidth]{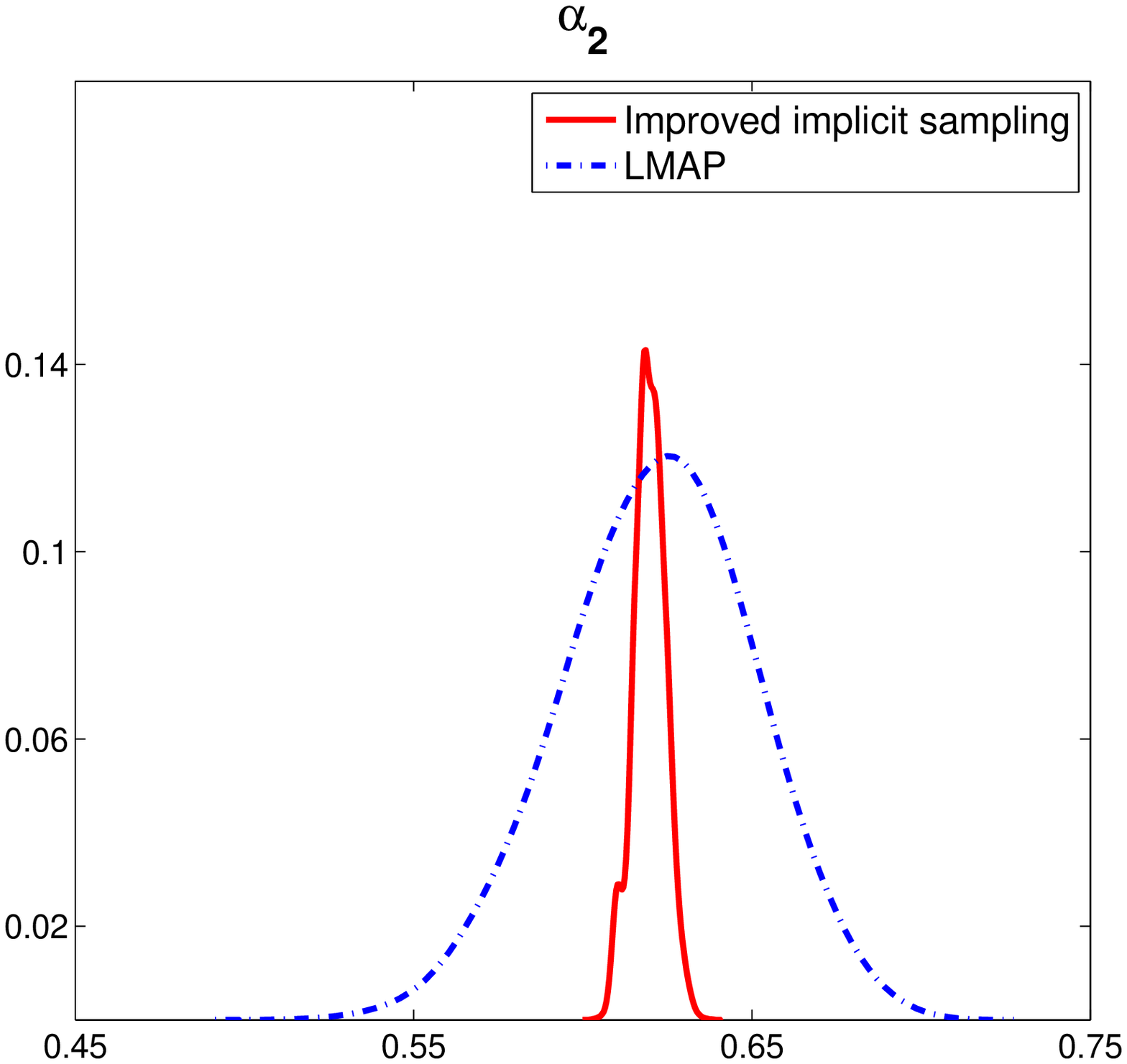}
}
\caption{The posterior marginal distribution for $\alpha_{1}$ (left) and $\alpha_{2}$ (right) by improved implicit sampling and LMAP.}\label{imp-lmap-2}
\end{figure}

\begin{table}
  \centering
\begin{tabular}{|c|c|c|}
\hline
 & Skewness  & Kurtosis\\
\hline
Improved  implicit sampling & (0.2841, -0.0662) & (-0.0392, 0.0230) \\
\hline
 LMAP              & (0.7403, -0.1804) & (-0.5048, -0.1120) \\
\hline
\end{tabular}
 \caption{\label{Tab-skewness-kurtosis}The skewness and kurtosis of $(\alpha_{1}, \alpha_{2})$  by improved implicit sampling and LMAP.}
\end{table}

\begin{figure}
\centering
\includegraphics[width=0.75\textwidth, height=0.6\textwidth]{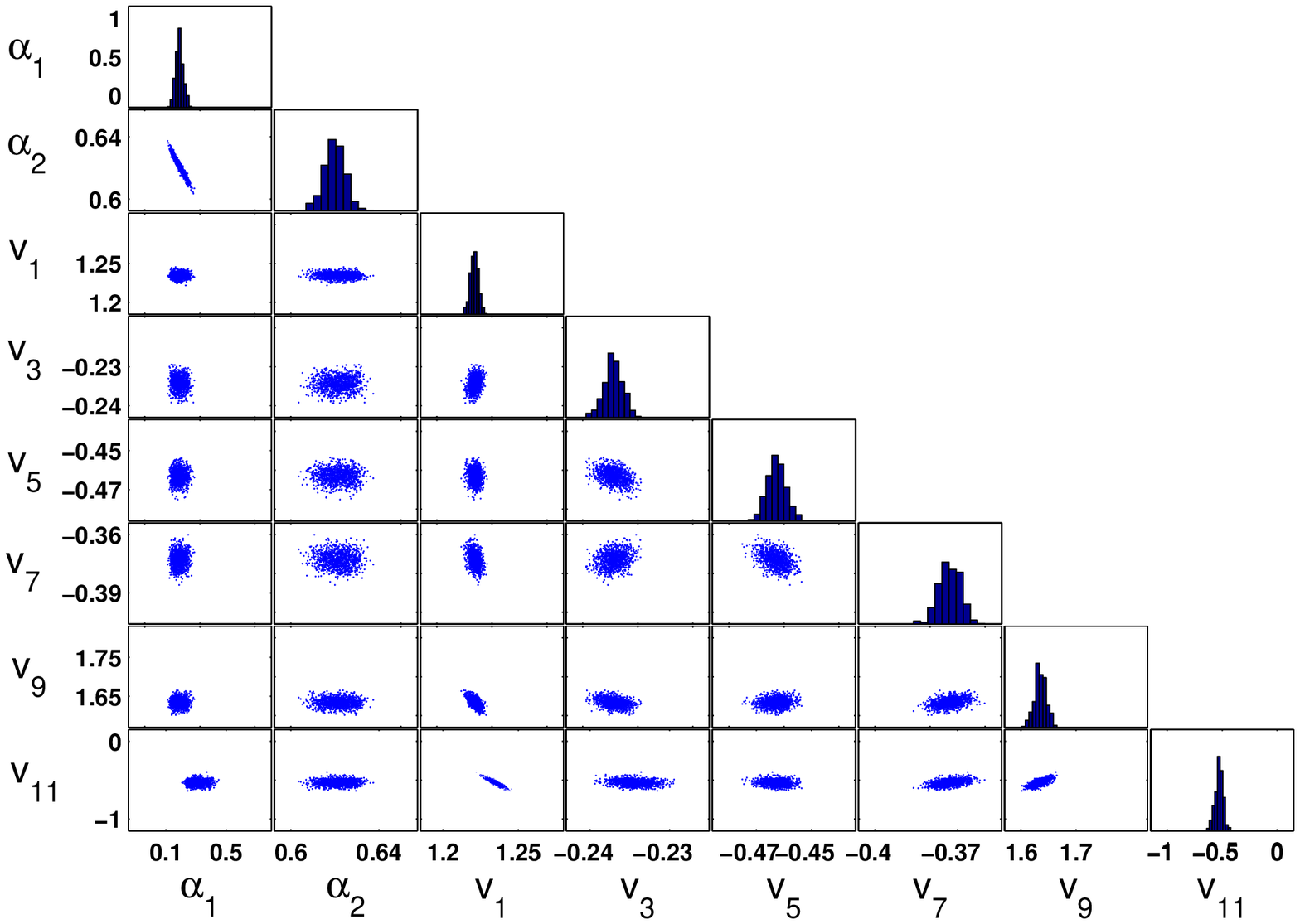}
\caption{The matrix plot for $[\alpha_{1},\alpha_{2},v_{1},v_{3},\cdots,v_{11}]$ by improved implicit sampling.}\label{matrixplot-2}
\end{figure}

Finally we make a comparison between improved implicit sampling and LMAP.
Figure \ref{imp-lmap-2} plots the marginal posterior for $\alpha_{1}$ and $\alpha_{2}$ using these two approaches, which
illustrates that samples  by improved implicit sampling are more concentrated on the region near MAP point than LMAP. To characterize the marginal posterior distribution of $(\alpha_{1}, \alpha_{2})$, we list the skewness and kurtosis of the two parameters based on samples  by the improved implicit sampling and LMAP method in Table \ref{Tab-skewness-kurtosis}. We find  that the skewness of $\alpha_{1}$ produced by LMAP method is positive and it implies that  the heavy tail is on the right, which is consistent with Figure \ref{imp-lmap-2} (a).
To describe  the correlation among these parameters, the matrix plot for $[\alpha_{1},\alpha_{2},v_{1},v_{3},\cdots,v_{11}]$ is shown  in Figure \ref{matrixplot-2}. It is observed from the pairwise joint sample plots in Figure \ref{matrixplot-2} that the multi-term fractional orders $\alpha_{1}$ and $\alpha_{2}$ are clearly negatively correlated, which is  consistent with section \ref{Ex_1}. We can also conclude that the fractional orders are almost independent of all the parameters in the diffusion  field.  The parameters in the diffusion  field are almost  independent except for a small correlation between $v_{1}$ and $v_{11}$.


\subsection{Inversion for reaction coefficient}\label{Ex_3}
\begin{figure}
\centering
\subfigure[]{
\includegraphics[width=0.35\textwidth, height=0.3\textwidth]{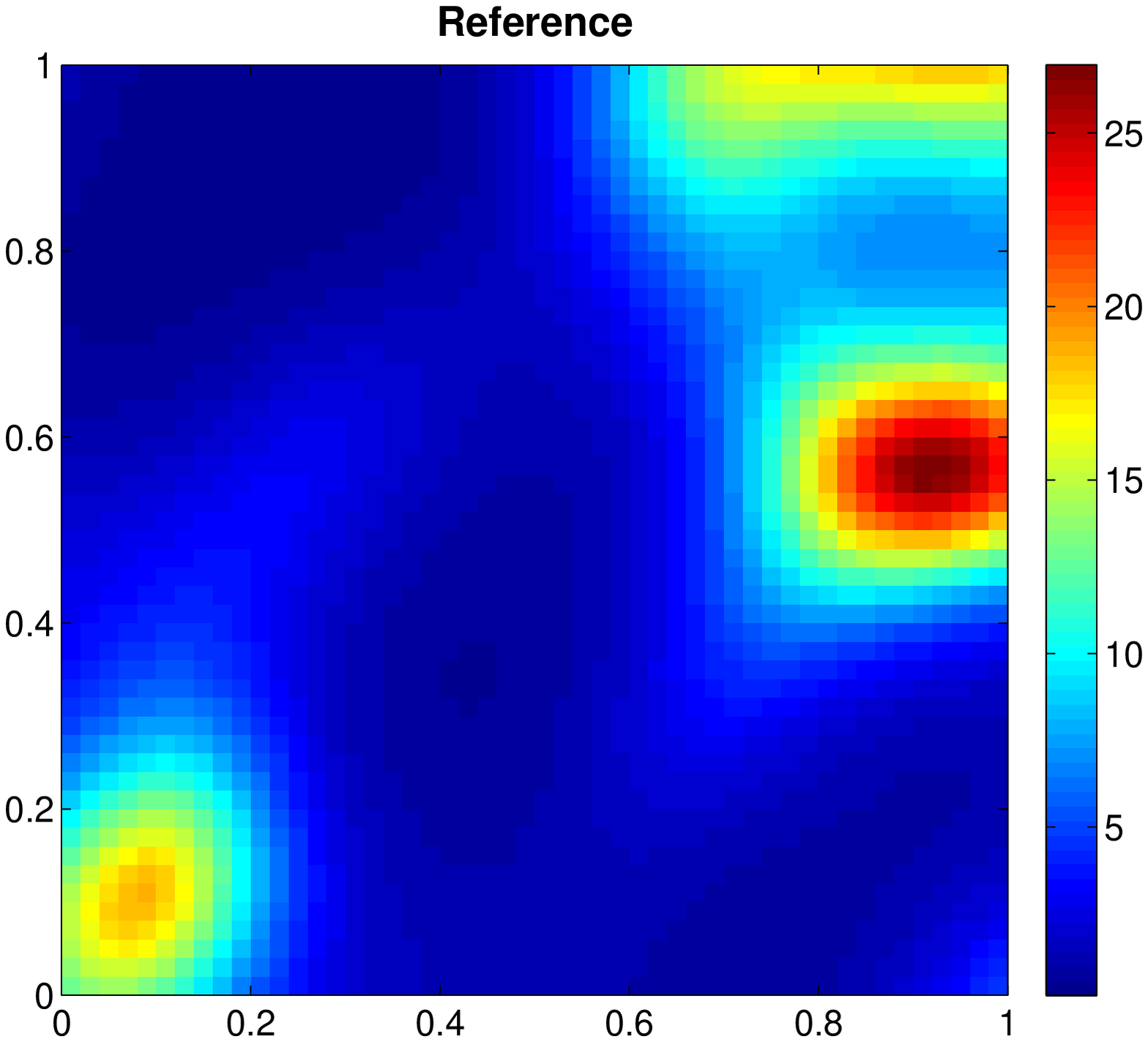}
}
\subfigure[]{
\includegraphics[width=0.35\textwidth, height=0.3\textwidth]{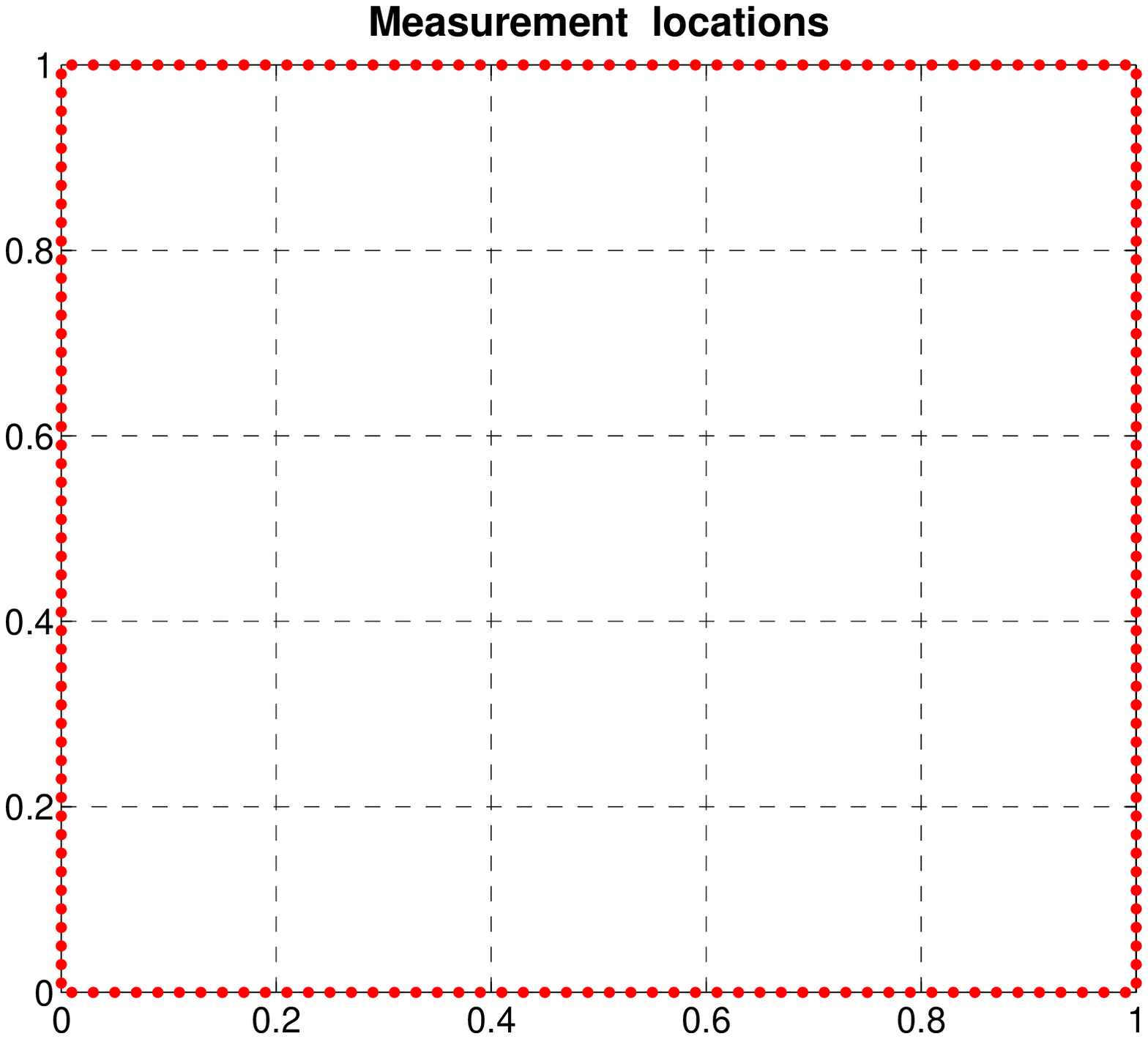}
}
\caption{The truth  profile of $q(x)$ and the measurement distribution  for Subsection \ref{Ex_3}.}\label{ref_3}
\end{figure}

In this subsection, we aim at  recovering the reaction coefficient $q(x)$. The truth of  $q(x)$ is depicted  in Figure \ref{ref_3} (a). As illustrated in Subsection \ref{Ex_2}, we assume  that the unknown reaction  field $q(x,\varrho)$ can be parameterized  by Karhunen-Lo$\grave{e}$ve expansion. For the covariance function of  $\log q(x,\varrho)$, we set $l_{1}=l_{2}=0.03$, $\rho=1$ and $\kappa=1$. Here, the first 34 terms are truncated
for the parametrization.  The forward model is defined  on $50\times 50$ uniform fine grid, and mixed GMsFEM model  is implemented on $5\times 5$ coarse  grid. The number of local multiscale basis functions is set to $5$. We choose the source term $f=10$ and the Dirichlet boundary condition $g=1-x_{1}$. The measurements are taken at $t=0.6$ from the whole boundary Neumann data, whose distribution is shown in Figure \ref{ref_3} (b).

\begin{figure}
\centering
\subfigure[]{
\includegraphics[width=0.22\textwidth, height=0.22\textwidth]{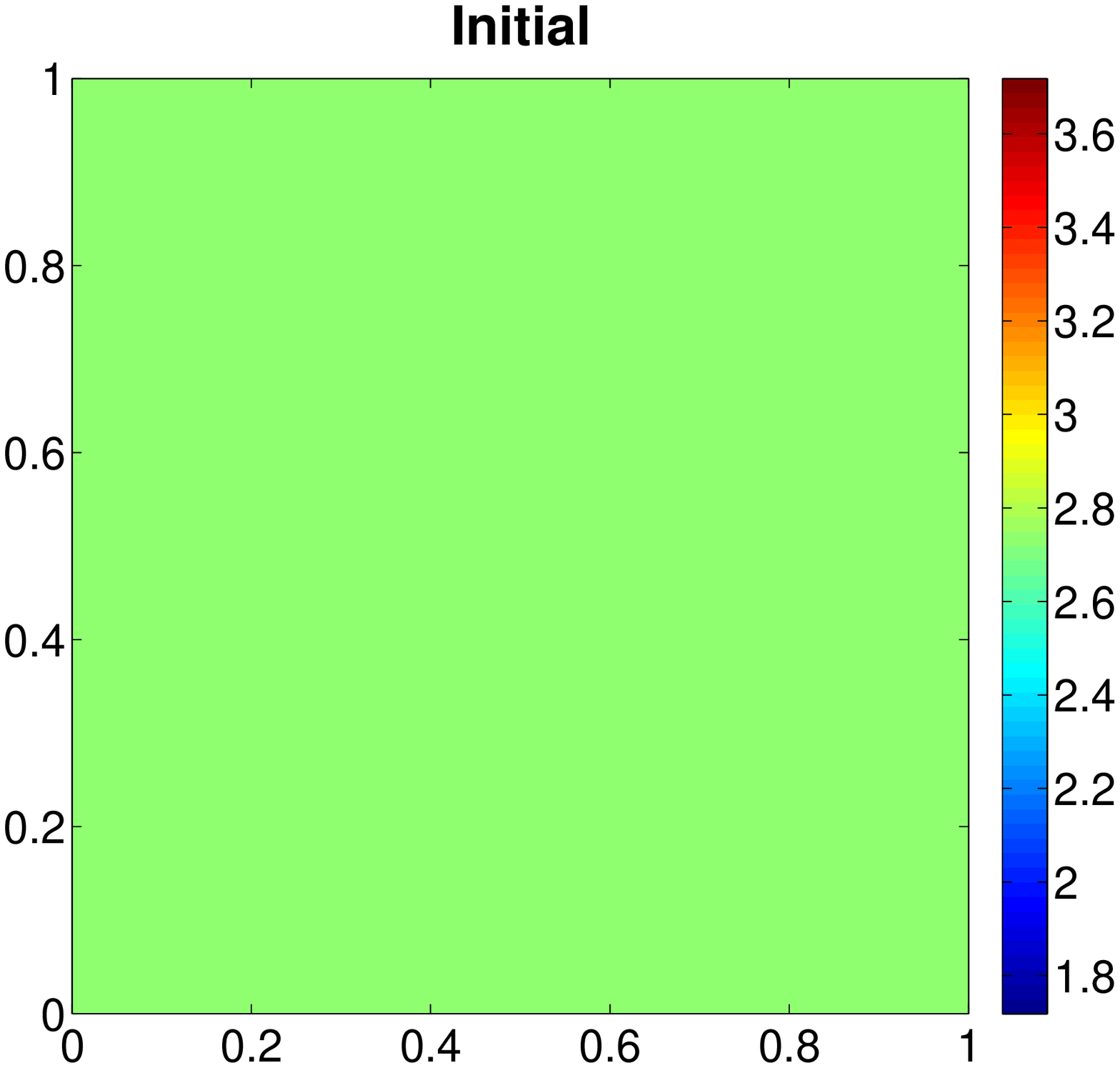}
}
\subfigure[]{
\includegraphics[width=0.22\textwidth, height=0.22\textwidth]{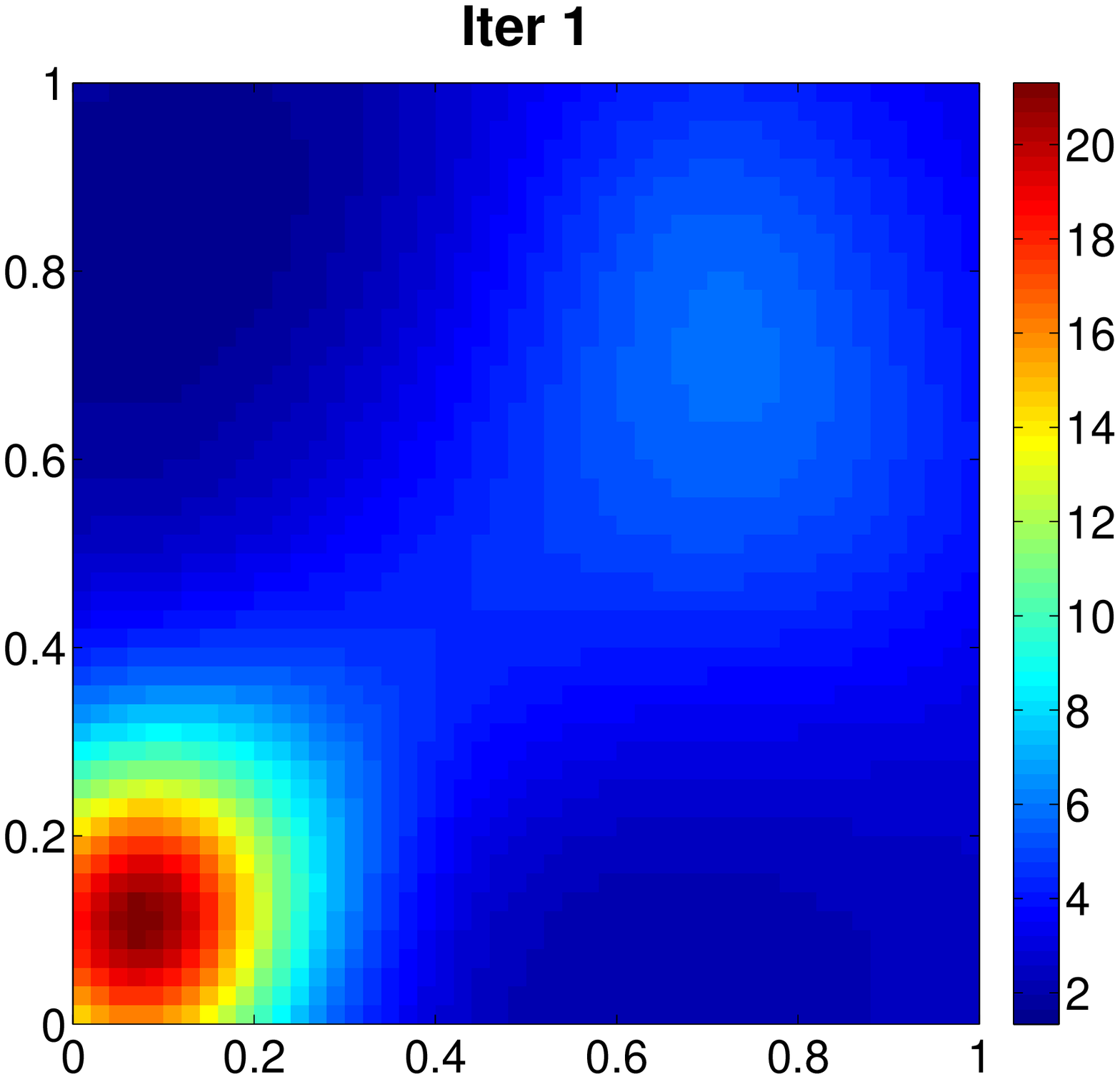}
}
\subfigure[]{
\includegraphics[width=0.22\textwidth, height=0.22\textwidth]{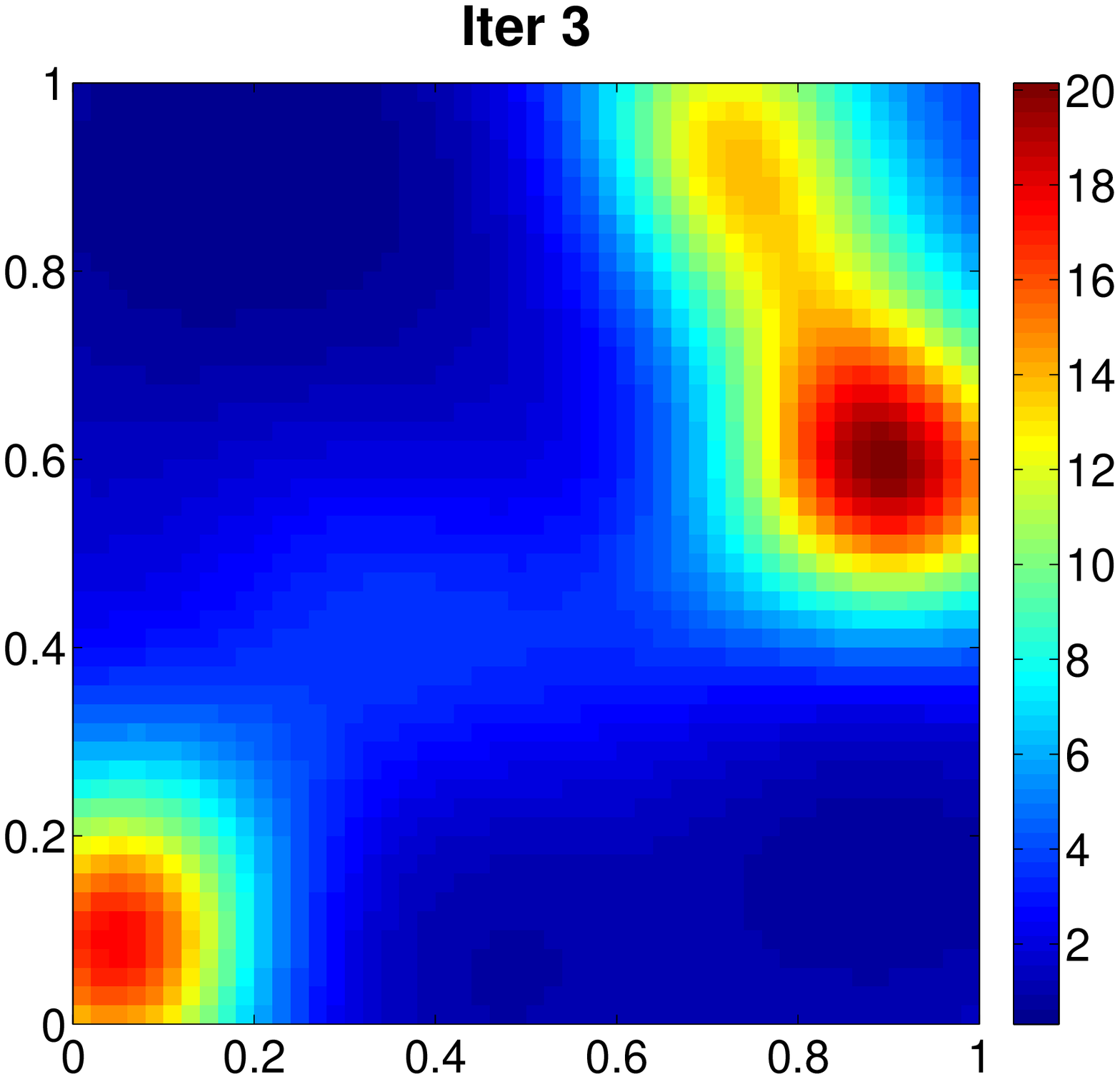}
}
\subfigure[]{
\includegraphics[width=0.22\textwidth, height=0.22\textwidth]{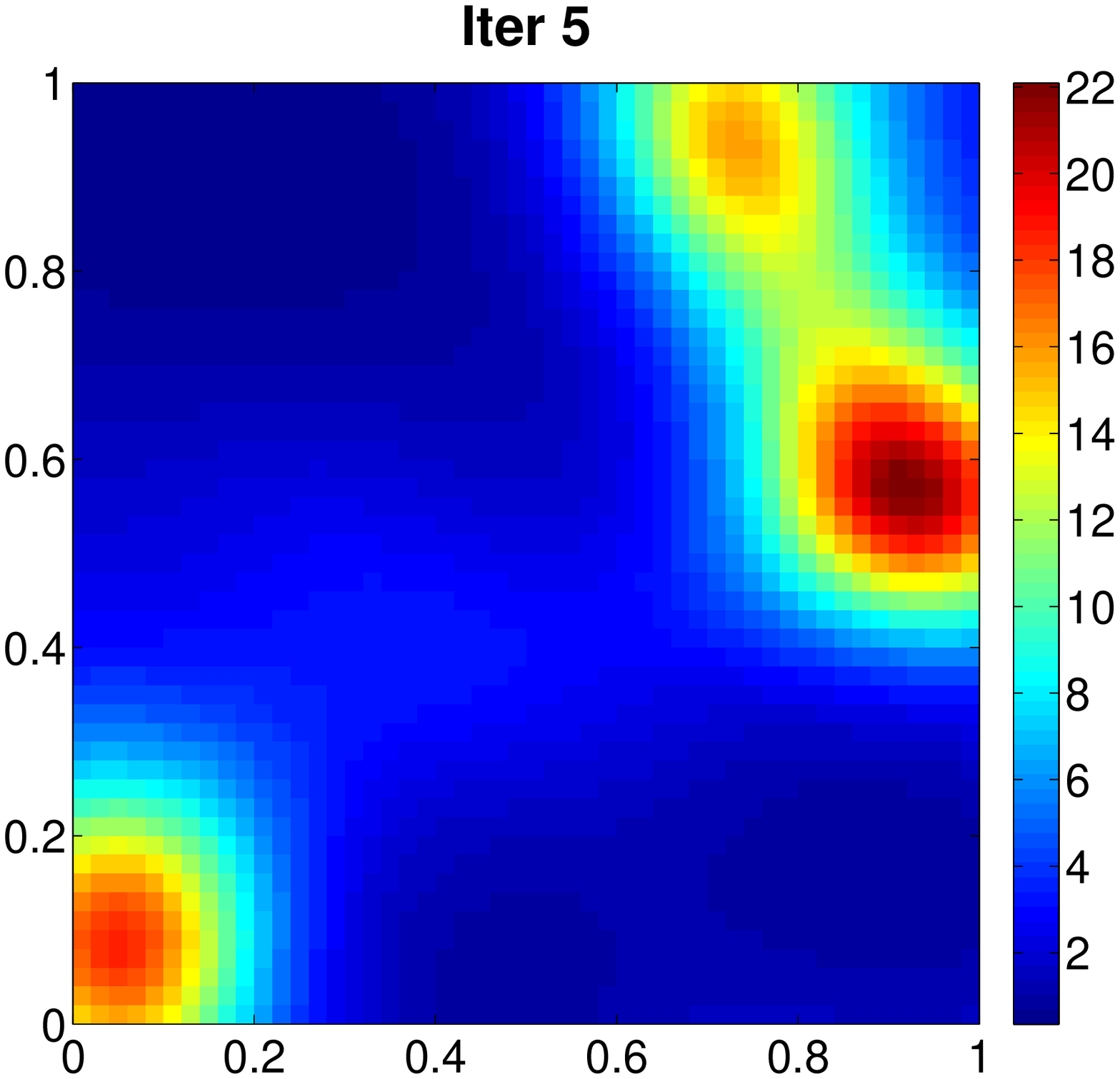}
}
\caption{The initial guess (a) and the reconstruction results (b-d) for $q(x,\varrho)$ in  the $1$-st, $3$-rd and $5$-th  iteration.}\label{L1}
\end{figure}

In this subsection, we attempt to use  sparse prior information to recover the unknown  $q(x,\varrho)$. Assume that the regularization parameter $\mu=0.005$ is fixed. In order to illustrate the performance of the iteratively reweighted approach, we plot  the reconstruction results in Figure \ref{L1}.
 By  Figure \ref{L1}, after five iterations, the inversion solution can  effectively capture the major feature of the truth profile although the measurement data are only collected on the boundary.

\begin{figure}
\centering
\subfigure[]{
\begin{minipage}[b]{0.17\textwidth}
\includegraphics[width=1\textwidth,height=1\textwidth]{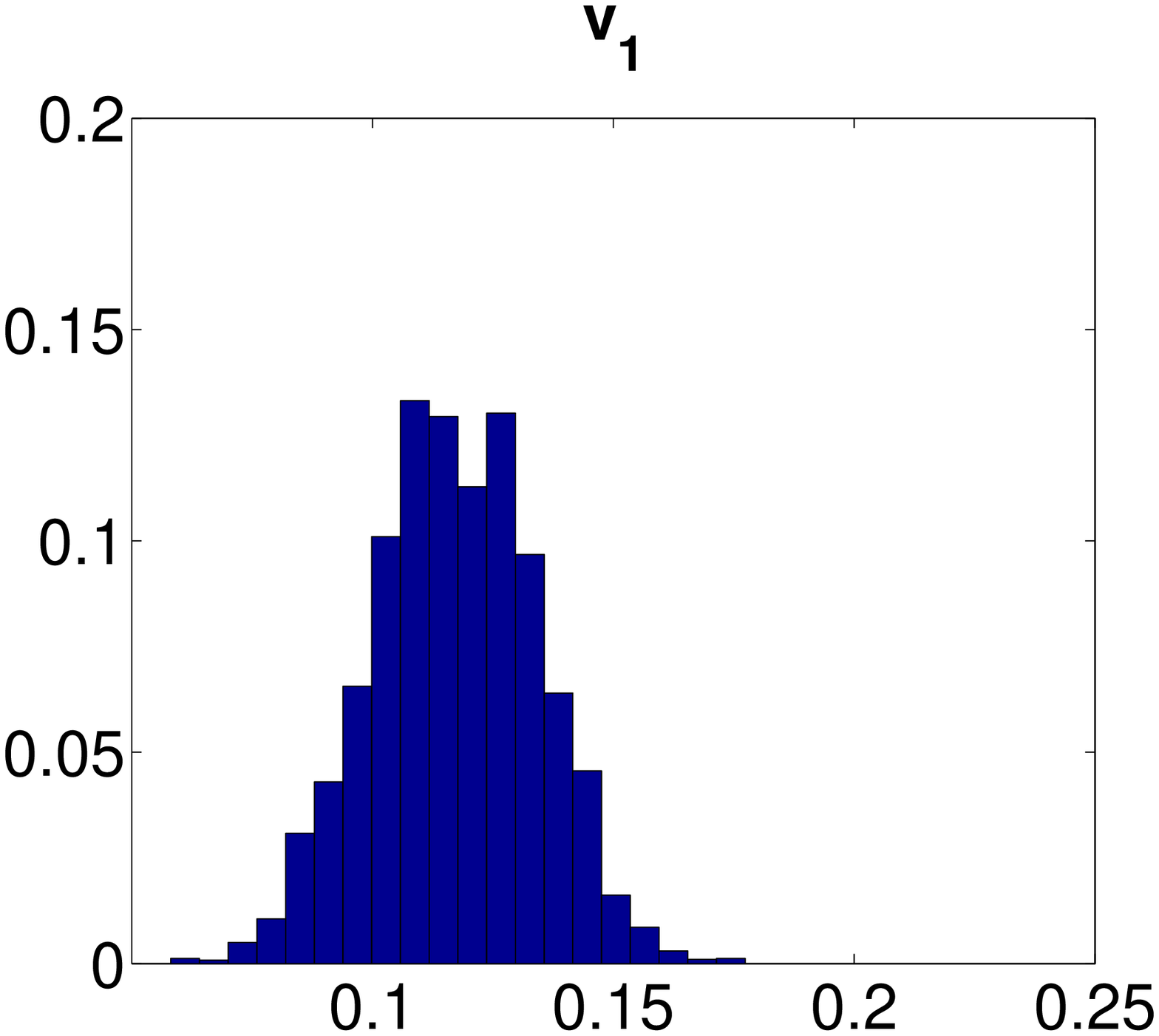} \\
\includegraphics[width=1\textwidth,height=1\textwidth]{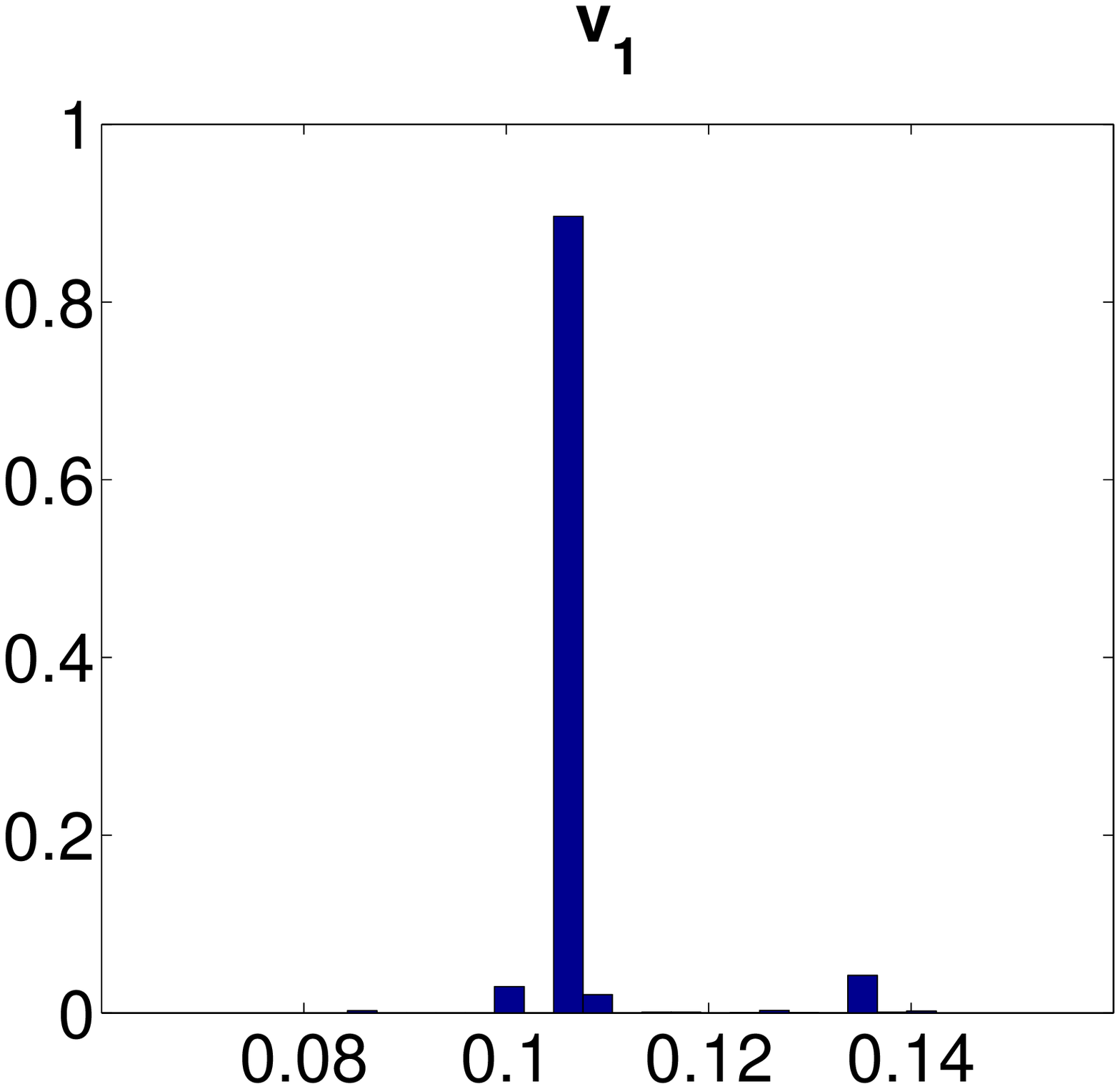}
\end{minipage}
}
\subfigure[]{
\begin{minipage}[b]{0.17\textwidth}
\includegraphics[width=1\textwidth,height=1\textwidth]{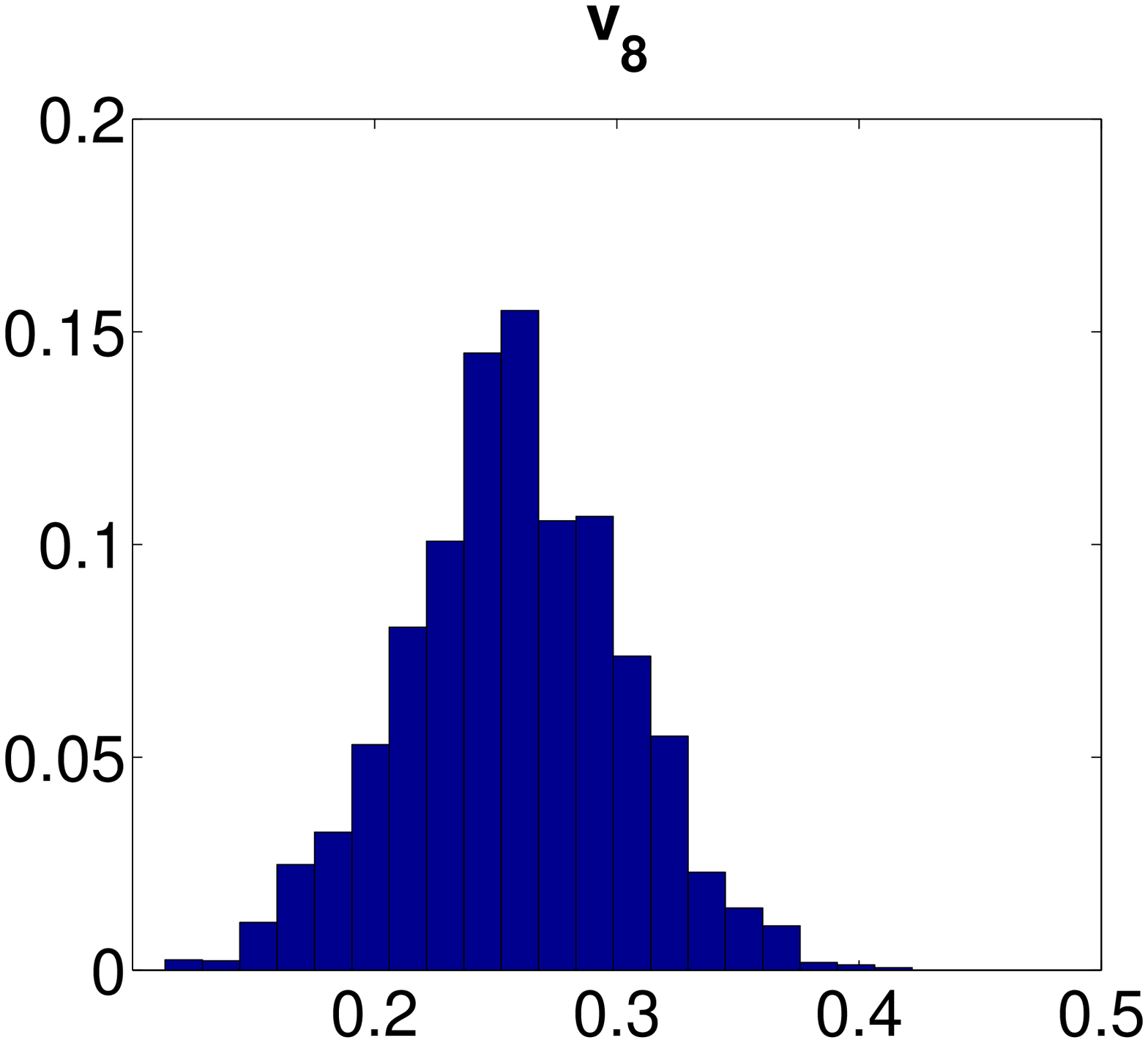} \\
\includegraphics[width=1\textwidth,height=1\textwidth]{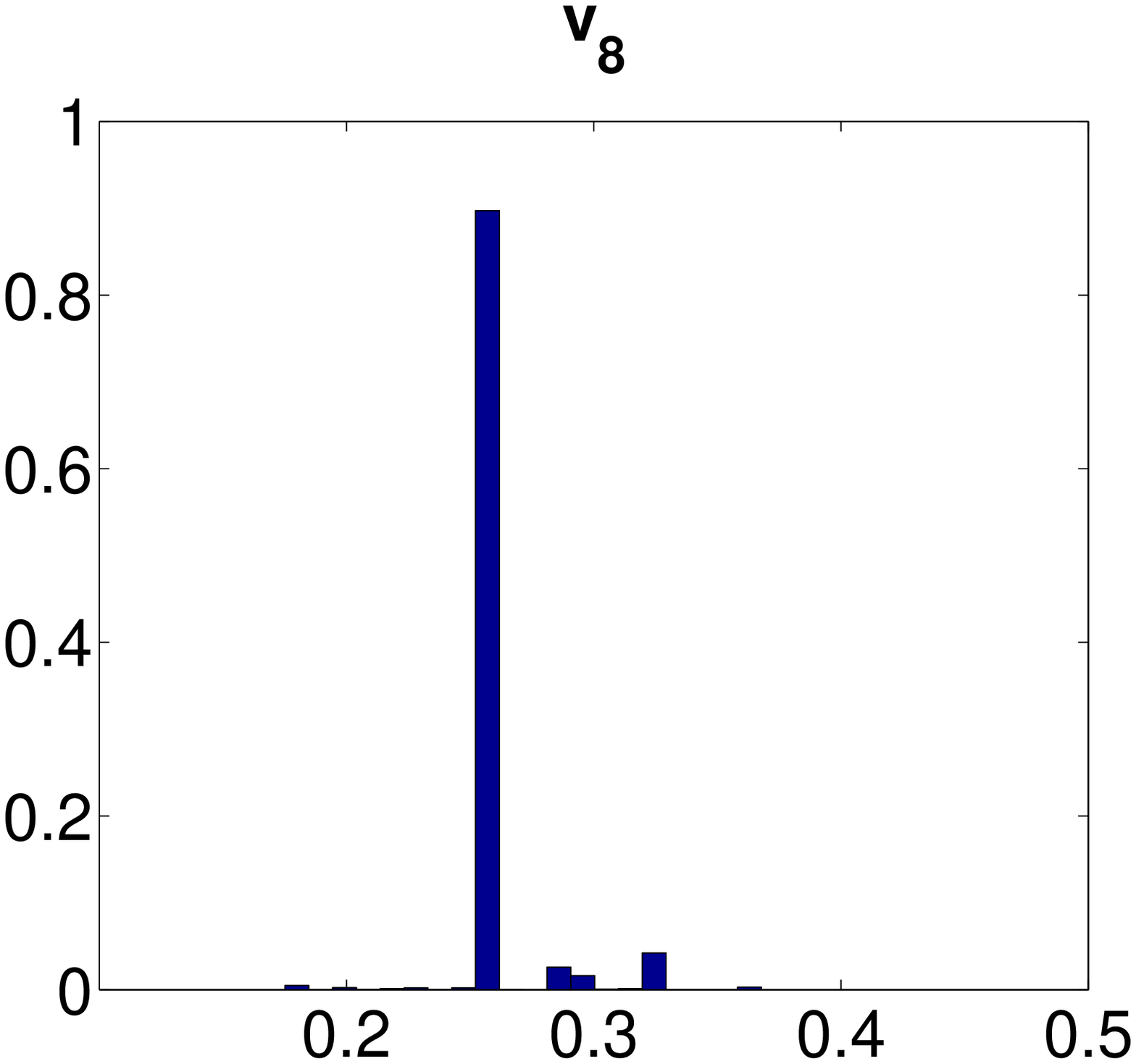}
\end{minipage}
}
\subfigure[]{
\begin{minipage}[b]{0.17\textwidth}
\includegraphics[width=1\textwidth,height=1\textwidth]{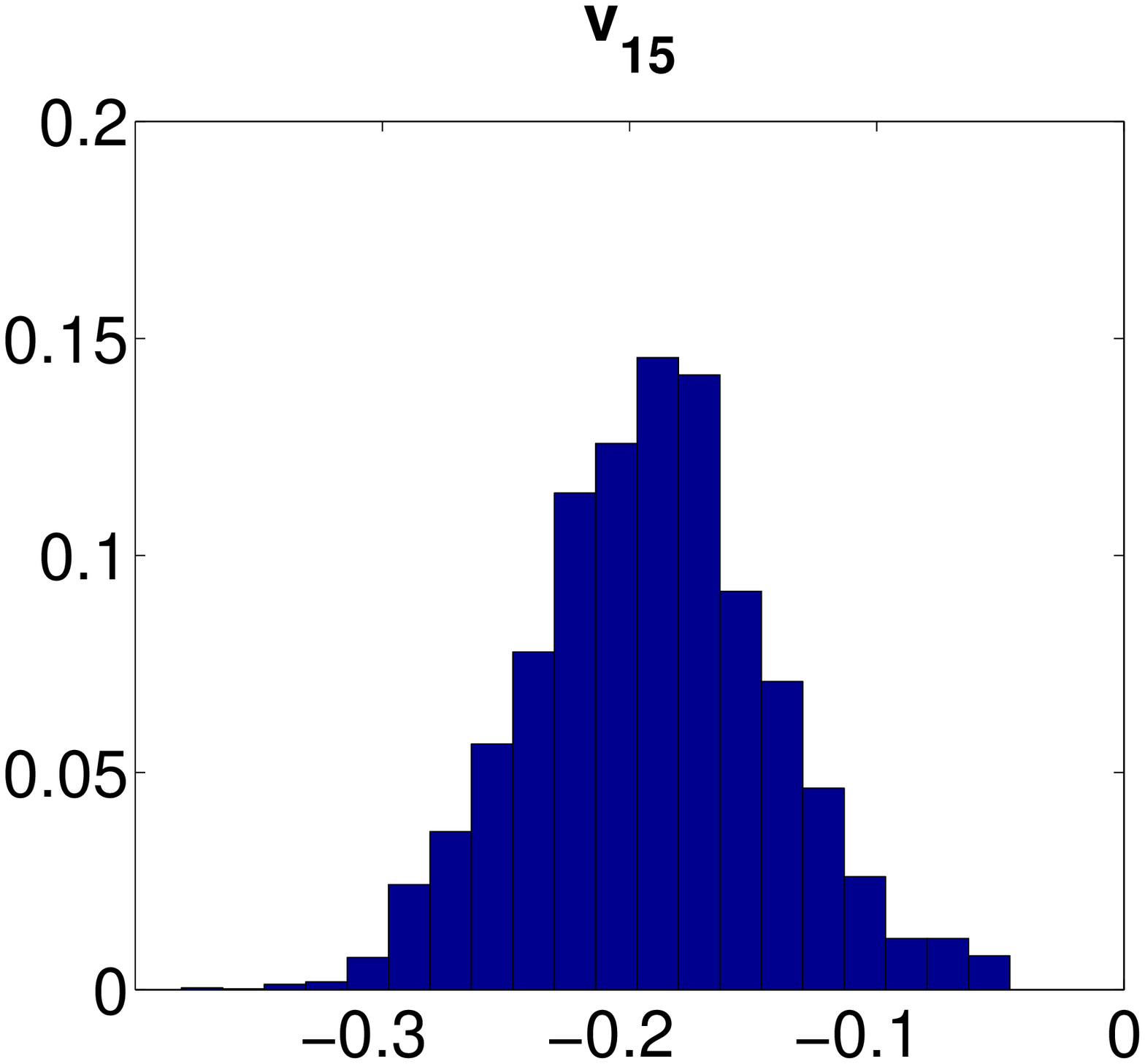} \\
\includegraphics[width=1\textwidth,height=1\textwidth]{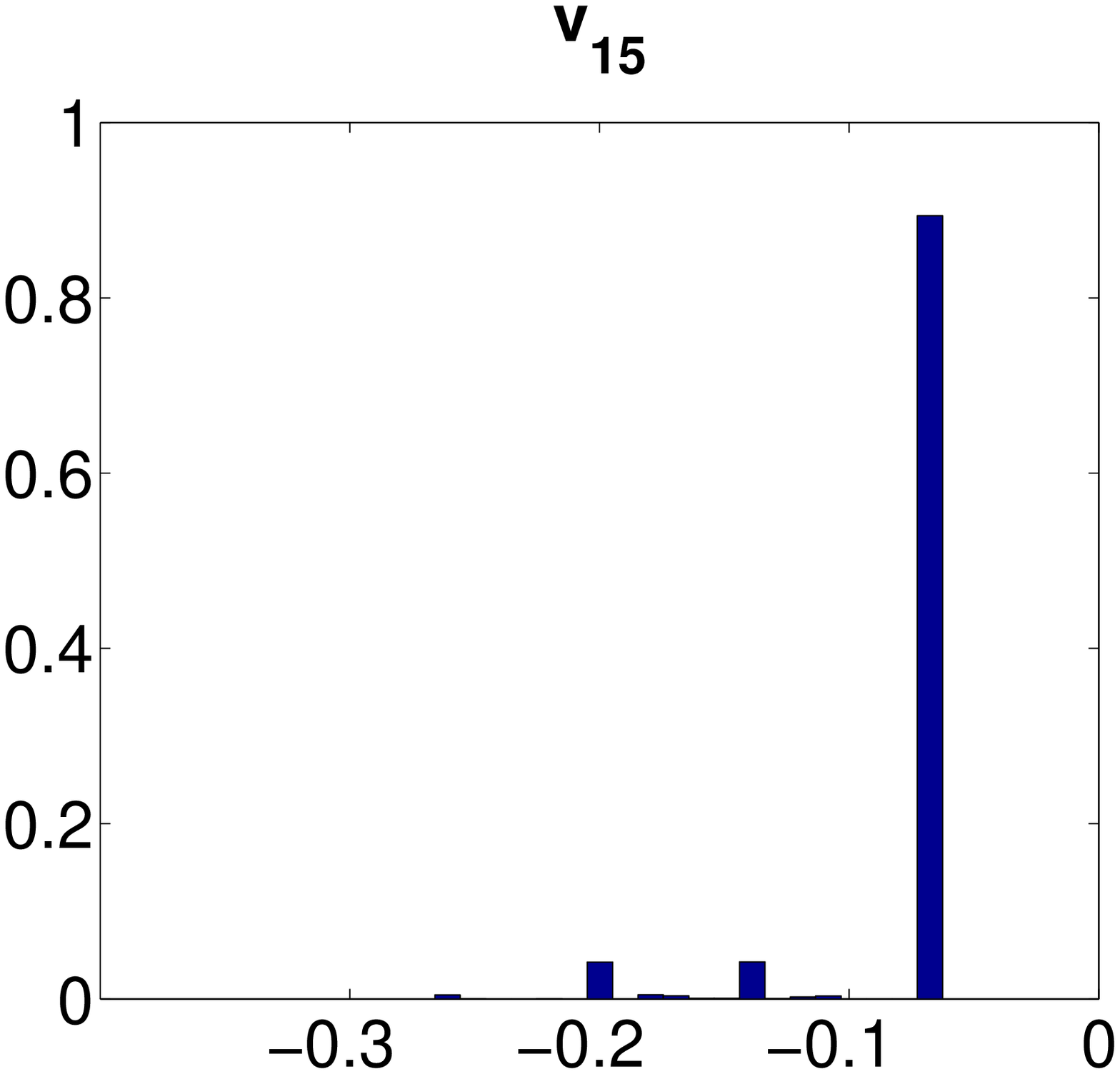}
\end{minipage}
}
\subfigure[]{
\begin{minipage}[b]{0.17\textwidth}
\includegraphics[width=1\textwidth,height=1\textwidth]{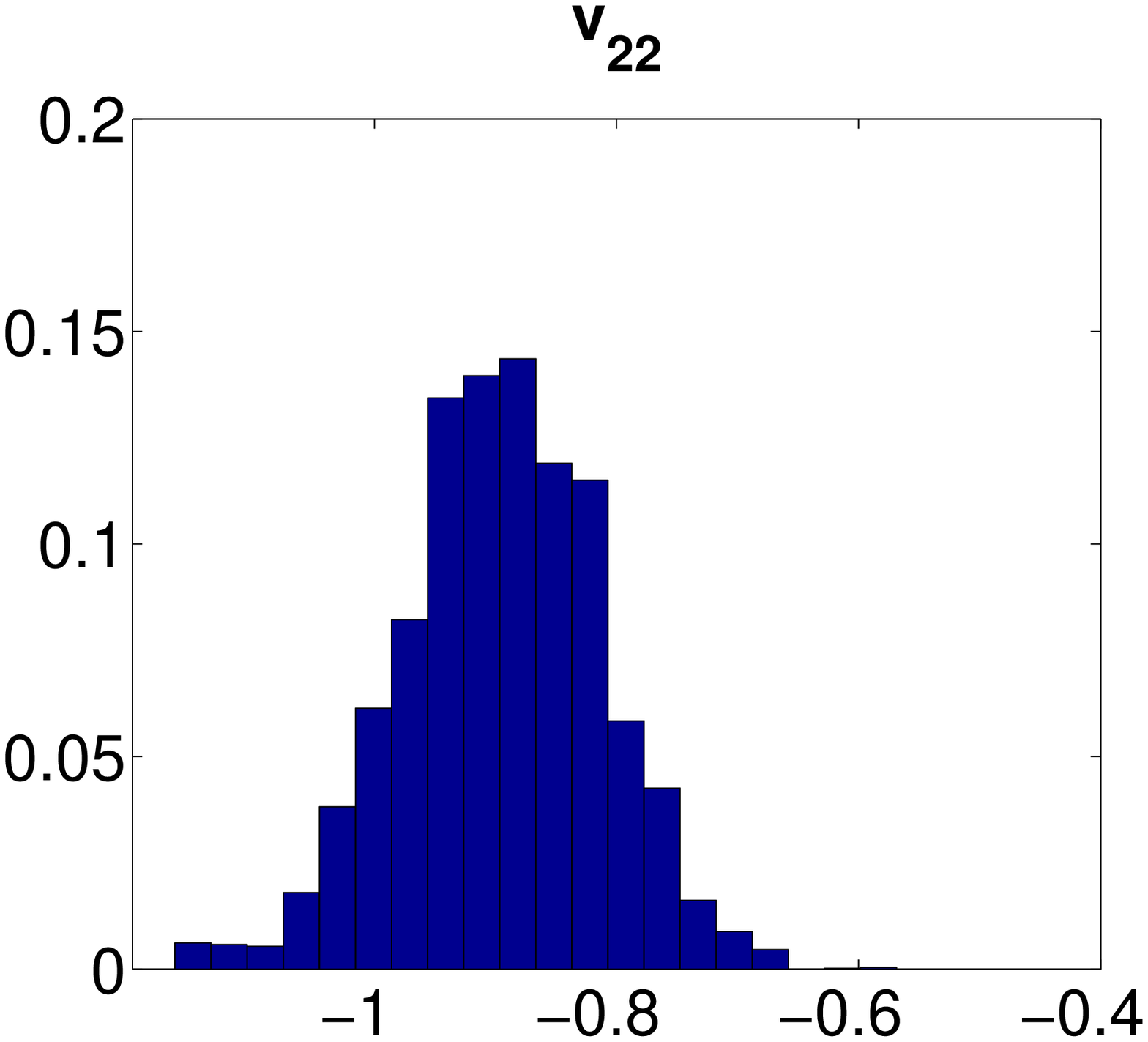} \\
\includegraphics[width=1\textwidth,height=1\textwidth]{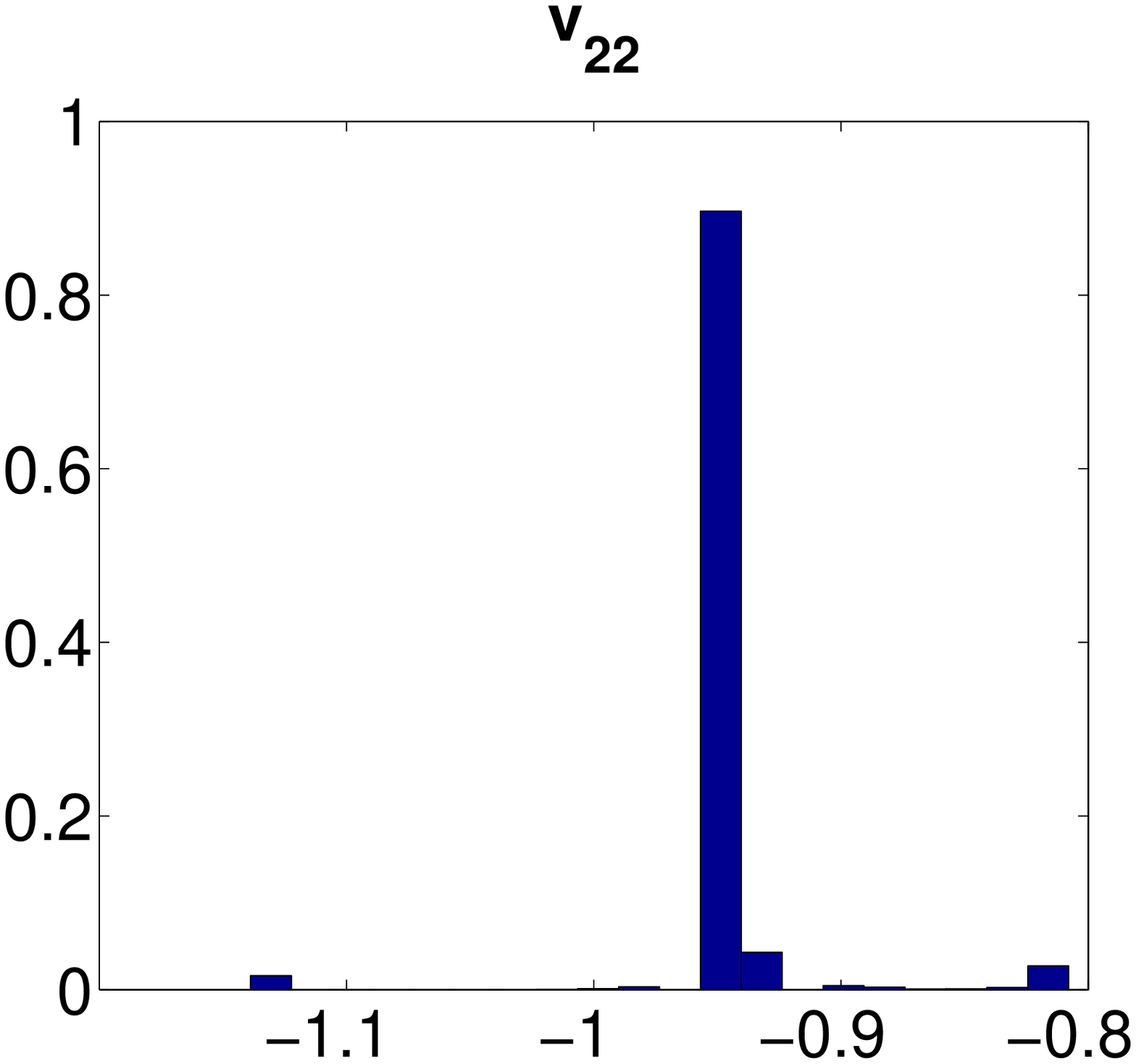}
\end{minipage}
}\subfigure[]{
\begin{minipage}[b]{0.17\textwidth}
\includegraphics[width=1\textwidth,height=1\textwidth]{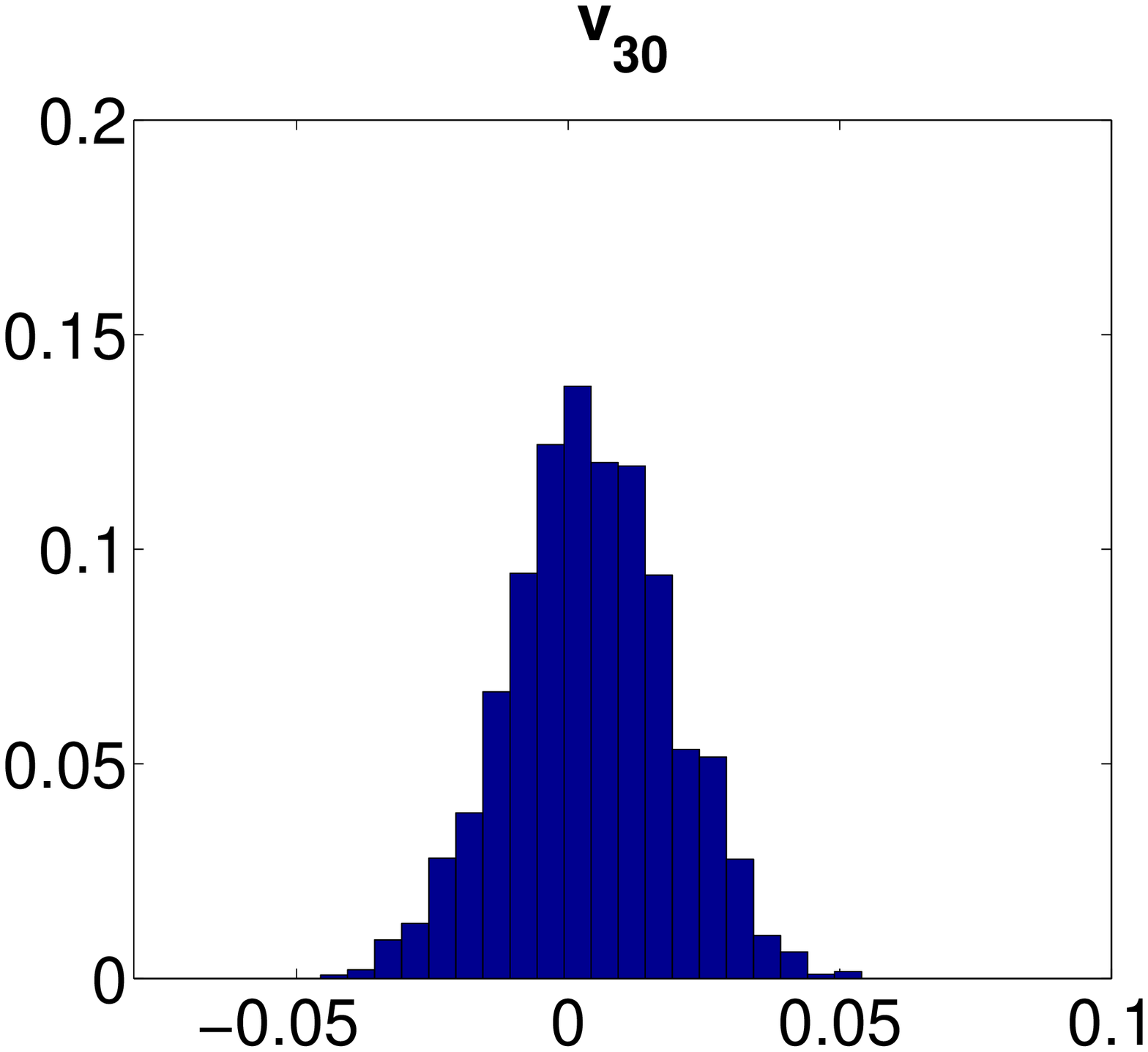} \\
\includegraphics[width=1\textwidth,height=1\textwidth]{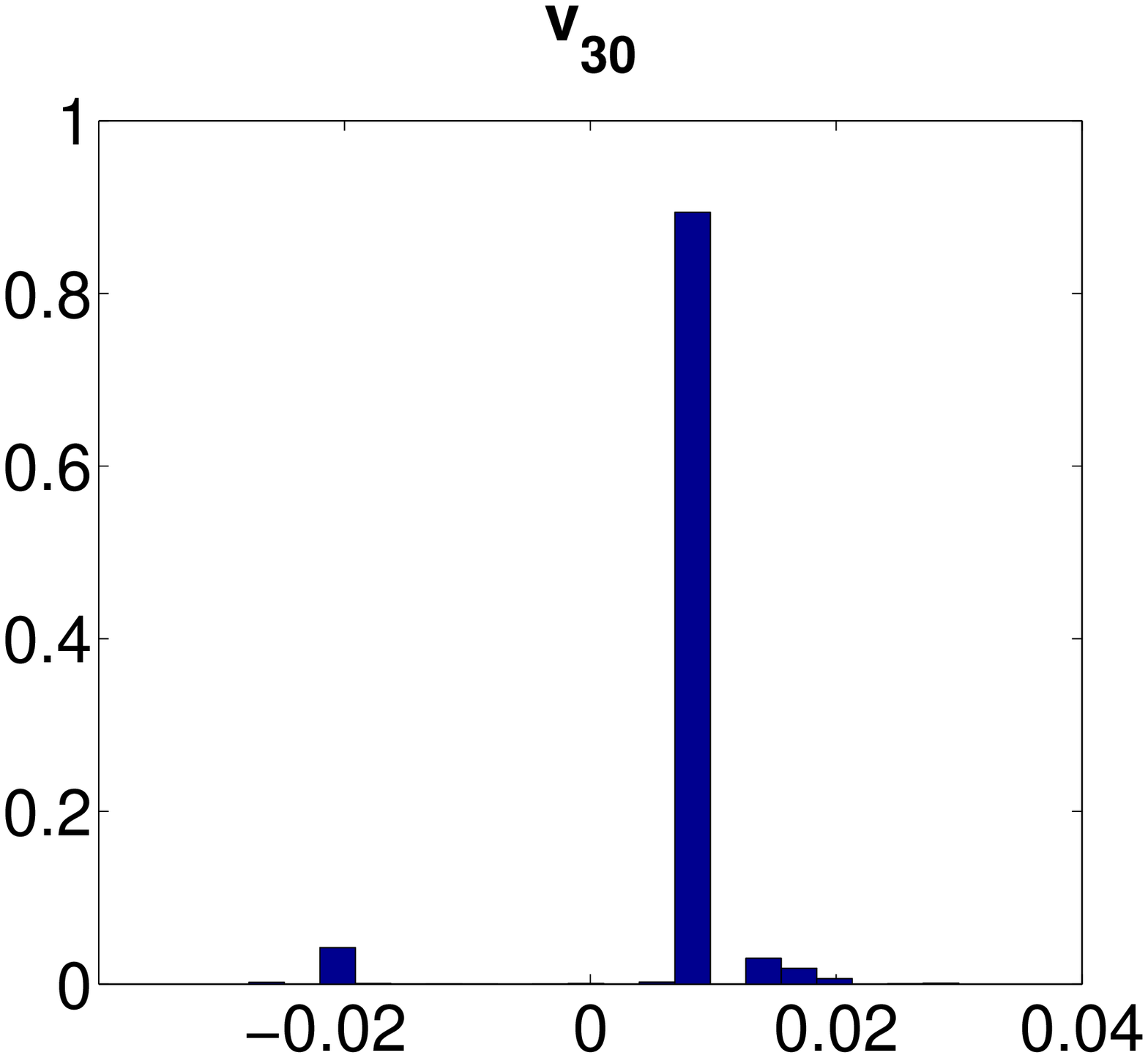}
\end{minipage}
}
\caption{The histogram of marginal posterior by improved  implicit sampling  (the first row) and conventional implicit sampling  (the second row) for $[v_{1}, v_{8}, v_{15}, v_{22}, v_{30}]$.}
\label{hist-3}
\end{figure}

\begin{figure}
\centering
\subfigure[]{
\includegraphics[width=0.35\textwidth, height=0.3\textwidth]{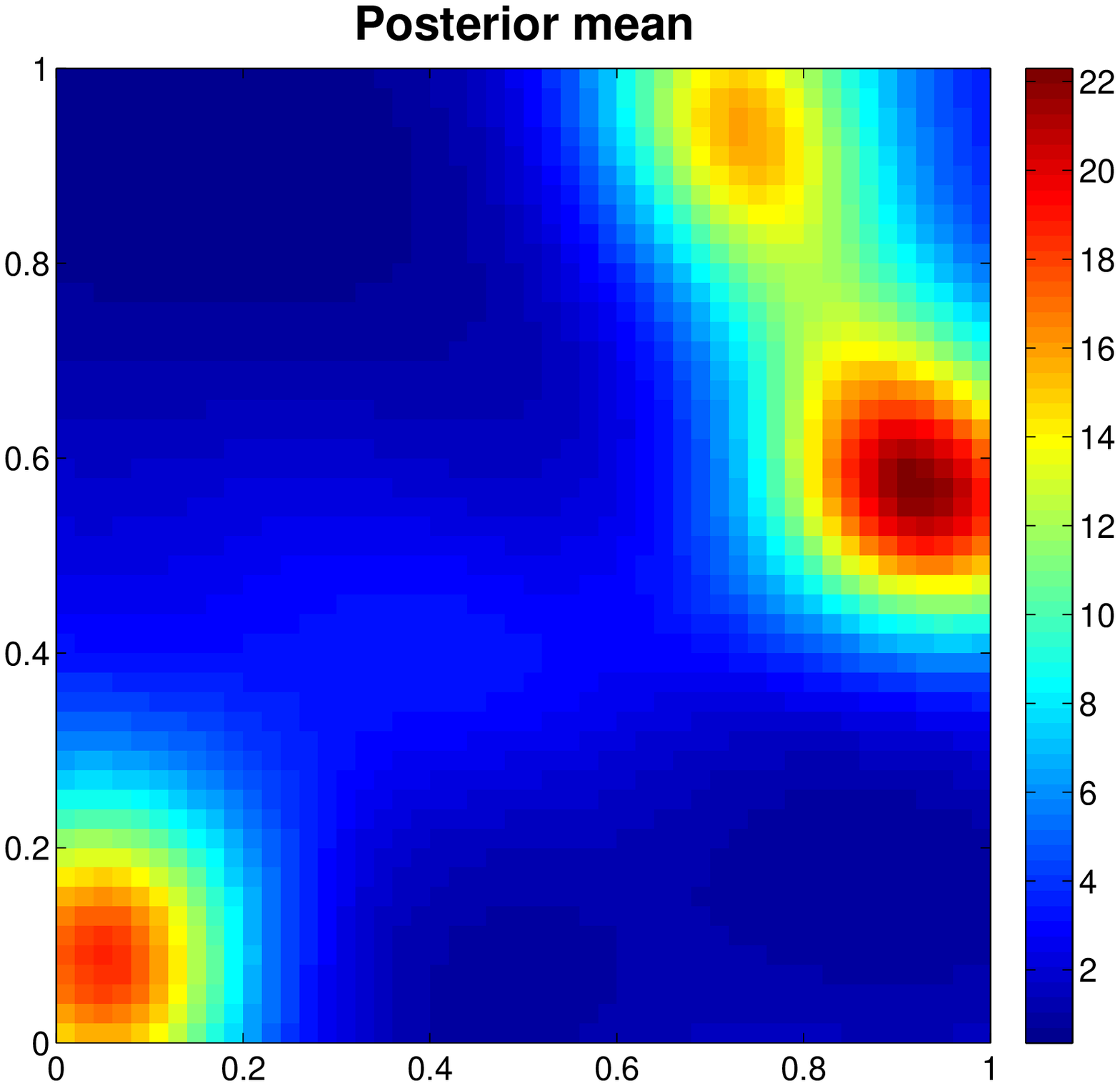}
}
\subfigure[]{
\includegraphics[width=0.35\textwidth, height=0.3\textwidth]{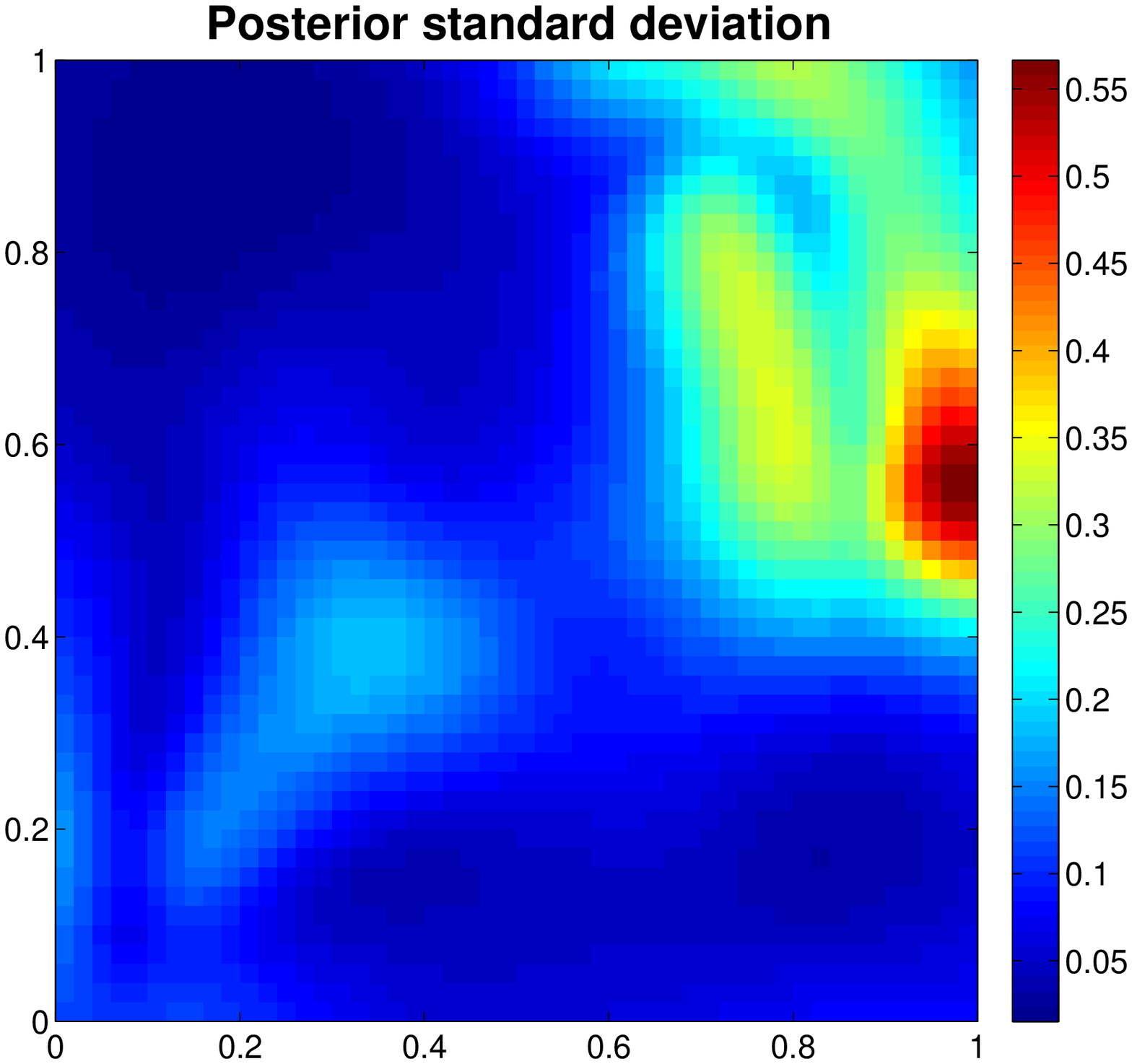}
}
\caption{The posterior mean (left) and posterior standard
deviation (right) of $q(x,\varrho)$ by improved implicit sampling. }\label{mean-std-ex3}
\end{figure}

In order to compute  the statistical information, $5000$  samples are produced by  implicit sampling. However, some of weights are  close to $1$ if we use the conventional implicit sampling.
To avoid this degenerative situation, we choose the scale parameter $\vartheta$ acting as the relaxation factor for improved implicit sampling.
 We can select the parameter $\vartheta$ dynamically based on the ESS in (\ref{ess_w}).
 In this example, we set $\vartheta=10$ to ensure the ESS to be  more than $1000$.
We list the distribution of weights corresponding  to samples by improved implicit sampling and conventional implicit sampling  in Table \ref{table-weights-3}.
It is obvious that the improved implicit sampling gives much better  distribution of weights of samples than the conventional implicit sampling.
The histogram of the marginal posterior by the two implicit sampling methods  for the parameter vector $[v_{1}, v_{8}, v_{15}, v_{22}, v_{30}]$ are plotted in Figure \ref{hist-3}, which agrees  with the results in Table \ref{table-weights-3}.
The posterior mean and standard derivation by the improved implicit sampling are plotted in Figure \ref{mean-std-ex3}. This figure  illustrates that the uncertainty mainly occurs around  the region with large  jump.

\begin{figure}[]
\centering
\subfigure[]{
\includegraphics[width=0.3\textwidth, height=0.3\textwidth]{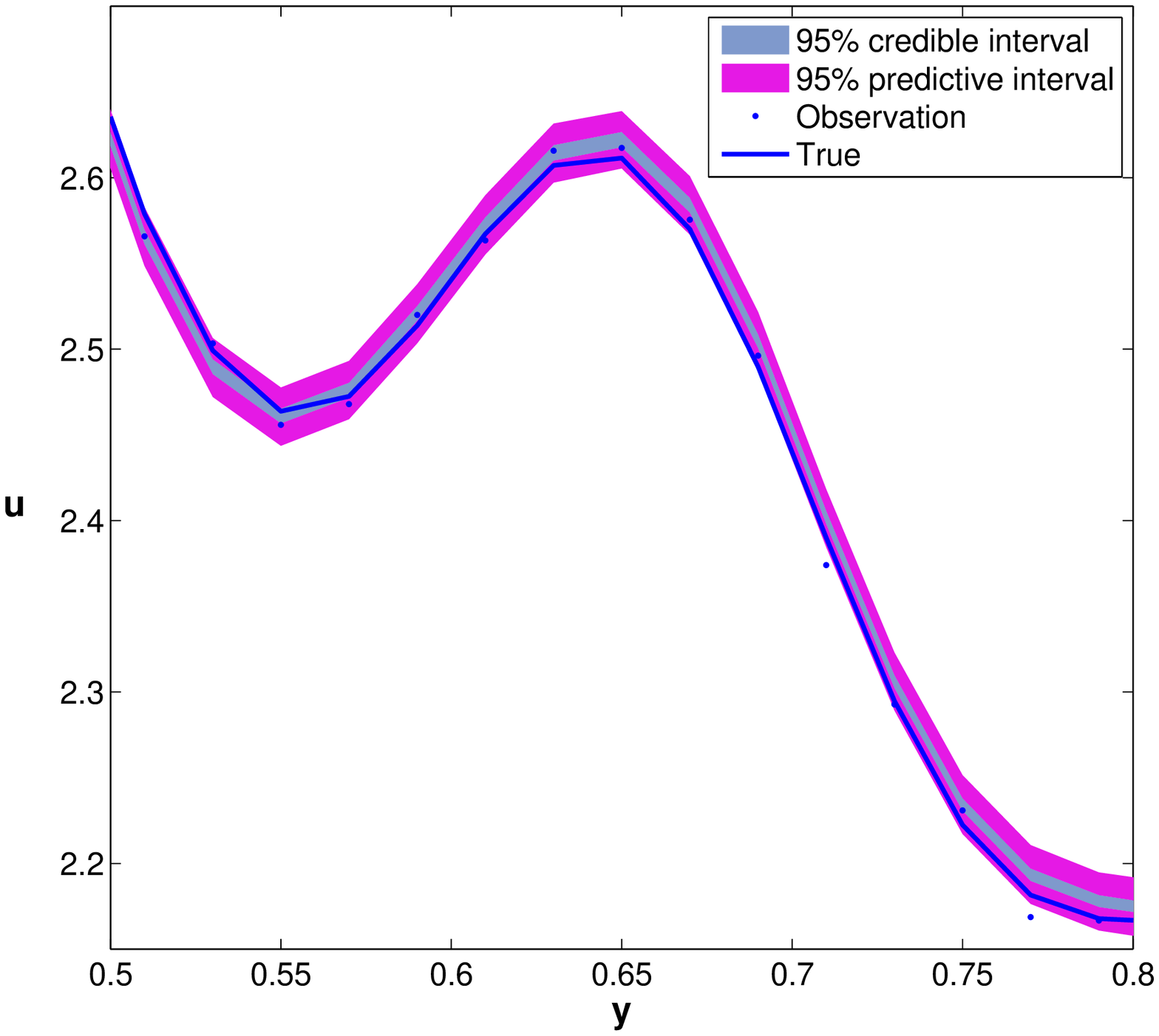}}
\quad
\subfigure[]{
\includegraphics[width=0.3\textwidth, height=0.3\textwidth]{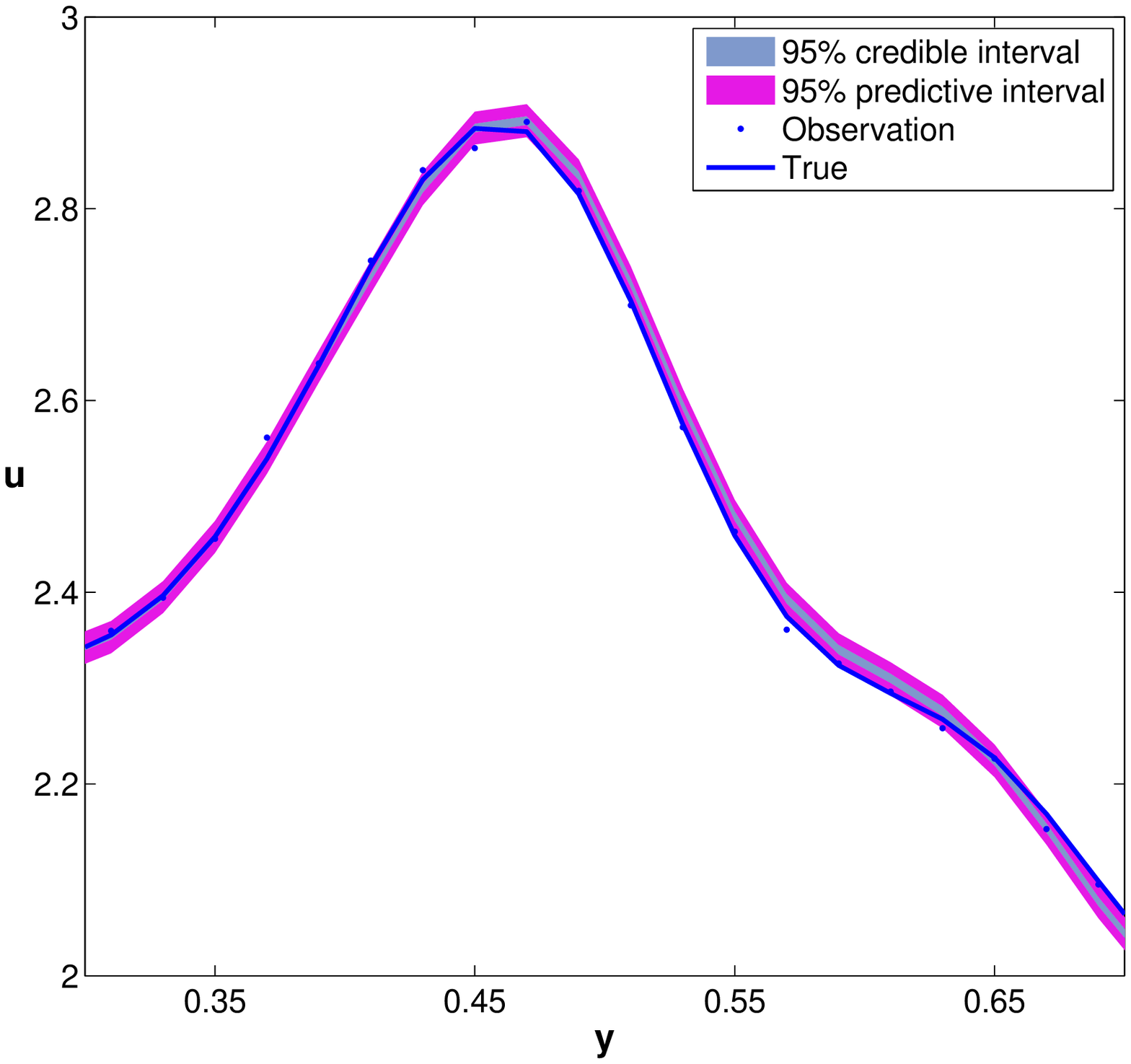}
}
\quad
\subfigure[]{
\includegraphics[width=0.3\textwidth, height=0.3\textwidth]{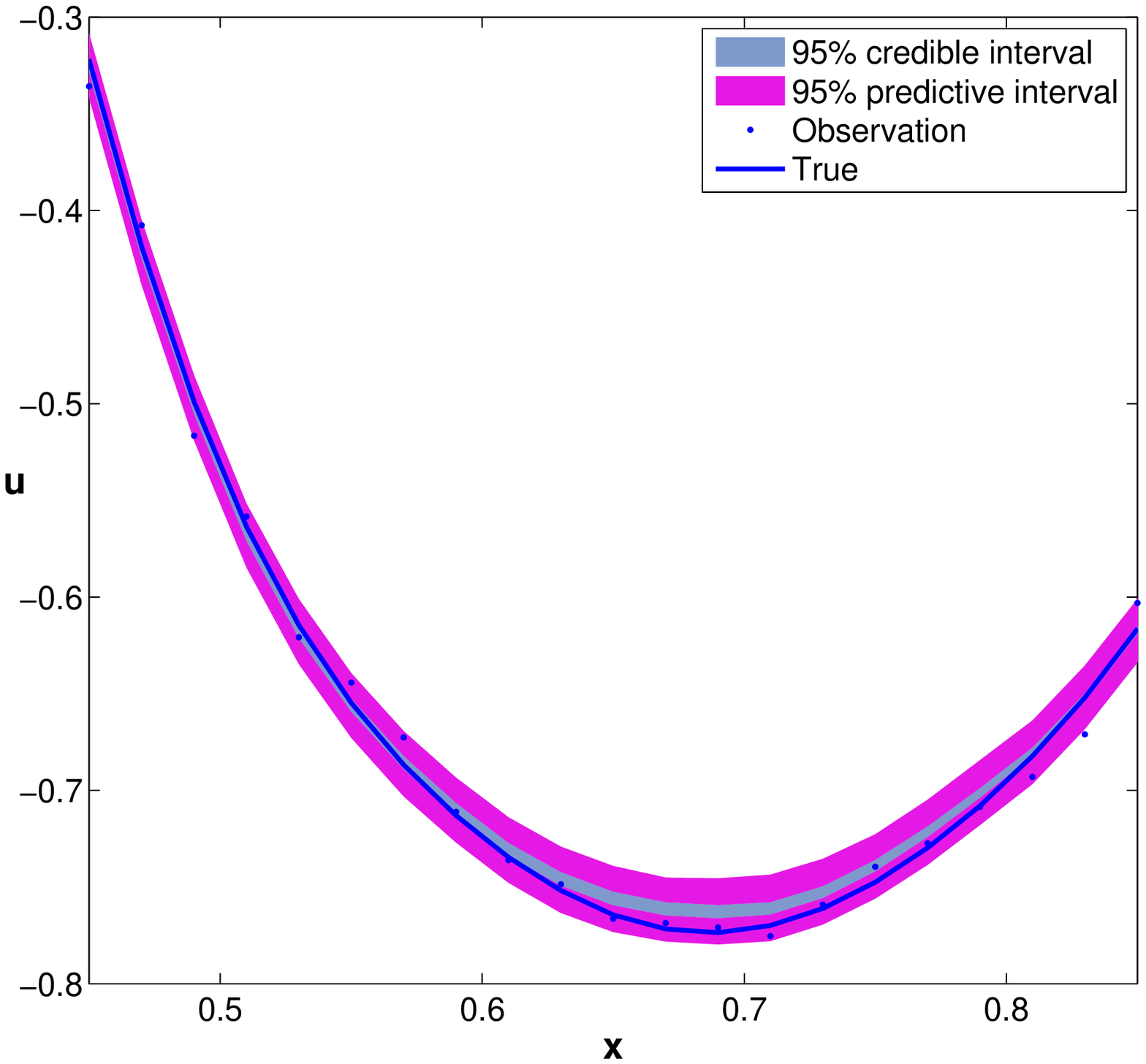}
}
\quad
\subfigure[]{
\includegraphics[width=0.3\textwidth, height=0.3\textwidth]{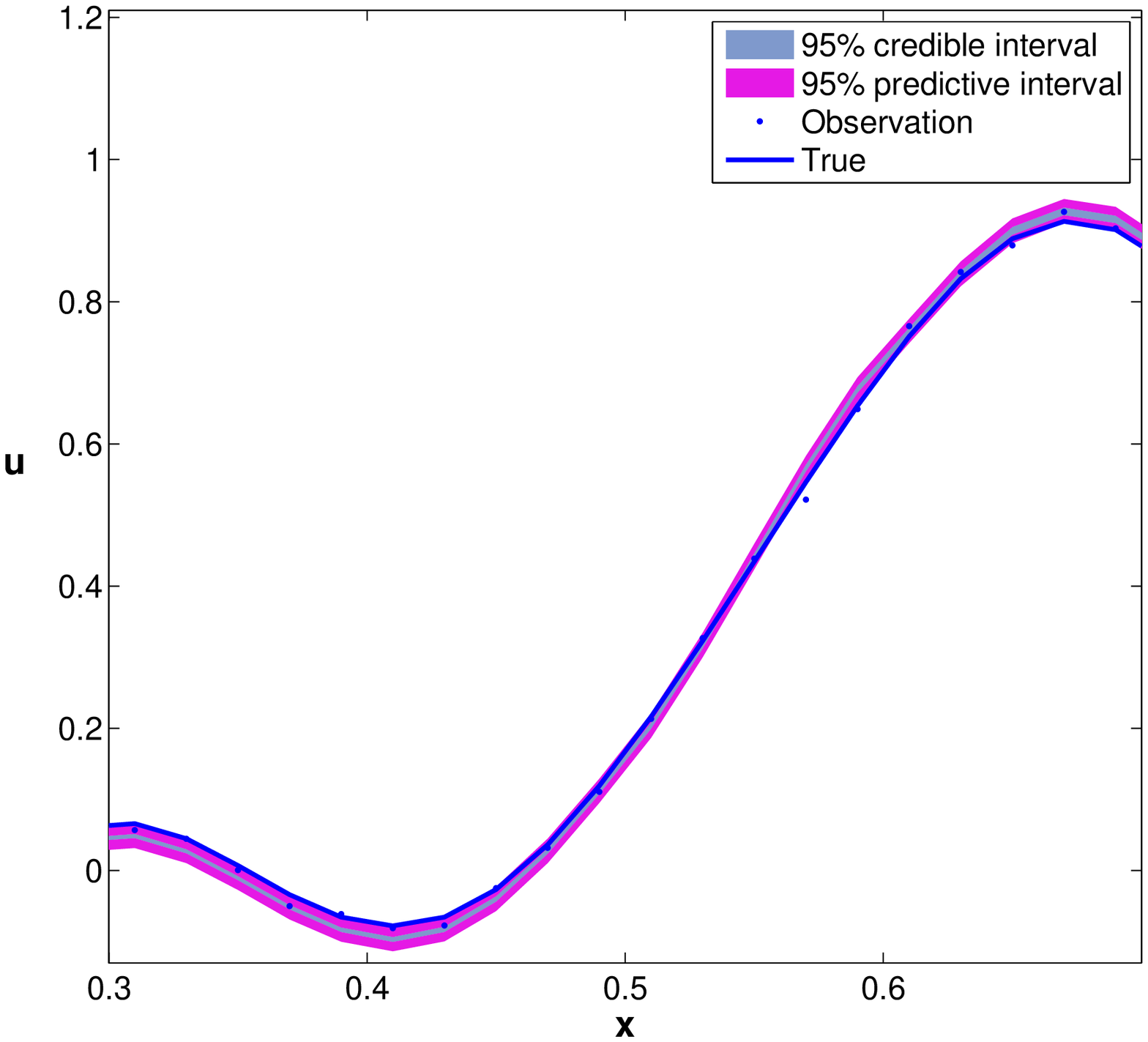}
}
\caption{Data, true, and $95\%$ credible interval and prediction interval  by improved implicit sampling for (a) $u(0,y; 0.6)$, (b) $u(1,y; 0.6)$, (c) $u(x,0; 0.6)$, (d) $u(x,1; 0.6)$.}
\label{cre_3}
\end{figure}

The $95\%$ credible intervals for model response at the measurement locations from the $5000$ realizations of the model are plotted in Figure \ref{cre_3}.
The prediction intervals are built by  incorporating the measurement error $\sigma$ and showed in Figure \ref{cre_3}, which shows  most measurement data lies  in the $95\%$ prediction interval.

\renewcommand\arraystretch{1.5}
\begin{table}
\centering
\resizebox{\textwidth}{!}{
   \begin{tabular}{cccccccc}
       \toprule
       \multirow{2}{*}{Method}
        & \multicolumn{7}{c}{Interval}\\
       \cmidrule{2-8}
       & [0, $10^{-6}$) & [$10^{-6}$, $10^{-5}$) & [$10^{-5}$,$10^{-4}$) & [$10^{-4}$,$10^{-3}$) & [$10^{-3}$,$10^{-2}$) & [$10^{-2}$,$10^{-1}$) & [$10^{-1}$,1]\\
       $\text{\uppercase\expandafter{I-IS}}$&2&186&2412&2262&138&0&0\\
       $\text{\uppercase\expandafter{C-IS}}$&4940&23&13&14&6&3&1\\
       \bottomrule
  \end{tabular}
  }
  \caption{The distribution  of the weights for improved implicit sampling (I-IS) and conventional implicit sampling (C-IS).}\label{table-weights-3}
\end{table}

\section{Conclusions}
In this work, we presented an improved implicit sampling method for  Bayesian inverse problems.
The approach  generated  independent samples around  the region with high probability, i.e. near the MAP point.
However, the weights of conventional  implicit sampling  may cause excessive concentration of samples and lead to ensemble collapse, which would result in the poor estimation of the posterior. In order to avoid this issue,
 we proposed  a new weight formulation for  implicit sampling.
 For more practical applications, we considered the improved implicit sampling for  hierarchical Bayesian models.
As an application, we focused on the multi-term time fractional multiscale diffusion equations.
To effectively  capture the multiscale and heterogeneity feature of the diffusion field,  we presented a mixed GMsFEM to build a reduced model  and  speed up the Bayesian inversion.
Some comparison were carried out among  improved implicit sampling,conventional sampling and LMAP method by a few numerical examples.
We found that the proposed improved implicit sampling is able to effectively improve samples quality and posterior inference than the conventional implicit sampling and LMAP.

\begin{ack}
L. Jiang acknowledges the support of Chinese NSF 11871378 and TZ2018001.
\end{ack}

\bibliographystyle{siam}
\bibliography{ref2}

 \end{document}